\documentclass[a4paper,10pt]{article}
\usepackage{a4}
\usepackage{amsmath, amsfonts, amssymb, amsthm, amscd, bezier}
\usepackage[brazil,english]{babel}
\usepackage[latin1]{inputenc}
\usepackage[T1]{fontenc}
\usepackage{graphicx}
\usepackage[all]{xy}
\usepackage{color}
\usepackage{fullpage}
\usepackage{cite}
\advance\oddsidemargin by -0.5cm \advance\textwidth by 1cm
\usepackage{mathrsfs}
\newcommand{\cqd}{\hfill $\square$\vspace{0.2cm}}
\newtheorem{obs}{Remark}[section]

\newtheorem{prop}[obs]{Proposition}
\newtheorem{teo}[obs]{Theorem}
\newtheorem{lem}[obs]{Lemma}
\newtheorem{cor}[obs]{Corollary}

\theoremstyle{definition}
\newtheorem{ex}{Example}[section]

\newtheorem{defin}[obs]{Definition}

\def\dem{\noindent{\bf Proof:}\hspace{.1in}}

\def\href{}

\DeclareMathOperator{\Ker}{Ker}
\DeclareMathOperator{\Ima}{Im}

\begin{document}

\begin{center}

{ \bf  \LARGE Cancellations for Circle-valued Morse Functions via Spectral Sequences}


\vspace{0.5cm}

\hspace{0.3cm} D.V.S. Lima\footnote{ Supported by FAPESP under grants 2014/11943-6 and 2015/10930-0.}  
\hspace{0.3cm} 
O. Manzoli  Neto \footnote{Supported by FAPESP under grant 2012/24249-5} 
\hspace{0.3cm}K.A. de Rezende\footnote{Partially supported by CNPq under grant 309734/2014-2 and by FAPESP under grant 2012/18780-0.}
\hspace{0.3cm}
 M. R. da
Silveira\footnote{Partially supported by FAPESP under grant 2012/18780-0.}    

\end{center}

\begin{abstract}
In this article, a spectral sequence analysis of a filtered Novikov complex $(\mathcal{N}_{\ast}(f),\Delta)$ over $\mathbb{Z}((t))$ is developed with the goal of obtaining results 
relating the algebraic and dynamical settings. Specifically, the unfolding of a spectral sequence of $(\mathcal{N}_{\ast}(f),\Delta)$ and the cancellation of its modules is associated to a one parameter family of circle valued Morse functions  on a surface and the dynamical cancellations of its critical points. The data of a spectral sequence computed for $(\mathcal{N}_{\ast}(f),\Delta)$ is encoded in a family of matrices $\Delta^r$ produced by the Spectral Sequence Sweeping Algorithm (SSSA), which has as its initial input the differential $\Delta$. As one ``turns the pages'' of the spectral sequence, differentials which are isomorphisms produce cancellation of pairs of modules. Corresponding to these cancellations, a family of circle-valued Morse functions $f^r$ is obtained by successively removing the corresponding pairs of critical points  of $f$. We also keep track of all dynamical information on the birth and death of connecting orbits between consecutive critical points,  as well as  periodic orbits that arise within a family of negative gradient flows associated to $f^r$.

\end{abstract}

{\small
}

\vspace{-.5cm}

\section*{Introduction }

Given $f:M\to\mathbb{R}$ a Morse
function, consider a filtered  Morse chain complex associate to it. Once the dynamics is presented at
the chain complex level, it is natural to use homology techniques to explore it. 
In \cite{BLMdRS,CdRS,FdRS,MdRS} we have chosen to work with a powerful homological tool, namely, a spectral sequence of a filtered Morse chain complex. Spectral sequence analysis has proven to be useful in relating  algebraic and  dynamical  information.
This relation has been explored for gradient flows generated by  Morse functions and their filtered  Morse chain complexes.

Given a Morse function $f: M \rightarrow \mathbb{R}$ on a smooth  closed $n$-manifold, the Morse chain complex of $f$, $(C_{\ast}(f),\partial)$, where $C_{\ast}(f)$ is a $\mathbb{Z}$-module generated by the critical points and graded by their indices, i.e., $$C_{k}(f) = \bigoplus_{x\in Crit_{k}(f)} \mathbb{Z} \langle x \rangle, $$
where $Crit_{k}(f)$ is the set of critical points of $f$ of index $k$. The differential $\partial$ is defined on the generators of $C_{k}$ by 
$$ \partial_{k} (x) = \sum_{y\in Crit_{k-1}(f)} n(x,y) \langle y \rangle, $$
where $n(x,y)$ is the intersection number of $x$ and $y$. See  \cite{Sa2} for more details. 

It was proven in \cite{CdRS, MdRS} that a spectral sequence $(E^{r},d^{r})$ of a filtered Morse chain complex $(C_{\ast}(f),\partial)$  can be retrieved from its differential $\partial$ via the Spectral Sequence Sweeping Algorithm - SSSA.
In  \cite{FdRS}, spectral sequences were considered over $\mathbb{Z}_{2}$ and it was proven that  changes of generators of the modules $E^{r}$ determine bifurcation behavior.  
In \cite{BLMdRS} and \cite{BLMdRS2},  the spectral sequence analysis was realized for an $n$-dimensional filtered Morse chain complex over $\mathbb{Z}$,  $n=2$ or $n\geq 5$. This analysis relates the algebraic cancellations caused by non-zero differentials $d^{r}$ to dynamical cancellations of pairs of  critical points. Consequently, as one ``turns''  the pages of the spectral sequence, dynamical cancellations occur successively. Thus, we obtain a global dynamical cancellation theorem which keeps track of  birth and death of connections of critical points.

Motivated by outlook, in this article, we undertake a completely new dynamical setup by considering circle-valued Morse functions $f:M \rightarrow S^{1}$ on a closed manifold $M$. Similar to the Morse complex, the Novikov complex $(\mathcal{N}_{\ast}(f),\partial)$ associated to $f$ is
a chain complex of free modules generated by the critical points of $f$, see Section \ref{novikovcomplex}.
However, instead of the  ring of integers, the Novikov complex is defined over $\mathbb{Z}((t))$, the ring  of Laurent series in one variable with integral coefficients and finite negative part. The Novikov boundary operator $\partial$ in this complex is related to the topological-dynamical features of a gradient flow of $f$. More specifically, $\partial$  ``counts''  orbits with signs connecting consecutive critical points over $\mathbb{Z}((t))$ and it is represented by a matrix denominated as a Novikov matrix $\Delta$.

In this article, we establish the groundwork for the association of a spectral sequence analysis of a filtered two dimensional Novikov chain complex with the dynamical exploration of connections in the negative gradient flow of $f$. This is done by applying the SSSA to $(\mathcal{N}_{\ast}(f),\partial)$.
More specifically, as one ``turns the pages'' of this spectral sequence, i.e.  considers progressively modules $E^r$, one  observes that each algebraic cancellation occurring among the $E^r$, tracked by the SSSA, corresponds to a dynamical cancellation of a pair of critical points of $f$. Moreover, as $r$ increases, the spectral sequence $(E^r,d^r)$ stabilizes and converges to the Novikov homology $H_{\ast}^{Nov}(M,f)$ of $f$. Although cancellation of critical points in the  Novikov setting has been done, e.g. in \cite{L}, the novelty in this article is a global dynamical cancellation theorem for circle-valued Morse functions via spectral sequence analysis. 

It is unquestionable that we have developed a phenomenal algebraic tool, SSSA, to compute the modules and differentials $d^r$ of a spectral sequence of a filtered Novikov chain complex. Although this result  has considerable algebraic value in itself, this effort is fully compensated by the algebraic-dynamical cancellation results herein.  A significant result uses the differentials $d^r$'s of a spectral sequence, that are responsible for algebraic cancellations, 
to understand critical point cancellations within a family of circle-valued Morse functions $f^r$ and the respective family of associated negative gradient flows $\varphi^r$. The dynamical significance of these algebraic cancellations, 
ultimately, provides a backstory of death and birth of orbits of $\varphi^r$ due to cancellations of consecutive critical points of $f^r$.
Also, as the spectral sequences analysis unfolds, it is possible to detect all  periodic orbits that arise within the family $\varphi^r$.

In Section 2, we recall the definition of a spectral sequence for a filtered chain complex and present the filtrations  considered throughout this paper. The main tool behind this correspondence  is  the Spectral Sequence Sweeping Algorithm (SSSA), which is presented in  Section  \ref{sweeping}. In Section \ref{caracterizacao}, Theorem  \ref{teo:primarypivots} addresses the fact that the SSSA can be applied  to a Novikov  complex $(\mathcal{N}_{\ast}(f),\Delta)$ in dimension $2$.  Also, within this section, it is shown that the SSSA  produces a sequence of matrices, which we call Novikov matrices, that are characterized in Theorems \ref{bloco1} and \ref{bloco2}. In Theorem \ref{lemalinhanula}, we prove  the surprising result that the last matrix produced by the SSSA has polynomial entries{\footnote{A polinomial  in $\mathbb{Z}((t))$ is a finite sum of powers of $t$, $t^{\ell}$, with $\ell \in \mathbb{Z}$.}} in $\mathbb{Z}((t))$, although the intermediate Novikov matrices may exhibit entries which are infinite series. In Section \ref{ultimo},  we prove in Theorems \ref{teo:inducaodoe} and \ref{teo:interpretacaodr}  that the modules and differentials $(E^r, d^r)$ of the spectral sequence may be retrieved from  the sequence of Novikov matrices produced by the SSSA. More specifically, the SSSA provides a system which spans $E^r$ in terms of the original basis of $\mathcal{N}_{\ast}(f)$ as well as  identifies all differentials $d^r_p:E^r_p\to E^r_{p-r}$.  Finding a system of $(r-1)$-cycles that span $E^{r}$ in terms of the original basis of $\mathcal{N}_{\ast}(f)$ is a non-trivial matter. In Section \ref{sec6} is to obtain in Theorem 6.1 a global dynamical cancellation theorem for circle-valued Morse functions  via spectral sequences analysis.

Hence, the emphasis in this article is on the spectral sequence analysis and its correlation to the dynamics. The importance of this method is that it establishes the foundation for future investigation on the algebraic implications to the dynamics in the Novikov setting in higher dimensions.

\section{Characterization of the Novikov Complex}\label{novikovcomplex}

The goal of this section  is to prove  a characterization of the Novikov differential in the case of orientable surfaces, which will be essential for both dynamical and algebraic further results. Initially we present some  background material  on circle-valued functions and Novikov complex over $\mathbb{Z}((t))$. Further details can be found in \cite{P}.  

Denote by $\mathbb{Z}[t,t^{-1}]$ the Laurent polynomial ring. Let $\mathbb{Z}((t))$ be the set consisting of all  Laurent series 
$$ \lambda= \sum_{i\in \mathbb{Z}}a_{i} t^{i}  $$
in one variable with coefficients $a_{i} \in \mathbb{Z}$, such that the part  of $\lambda$ with negative exponents  is finite, i.e., there is $n=n(\lambda)$ such that $a_k=0$ if $k<n(\lambda)$. In fact, $\mathbb{Z}((t))$ has a natural ring structure such that the inclusion $\mathbb{Z}[t,t^{-1}] \subset \mathbb{Z}((t))$ is a homomorphism. Moreover, $\mathbb{Z}((t))$ is a Euclidean ring. 

Let $M$ be a closed connected manifold  and $f:M\rightarrow S^{1}$ be a smooth map from $M$ to the one-dimensional sphere{\footnote{ $S^{1}$ is viewed as the submanifold $\{(x,y) \in \mathbb{R}^{2} \mid x^{2}+y^{2}=1\}$ and is endowed with the corresponding smooth structure.}}. Given a point $x\in M$ and a neighbourhood $V$ of $f(x)$ in $S^{1}$ diffeomorphic to an open interval of $\mathbb{R}$, the map  $f|_{f^{-1}(V)}$ is identified to a smooth map from $f^{-1}(V)$ to $\mathbb{R}$. Therefore, in this context  one can define non-degenerate critical points and Morse index as in the classical case of  smooth real valued function. A smooth map $f:M\rightarrow S^{1}$ is called a {\it circle-valued Morse function} if its critical points are non-degenerate. The set of critical points of $f$ will be denoted by $Crit(f)$ and $Crit_{k}(f)$ is the set of critical points of $f$ of index $k$.

 Considering the exponential function $Exp: \mathbb{R} \rightarrow S^{1}$ given by $t\mapsto \mathrm{e}^{2\pi i t}$ and the infinite cyclic covering $E:\overline{M}\rightarrow M$ induced by the map $f:M\rightarrow S^{1}$ from the universal covering $Exp$,
  there exists a map $F:\overline{M}\rightarrow \mathbb{R}$ which makes the following diagram commutative:
  \vspace{-0.3cm}
 \begin{equation}\nonumber
   \xymatrix@R+2em@C+2em{
  \overline{M} \ar[r]^-F \ar[d]_E & \mathbb{R} \ar[d]^{Exp} \\
  M \ar[r]^-f  & \mathbb{S}^1
  }
 \end{equation}
Moreover, $f$ is a { circle-valued Morse function} if and only if $F$ is a real valued  Morse function.  Observe that if $Crit(F)$ is non empty then it  has  infinite cardinality. If $\overline{M}$ is non compact, one can not apply the classical Morse theory to study $F$, however, one can restrict $F$  to a fundamental cobordism $W$ of $\overline{M}$, which is compact, and apply the techniques of Morse theory. The {\it fundamental cobordism} $W$ is defined as 
$ W = F^{-1}([a-1,a]),  $
 where $a$ is a regular value of $F$. The cobordism $W$ can be viewed  as the manifold $M$ obtained by cutting along the submanifold $V=f^{-1}(\alpha)$, where $\alpha = Exp(a)$. Hence, $W$ is a cobordism with both boundary components diffeomorphic to $V$.

Given a circle-valued Morse function $f$, consider the vector field $v=-\nabla f$. One says that $f$ satisfies the transversality condition  if the lift of $v$ to $\overline{M}$ satisfies the classical transversality condition on the unstable and stable manifolds.

From now on, we consider circle-valued Morse functions $f$ such that $v=-\nabla f$ satisfies the transversality condition. Denote by $\overline{v}$ the lift of $v$ to $\overline{M}$ and arbitrarily choose orientations for all unstable manifolds $W^{u}(p)$ of critical points of $f$.

Given $p \in Crit_{k}(f)$ and $q \in Crit_{k-1}(f)$, {\it the Novikov incidence coefficient} between $p$ and $q$ is defined as 
$$ N(p,q;f) = \sum_{\ell \in \mathbb{Z}} n(p,t^{\ell}q;\overline{v}) t^{\ell},  $$
where $n(p,t^{\ell}q;\overline{v})$ is the intersection number between the critical points $p$ and $t^{\ell}q$ of $F$, i.e., the number obtained by counting with signs the flow lines of $\overline{v}$ from $p$ to $t^{l}q$, when one considers the  orientations on the unstable manifolds $W^{u}(p)$ and $W^{u}(t^{\ell}q)$ according to the previous fixed orientations in $W^{u}(p)$ and $W^{u}(q)$. For more details, see \cite{BH1} and \cite{P}. 

Let $\mathcal{N}_{k}$ be  the $\mathbb{Z}((t))$-module freely generated by the critical points of $f$ of index $k$. Consider the $k$-th boundary operator $\partial_{k}: \mathcal{N}_{k} \rightarrow \mathcal{N}_{k-1}$ which is  defined on a generator $p \in Crit_{k}(f)$ by 
$$ \partial_{k}(p) = \sum_{q\in Crit_{k-1}(f)} N(p,q;f)q $$
and extended to all chains by linearity. In \cite{P} it is proved that $\partial_{k}\circ\partial_{k+1}=0$, hence $(\mathcal{N}_{\ast}(f),\partial)$ is a chain complex which is called the {\it Novikov complex} associated to the pair $(f,v)$.

\begin{figure}[!ht]
\centering
\includegraphics[width=6.7cm]{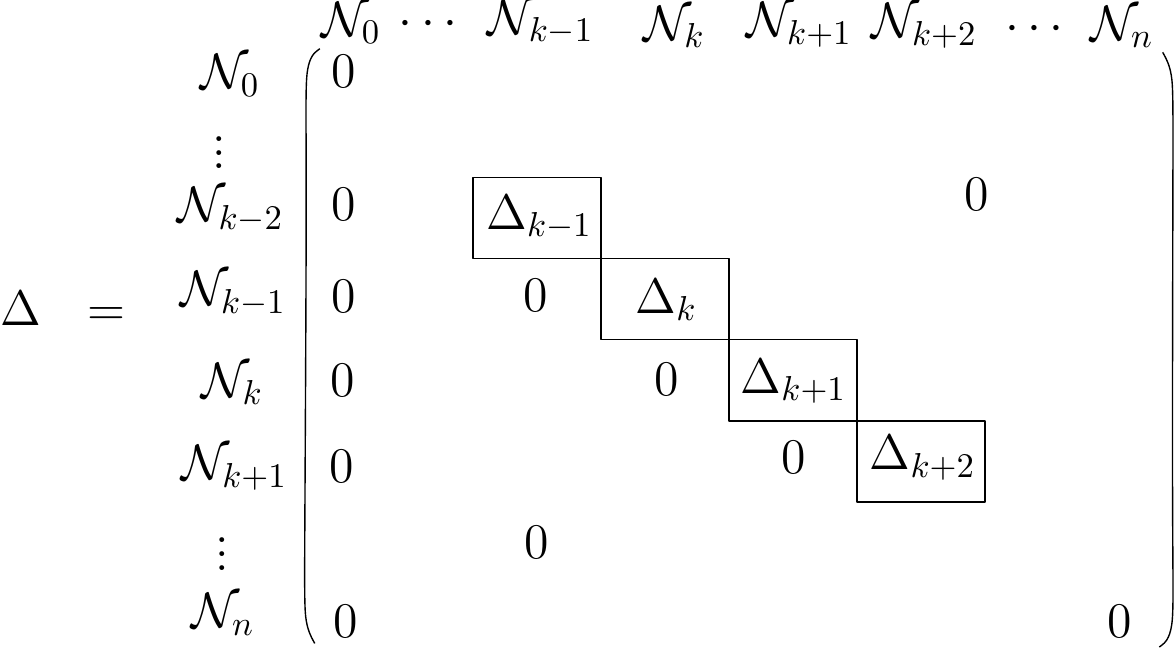} \\ \vspace{-0.5cm}
\caption{Novikov differential viewed as a matrix, where $\Delta_{k}$ is the matricial representation of $\partial_{k}$.}\label{fig:differential}
\end{figure}

One can view the Novikov differential $\partial$ as a matrix $\Delta$ where each column corresponds to generators  $p,q \in Crit(f)$  and the entries are the coefficients $N(p,q;f)$ of the Novikov differential $\partial$. Moreover, without loss of generality,  one  assumes  that the columns of $\Delta$ are ordered with respect to the Morse indices of the critical points, e.g., in increasing order of Morse index. See Figure \ref{fig:differential}.

\begin{figure}[!ht]
\centering
\includegraphics[width=5cm]{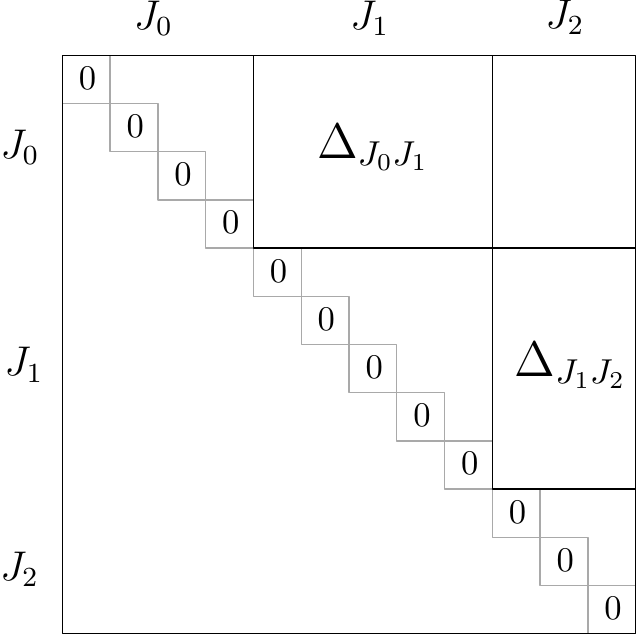} \\ \vspace{-0.3cm}
\caption{ Novikov differential, with $J_{0}=\{1,2,3,4\}$, $J_{1}=\{5,6,7,8,9\}$ and $J_{0}=\{10,11,12\}$.}\label{fig:cmsurface}
\end{figure}

In the case when $M$ is a surface, the columns of $\Delta$ may be partitioned into subsets $J_{0}, J_{1}, J_{2}$ such that  $J_{s}$ are the columns  associated with critical points of index $s$, i.e., the generators of $\mathcal{N}_{s}$. Hence, the non-zero entries of the matrix $\Delta$ are in the block $J_{0}\times J_{1}$, which corresponds to connections from saddles to sinks, and in the block $J_{1}\times J_{2}$, which corresponds to  connections from sources to saddles. The block $J_{0}\times J_{1}$ (respectively, $J_{1}\times J_{2}$) is referred to  as the first block (respectively, second block) of the matrix $\Delta$.  Figure \ref{fig:cmsurface} illustrates a possible structure for a Novikov differential associated to a circle-valued Morse function on a surface.

From this point on, 
 we will use the notation $\Delta$  to denote interchangeably the Novikov differential $\partial$ and its matrix representation.

 The following theorem presents   special characteristics of the Novikov differential  $\Delta$ by describing the Novikov incidence coefficients. In order to do this, we must define a fundamental domain. A cobordism $W=F^{-1}([a-1;a+\lambda])$, where $a$ is regular value of $F$ and $\lambda \in \mathbb{N}$, is said to be a {\it fundamental
domain} for $(M;f)$ if the following property is satisfied: given $p\in Crit_{k}(f)$  and $q\in Crit_{k-1}(f)$, $W$ contains a
lift of each orbit of the flow $v$ from $p$ to $q$.

\begin{teo}\label{pro:initialmatrix} 
 Let $M$ be an orientable surface and $(\mathcal{N}_{\ast}(f),\Delta)$ be the Novikov complex associated to a circle-valued Morse function $f: M \rightarrow S^{1}$. The Novikov incidence coefficient $N(p,q;f)$ is  either zero, a monomial $\pm t^{\ell}$ or a binomial $ t^{\ell_{1}} -t^{\ell_{2}} $.
\end{teo}

\dem
Given an orientable surface $M$, let $W_{\lambda} = F^{-1}([a-1,a+\lambda])$ be a fundamental domain for $(M,f)$.  Then $W_{\lambda}$ is an orientable compact surface with boundary $\partial W_{\lambda}$, possibly empty. If $\partial W_{\lambda} =\emptyset $, define $\widetilde{W}_{\lambda} = W_{\lambda}$. In the case that the boundary $\partial W_{\lambda}$ is non empty, it is the disjoint union of $\partial W_{\lambda}^{-} = F^{-1}(a-1)$ and $\partial W_{\lambda}^{+} = F^{-1}(a+\lambda)$.
Let $\widetilde{W}_{\lambda}$ be the closed surface obtained from $W_{\lambda}$ by gluing a  $2$-dimensional disk along its boundary to each connected component of $\partial W_{\lambda}$.
 Moreover, one can assume that each of its disks  contains a singularity, more specifically,  a source if the disk is glued to $\partial W_{\lambda}^{+}$ and a sink if the disk is glued to  $\partial W_{\lambda}^{-}$.  This procedure extends the Morse function $F:W_{\lambda} \rightarrow \mathbb{R}$ to a classical Morse function $\widetilde{F}:\widetilde{W}_{\lambda} \rightarrow \mathbb{R}$ on a closed surface.

Given $p\in Crit_{k}(f)$ and $q\in Crit_{k-1}(f)$, the Novikov incidence coefficient $N(p,q;f)$ counts the number of flow lines from $p$ to $q$ with signs. Since, each of these flow lines has a lift in ${W}_{\lambda}$, then $N(p,q;f)$ can be obtained by analysing $W_{\lambda}$. On the other hand, the intersection number $n(p,t^{\ell}q;\overline{v})$ between critical points $p$ and $t^{\ell}q$ in $W_{\lambda}$ is the same as when  considered  in the surface $\widetilde{W}_{\lambda}$. Since $\widetilde{W}_{\lambda}$ is closed,  $n(p,t^{\ell}q;\overline{v})$ is zero, when there are two flow lines from $p$ to $t^{\ell}q$. It is $-1$ or $+1$, when there is  one flow line from $p$ to $t^{\ell}q$. See \cite{BLMdRS}.
 Therefore, the Novikov incidence coefficient $N(p,q;f)$ is:
 \begin{enumerate}
 \item[(a)] $0$, if there are two flow lines from $p$ to $q$ in $v$ which intersect the level set $f^{-1}(a)$ the same number of times;
 \item[(b)]  $\pm t^{\ell}$, if there is only one flow line from $p$ to $q$ which intersects $\ell$ times the level set $f^{-1}(a)$; 
 \item[(c)] $ t^{\ell_{1}} - t^{\ell_{2}}$, if there are two flow lines from $p$ to $q$ in $v$, one intersecting  $\ell_{1}$ times the level set $f^{-1}(a)$ and the other  one intersecting $\ell_{2}$ times this level set.\cqd
 \end{enumerate}

\begin{cor}[Characterization of the Novikov differential on orientable surfaces]\label{cor:charac}
Let $M$ be an orientable surface and $(\mathcal{N}_{\ast}(f),\Delta)$ be the Novikov complex associated to a circle-valued Morse function $f: M \rightarrow S^{1}$  such that the chosen orientation on the unstable manifold of each critical point of index $2$  is the same. Then there are three possibilities for either a column or a row $j \in J_{1}$ of $\Delta$:
\begin{enumerate}
\item[(1)]  all entries are null. 
\item[(2)] exactly one non zero entry which is a binomial $t^{\ell_1} - t^{{\ell_{2}}}$, for some $\ell_{1},\ell_{2} \in \mathbb{Z}$.
\item[(3)] exactly two non zero entries which are monomials  $t^{\ell_{1}}$ and $-t^{\ell_{2}}$, for some $\ell_{1},\ell_{2} \in \mathbb{Z}$.
\item[(4)] exactly one non zero entry which is a monomial $t^{\ell}$, for some $\ell \in \mathbb{Z}$.
\end{enumerate}
\end{cor}
\dem
 If $f$ does not have a critical point of index $1$, then the Novikov matrix is null. Suppose that $f$ has at least one critical point of index $1$, i.e., a saddle. In this case, it is clear that, given a row (respectively, column) $j\in J_{1}$ of $\Delta$, there are at most two non zero entries in this row (respectively, column). In fact, there are exactly two flow lines whose $\omega$-limit (respectively, $\alpha$-limite) sets are the same  saddle. By \cite{BLMdRS}, choosing the same orientation for each unstable manifold of  critical points of index $2$, the signs on flow lines associated to the stable (respectively, unstable) manifold  of a saddle are opposite. 
Therefore, if there are two flow lines from a source  (respectively, saddle) $p$ to a saddle (respectively, sink) $q$ intersecting a regular level set $f^{-1}(a)$ $\ell$ times, then $N(p,q;f) = t^{\ell} -t^{\ell }=0$. On the other hand, if there is a flow line from a source (respectively, saddle) $p$ to a saddle (respectively, sink) $q$ intersecting $\ell_{1}$ times a regular level set $f^{-1}(a)$ and a flow line from a source (respectively, saddle) $p'$ to the saddle (respectively, sink) $q$ intersecting $\ell_{2}$ times the same regular level set, then $N(p,q;f) = \pm t^{\ell_{1}} $ and $N(p',q;f)=\mp t^{\ell_{2} }$. In particular, if $p'=p$ then $N(p,q;f) = \pm t^{\ell_{1}} \mp t^{\ell_{2} }$. Case $(4)$ occurs in two situations. First  when, given a saddle $q$, one of the flow lines in $W^{s}(q)$ has as $\alpha$-limit  a source $p$ and the other flow line in $W^{s}(q)$ has as $\alpha$-limit  a periodic orbit. 
Hence $N(p,q;f) = \pm t^{\ell} $. Secondly, given a saddle $p$,  one of the flow lines in $W^{u}(p)$ has as $\omega$-limit  a sink $q$ and the other flow line in $W^{u}(p)$ has as $\omega$-limit  a  periodic orbit.
Hence $N(p,q;f) = \pm t^{\ell}$. 
\cqd

\section{Spectral Sequence for a Novikov Complex}\label{spectralsequence}

In this section, we present the basic results on spectral sequences associated to filtered chain complexes. This section is based on the references \cite{D} and \cite{Sp}.

A {\it spectral sequence} $E$ over $\mathbb{Z}((t))$ is a sequence $\{E^{r},\partial^{ r}\}$, for $r \geq 0 $, such that:  $E^{r}$ is a bigraded module over $\mathbb{Z}((t))$, i.e., an indexed  collection of $R$-modules $E^{r}_{p,q}$, for all $p,q\in\mathbb{Z}$;
 $d^{r}$ is a differential of degree $(-r,r-1)$ on $E^{r}$, i.e., an indexed collection of homomorphisms  $d^{r}:E^{r}_{p,q}\rightarrow E^{r}_{p-r,q+r-1}$, for all $p,q\in\mathbb{Z}$, and $(d^{r})^{2}=0$. Moreover, 
 there exists an isomorphism 
$H(E^r)\approx E^{r+1}$, for all $r\geq 0$, where $$H_{p,q}(E^r)=\frac{\mbox{Ker}
d^r:E^{r}_{p,q}\to E^r_{p-r,q+r-1}}{\mbox{Im} d^r:E^{r}_{p+r,
q-r+1}\to E^r_{p,q}}. $$

Given a chain complex $(C_{\ast},\partial)$, an \emph{increasing filtration} $\mathcal{F}$ in $C_{\ast}$ is a sequence of submodules $\mathcal{F}_pC$ for all integers $p$ such that $\mathcal{F}_pC\subset \mathcal{F}_{p+1}C$. Since $C$ is a graded module, the filtration must be compatible with the gradation, i.e. $\mathcal{F}_pC$ is graded by
$\{\mathcal{F}_pC_q\}$.

Let $(\mathcal{N}_{\ast}(f),\Delta)$ be a Novikov complex and $\mathcal{F}=\{\mathcal{F}_{p}\mathcal{N}\}$ be a filtration on this complex defined by 
\begin{equation}\label{eq:filtration}
F_{p}\mathcal{N} = \mathbb{Z}((t))\Big[h^1_{k_1},h^{2}_{k_2}\ldots , h^{p+1}_{k_{p+1}}\Big].
\end{equation}
Note that, for each $p\in \mathbb{Z}$, there is only one singularity in $F_p\mathcal{N}\setminus F_{p-1}\mathcal{N}$, hence the filtration $F$ is called a {\it finest filtration}.  The filtration $F$ is \emph{convergent}, i.e. $\cap_p F_{p} \mathcal{N}=0$ and $\cup F_{p} \mathcal{N}=\mathcal{N}$. In fact, $F$ is \emph{finite}, that is, $F_p \mathcal{N}=0$ for some
$p$ and $F_{p^{\prime}}\mathcal{N}=\mathcal{N}$ for some $p^{\prime}$. Moreover, the filtration $F$ is \emph{bounded below}
(i.e. given $q$, there exists $\epsilon(q)$ such that $F_{\epsilon(q)}\mathcal{N}_q=0$) and hence, the induced filtration on $H_*(\mathcal{N})$ is also bounded below.

Since the filtration $F$  on the chain complex $(\mathcal{N}_{\ast}(f),\Delta)$ is  bounded below and  convergent, then there exists a convergent spectral sequence with
$$E^{0}_{p,q}=F_p\mathcal{N}_{p+q}/F_{p-1}\mathcal{N}_{p+q}=G(\mathcal{N})_{p,q}$$ $$E^1_{p,q}\thickapprox H_{(p+q)}(F_p\mathcal{N}_{p+q}/F_{p-1}\mathcal{N}_{p+q})$$ and $E^{\infty}$ is isomorphic to the module $GH_*(\mathcal{N})$. This result  associates a spectral sequence to a filtered chain complex and its proof can be found in \cite{Sp}.
The algebraic formulas for the modules are $$E^{r}_{p,q}=Z^{r}_{p,q}/(Z^{r-1}_{p-1,q+1}+\partial
Z^{r-1}_{p+r-1,q-r+2}),$$ where  $$Z^r_{p,q}=\{c\in F_p\mathcal{N}_{p+q}\,|\,\,
\partial c\in F_{p-r}\mathcal{N}_{p+q-1}\}.$$

Whenever the filtration considered is a finest filtration $F$, the only $q$ such that $E^r_{p,q}$ is non-zero is $q=k-p$. Hence, we omit reference to $q$, i.e. $E^r_p$ is in fact $E^r_{p,k-p}$.

Note that, $E^{\infty}$ does not determine $H_*(\mathcal{N})$ completely, but
$$ E^{\infty}_{p,q}\approx GH_*(\mathcal{N})_{p,q}=\frac{F_pH_{p+q}(\mathcal{N})}{F_{p-1}H_{p+q}(\mathcal{N})}.$$

However, it is a well known fact \cite{D} that whenever $GH_*(\mathcal{N})_{p,q}$ is free and the filtration is bounded, \begin{equation}\label{eq:seqspec}\displaystyle\bigoplus_{p+q=k}GH_*(\mathcal{N})_{p,q}\approx H_{p+q}(\mathcal{N}).\end{equation}

In the more general setting of a circle valued function  on an $n$-dimensional manifold, (\ref{eq:seqspec}) need not be true. However, as a consequence of Theorem \ref{teo:primarypivots}, which will be proved in Section \ref{caracterizacao},
in the $2$-dimensional setting, $GH_*(\mathcal{N})_{p,q}$ is free for all $p$ and $q$ and thus (\ref{eq:seqspec}) holds.
In this case, it follows that the spectral sequence associated to the filtered Novikov complex $(\mathcal{N},\Delta)$, defined by the pair $(f,M)$, converges to the Novikov homology of $M$:
$$ E^{\infty}_{p,q}\approx GH_*(\mathcal{N})_{p,q}=\frac{F_pH_{p+q}(\mathcal{N})}{F_{p-1}H_{p+q}(\mathcal{N})} \approx H^{Nov}_ {\ast}(M,f).$$

\section{Spectral Sequence Sweeping Algorithm for a Novikov Complex}\label{sweeping}

A square matrix will be called a {\it Novikov matrix} if it  is a strictly upper triangular matrix  with square zero and entries in the ring $\mathbb{Z}((t))$. In particular, a Novikov differential is a Novikov matrix.

In this section,  we introduce the Spectral Sequence Sweeping Algorithm (SSSA) for a Novikov matrix associated to a Novikov complex $(\mathcal{N}_{\ast}(f),\Delta)$ on an orientable surface.
 The SSSA constructs
  a family of Novikov matrices $\{\Delta^r, r\geq 0\}$ recursively, where
$\Delta^0=\Delta$, by considering at each stage the $r$-th
diagonal.

\vspace{0.3cm}

{\bf Spectral Sequence Sweeping Algorithm - SSSA}\\
For a fixed $r$-th diagonal\footnote{By $r$-th diagonal one means the collection of entries $\Delta_{ij}$ of $\Delta$ such that $j-i=r$.}, the method described below must be applied for all $\Delta_k$ for $k=0,1,2$ simultaneously.

\subsection*{A - Initial step} 
\begin{itemize}
\item[]\begin{enumerate}

\item Without loss  of generality, we assume that  the first diagonal of $\Delta$ contains non-zero entries $\Delta_{{i,j}}$ where $j\in J_k$ and $i\in J_{k-1}$. 
Whenever the first diagonal contains only zero entries, we define $\Delta^{1}=\Delta$ and we repeat this step until we reach a diagonal of  $\Delta$ which contains non-zero entries. 

The non-zero entries  $\Delta_{{i,j}}$ of the first diagonal  are called  \emph{index $k$ primary pivots}. It follows that    
 the entries $\Delta_{{s ,j}}$ for
$s>i$ are all zero.

We end this first step by defining $\Delta^{1}$ as $\Delta$ with
the index $k$ primary pivots on the first diagonal
marked.

\item Consider the matrix $\Delta^{1}$ and let
$\Delta_{{i,j}}^{1}$ be the entries in $\Delta^{1}$ where
the $i\in J_{k-1}$ and $j\in J_k$. 
Analogously to step one, we assume without loss of generality that the second diagonal contains non-zero entries $\Delta_{{i,j}}^{1}$. We now construct a
matrix $\Delta^{2}$ following the procedure:

Given a non-zero entry $\Delta_{{i,j}}^{1}$ on the second
 diagonal of $\Delta^{1}$
\begin{enumerate}
\item if there are no
primary pivots in row $i$ and column $j$, mark it as
an index $k$ primary pivot and the numerical value of the entry
remains the same, i.e.
$\Delta_{{i,j}}^{2}=\Delta_{{i,j}}^{1}$.

\item if this is not the case, consider the entries in column $j$ and in a row $s$
with $s>i$ in $\Delta^{1}$.

\begin{itemize}
\item[(b1)] If there is an index $k$ primary pivot in an entry in 
column $j$ and in a row $s$, with $s>i$, then the numerical value
remains the same and the entry is left unmarked, i.e.
$\Delta_{{i,j}}^{2}=\Delta_{{i,j}}^{1}$.

\item[(b2)] If there are no primary pivots in 
column $j$ below $\Delta_{{i,j}}^{1}$ then there is an index $k$
primary pivot in row $i$, say a column $u$ of
$\Delta^{1}$, with $u<j$. In this case, we define  $\Delta_{{i,j}}^{2}=\Delta_{{i,j}}^{1}$  and this  entry is marked as a \emph{change of basis pivot}.

Note that we have defined a matrix $\Delta^{2}$ which is
actually equal to $\Delta^{1}$ except that the second
diagonal is marked with primary and change of basis pivots.
\end{itemize}
\end{enumerate}
\end{enumerate}
\end{itemize}

\subsection*{B - Intermediate step}
\begin{itemize}
\item[]In this step we consider a matrix $\Delta^r$ with the primary and change of basis pivots marked on the
$\xi$-th diagonal for all $\xi\leq r$. We now
describe how $\Delta^{r+1}$ is defined.  If there does not exist a
change of basis pivot on the $r$-th diagonal we go
directly to step B.2, that is, we define $\Delta^{r+1}=\Delta^r$ with 
diagonal $(r+1)$ marked with primary and change of
basis pivots as in B.2.
\end{itemize}
\subsubsection*{B.1 - Change of basis}
\begin{itemize}
\item[]Suppose $\Delta^r_{{i,j}}$
is a change of basis pivot on the $r$-th diagonal. Since
we have a change of basis pivot in row $i$, there is a
column, namely $u$-th column, associated to a $k$-chain such that
$\Delta^r_{{i,u}}$ is a primary pivot. Then, perform a change of
basis on $\Delta^{r}$ in order to zero out the entry
$\Delta_{{i,j}}^{r}$ without introducing non-zero entries in
$\Delta^{r}_{{s,j}}$ for $s>i$. 
 We will prove in Theorem \ref{teo:primarypivots}  that all the entries in $\Delta^r$ which are primary pivots are equal to $\pm t^{l_{1}} \pm t^{l_{2}}$ and, since these entries are invertible in $\mathbb{Z}((t))$, it is always possible choosing a particular change of basis using only columns $j$ and $u$ of $\Delta^r$.

Once this is done, we obtain a $k$-chain associated to 
column $j$ of $\Delta^{r+1}$. It is a linear combination over
$\mathbb{Z}((t))$ of columns $u$ and $j$ of $\Delta^{r}$  such that $\Delta_{{i,j}}^{r+1}=0$.
It is also a particular linear combination of the columns of $\Delta$ in $J_k$ on and to the left of column $j$.

Let $k_1$ be the column of $\Delta^{r}$ which is associated to a $k$-chain. We denote by
$\sigma^{j,r}_{k}$  the Morse index $k$-chain corresponding to column $j$ of $\Delta^r$. We have $$\sigma^{j,r}_{k}=\sum_{\ell=k_1}^{j}c^{j,r}_{\ell}h_k^{\ell}, $$
and  
 column $j$ of $\Delta^{r+1}$ corresponds to 
\begin{equation}\label{eq:algorithm}
\sigma^{j,{r+1}}_{k}=\sigma^{j,r}_{k} - \Delta^{r}_{{i,j}} (\Delta^{r}_{{i,u}})^{-1}  \sigma^{u,r}_{k}
=c^{j,r+1}_{k_1}h_k^{k_1} 
+
\cdots
+c^{j,r+1}_{j-1}h_k^{j-1}+ c^{j,r+1}_{j}h_k^{j}\end{equation}
where $c^{j,r+1}_{\ell}\in \mathbb{Z}((t))$ and $c^{j,r+1}_{j}=1$.

Therefore the matrix $\Delta^{r+1}$ has entries determined
by a change of basis over $\mathbb{Z}((t))$ of $\Delta^r$. In
particular, all the change-of-basis pivots on the $r$-th
diagonal of $\Delta^r$ are zero in $\Delta^{r+1}$.

Once the above procedure is done for all change-of-basis pivots of the $r$-th diagonal of $\Delta^{r}$, we can define a change-of-basis
matrix $T^{r}$ such that $\Delta^{r+1}=(T^{r})^{-1}\Delta^{r} T^{r}$.

\end{itemize}
\subsubsection*{B.2 - Marking the  diagonal $(r+1)$ of $\Delta^{r+1}$}

\begin{itemize}\item[]Consider the matrix $\Delta^{r+1}$ defined in the previous step. We mark  diagonal $(r+1)$ with primary and
change-of-basis pivots as follows:

Given a non-zero entry $\Delta_{{i,j}}^{r+1}$
\begin{enumerate}
\item If there are no
primary pivots in row $i$ and  column $j$, mark it as
an index $k$ primary pivot.

\item If this is not the case, consider the entries in column $j$ and in a row $s$
with $s>i$ in $\Delta^{r+1}$.

\begin{itemize}
\item[(b1)] If there is an index $k$ primary pivot in the entries in column $j$
 below $\Delta_{{i,j}}^{r+1}$ then leave the entry unmarked.

\item[(b2)] If there are no primary pivots in 
column $j$ below $\Delta_{{i,j}}^{r+1}$ then there is an index $k$
primary pivot in row $i$, say in the column $u$ of
$\Delta^{r+1}$, with $u<j$. In this case, mark it as a change-of-basis pivot.
\end{itemize}

\end{enumerate}
\end{itemize}

\subsection*{C - Final step}
\begin{itemize}\item[]We repeat the above procedure until all diagonals
have been considered.
\end{itemize}

\vspace{0.5cm}

Note that in the SSSA the columns of the matrix
$\Delta$ are not necessarily ordered with respect to $k$, or equivalently, that the
singularities $h_k$ are not ordered with respect to the filtration. In this work, without loss of generality,  we consider the singularities to be ordered with respect to the Morse index for the sole reason of  simplifying  notation.

In order to perform  the particular change of basis (\ref{eq:algorithm}) in step $B.1$ of the SSSA,  the primary pivots must be invertible polynomials in the ring $\mathbb{Z}((t))$. Otherwise, the change of basis in $(\ref{eq:algorithm})$ is not well defined. The example below shows a Novikov differential for which  all changes of bases are well defined and hence the SSSA is correct. 

\begin{ex}\label{exemplo}

Figure \ref{fig:exemplo} illustrates a flow on the torus $T^{2}$ associated to a circle-valued Morse function $f$ defined on $T^{2}$, where $a$ is a regular value of $f$. 
The Novikov chain groups are generated by the critical points as follows:
$\mathcal{N}_{0}=\{h_{0}^{1},h_{0}^{2}\}$,  
$\mathcal{N}_{1}=\{h_{1}^{3},h_{1}^{4},h_{1}^{5},h_{1}^{6}\}$
and $\mathcal{N}_{2}=\{h_{2}^{7},h_{2}^{8}\}$.
Choosing orientations for the unstable manifolds of the critical points of $f$ as indicated in Figure \ref{fig:exemplo}, the Novikov matrix associated to $\partial$ is presented in Figure \ref{fig:delta0ex2}.

\begin{figure}[!htb]
\centering
\begin{minipage}[b]{0.5\linewidth}
\includegraphics[width=7.5cm]{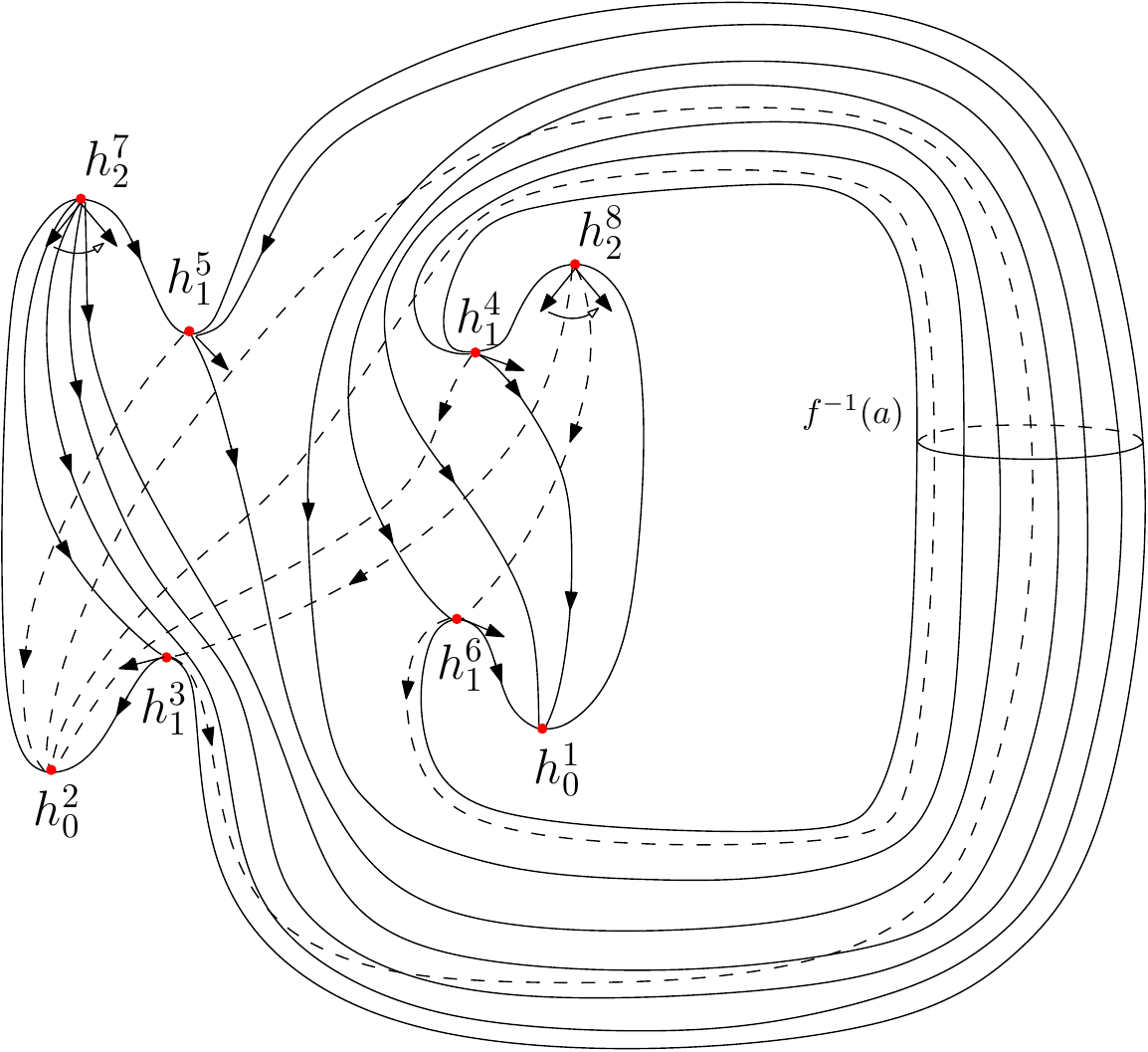}
\caption{Circle-valued Morse function on the torus.}
\label{fig:exemplo}
\end{minipage} 
\begin{minipage}[b]{0.45\linewidth}
\includegraphics[width=7cm]{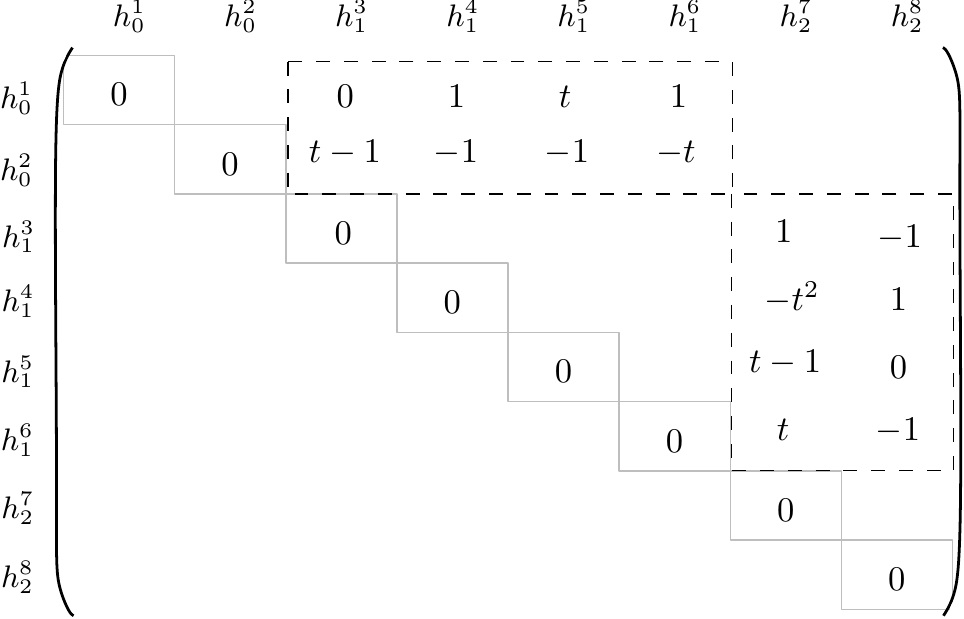}
\vspace{0.3cm}
\caption{Novikov matrix.}
\label{fig:delta0ex2}
\end{minipage}
\end{figure}

Applying the SSSA to the Novikov matrix $\Delta$ in Figure \ref{fig:delta0ex2}, one  obtains the sequence of Novikov matrices $\Delta^{1},\cdots,\Delta^{6}$ presented in Figures \ref{fig:delta1ex}, $\cdots$,  \ref{fig:delta6ex}, respectively. In these figures,  the markup process at the $r$-th iteration is done as follows: primary pivots are encircled and change-of-basis pivots are
encased in boxes. 

Note that, in this example, each marked primary pivot in $\Delta^{r}$ is invertible in the ring $\mathbb{Z}((t))$, making it possible to apply the SSSA and obtain the next Novikov matrix $\Delta^{r+1}$.

\begin{figure}[!htb]

\begin{minipage}[t]{0.5\linewidth}
 \includegraphics[width=7.5cm]{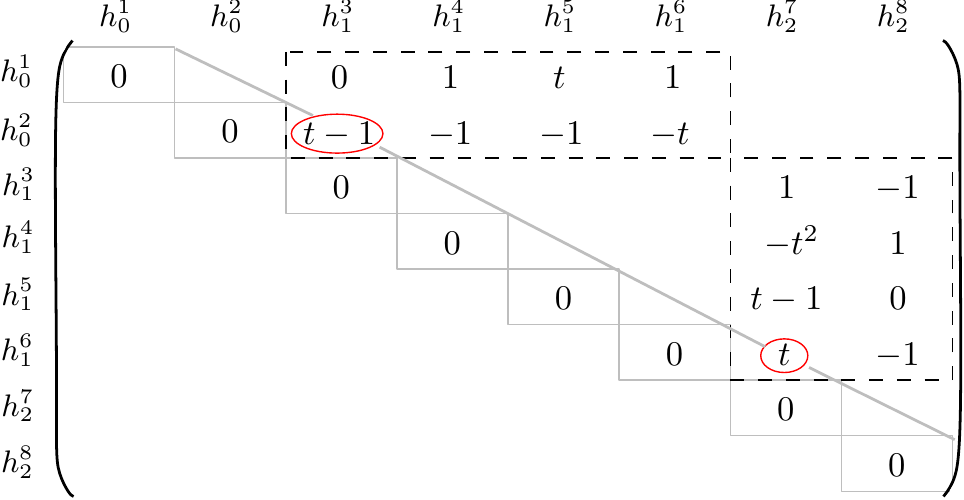}\\
{\footnotesize  $\sigma^{j,1}_{0}\!=\!h^j_{0}$ for $\!j\in \! J_{0}$;   $\sigma^{j,1}_{1}\!=\!h^j_{1}$  for $j\!\in \! J_{1}$;   $\sigma^{j,1}_{2}\!=\! h^j_{2}$  for $j\!\in \! J_{2}$.}\vspace{-0.3cm}
\caption{$\Delta^{1}$, marking primary pivots.}
\label{fig:delta1ex}
\end{minipage} 
\begin{minipage}[t]{0.5\linewidth}
 \includegraphics[width=7.5cm]{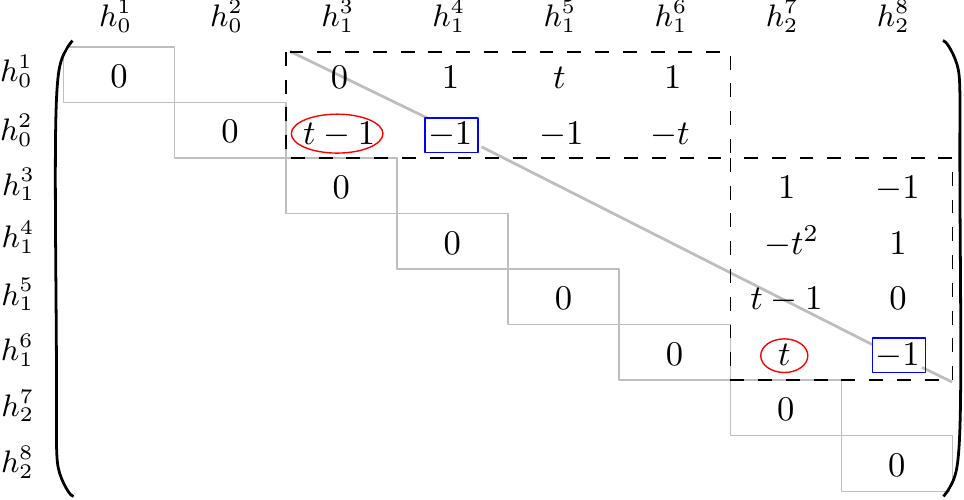}\\
  {\footnotesize    $\sigma^{j,2}_{0}\!=\!h^j_{0}$ for $\!j\in \! J_{0}$;   $\sigma^{j,2}_{1}\!=\!h^j_{1}$  for $j\!\in \! J_{1}$;   $\sigma^{j,2}_{2}\!=\! h^j_{2}$  for $j\!\in \! J_{2}$.}\vspace{-0.3cm}
\caption{$\Delta^{2}$, marking change-of-basis pivot.}
\label{fig:delta2ex}
\end{minipage}
\end{figure}
\begin{figure}[!htb]
\begin{minipage}[t]{0.5\linewidth}
\includegraphics[width=7.5cm]{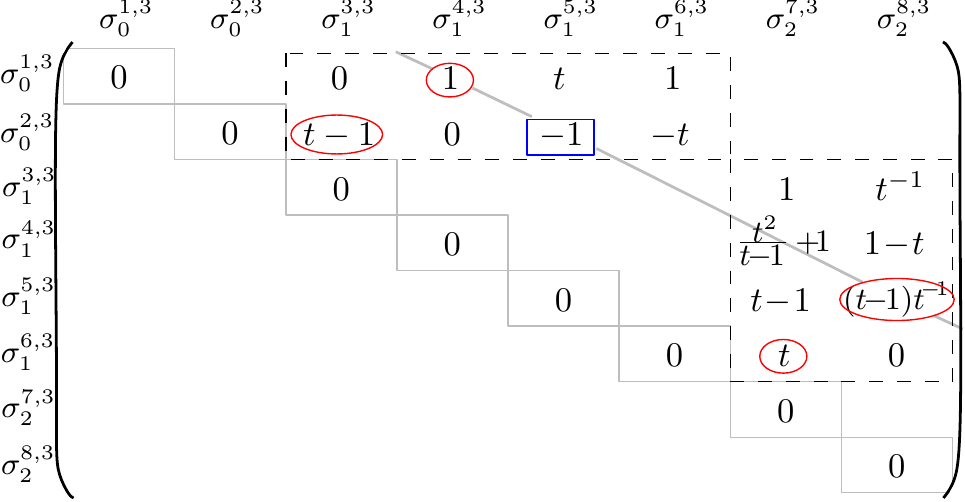}\\
 {\footnotesize    $\sigma^{4,3}_{1}\!=\!h^4_{1} +(t\!-\!1)^{-1}h^{3}_{1}$;    $\sigma^{8,3}_{1}\!=\!h^{8}_{2} + t^{-1}h^{7}_{2}$. \\
  $\sigma^{j,3}_{k} = \sigma^{j,2}_{k}$ for all the  remaining  $\sigma$'s.}\vspace{-0.3cm}
\caption{$\Delta^{3}$, marking pivots.}
\label{fig:delta3ex}
\end{minipage} 
\begin{minipage}[t]{0.5\linewidth}
\includegraphics[width=7.5cm]{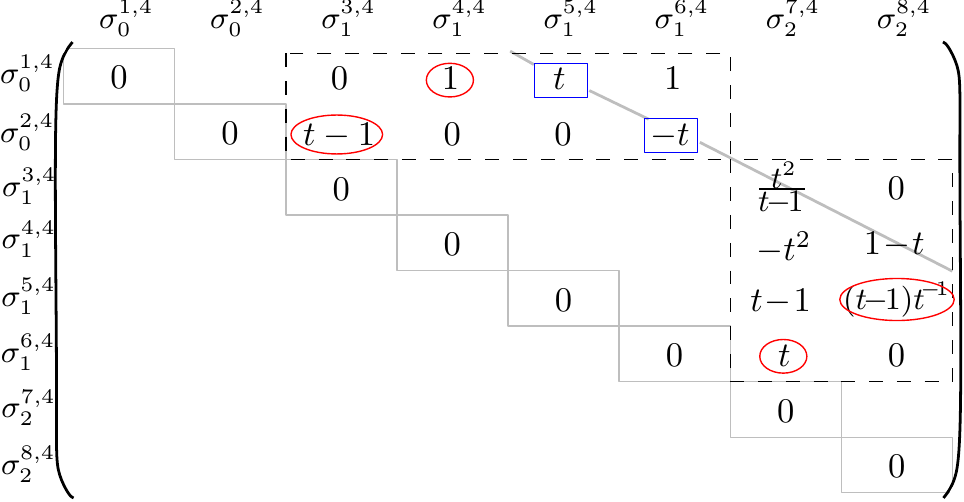}\\
   {\footnotesize    $\sigma^{5,4}_{1}\!=\!h^{5}_{1} + (t\!-\! 1)^{-1}h^{3}_{1}$. \\
  $\sigma^{j,4}_{k} = \sigma^{j,3}_{k}$ for all the  remaining  $\sigma$'s.}\vspace{-0.3cm}
\caption{$\Delta^{4}$, marking change-of-basis pivot.}
\label{fig:delta4ex}
\end{minipage}
\end{figure}
\begin{figure}[!htb]
\begin{minipage}[t]{0.5\linewidth}
\includegraphics[width=7.5cm]{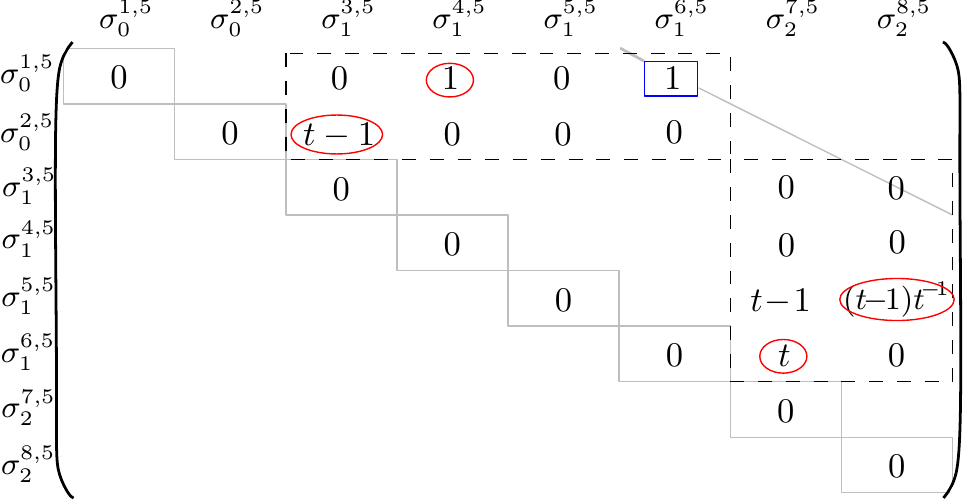}\\
     {\footnotesize    $\sigma^{5,5}_{1}\!=\!h^{5}_{1} -t \ h^{4}_{1}$;    $\sigma^{6,5}_{1}\!=\!h^{6}_{1} + t(t\!-\!1)^{-1}h^{3}_{1}$;   \\
  $\sigma^{j,5}_{k} = \sigma^{j,4}_{k}$ for all the  remaining  $\sigma$'s.}\vspace{-0.3cm}
\caption{$\Delta^{5}$, marking change-of-basis pivot.}
\label{fig:delta5ex}
\end{minipage} 
\begin{minipage}[t]{0.5\linewidth}
\includegraphics[width=7.5cm]{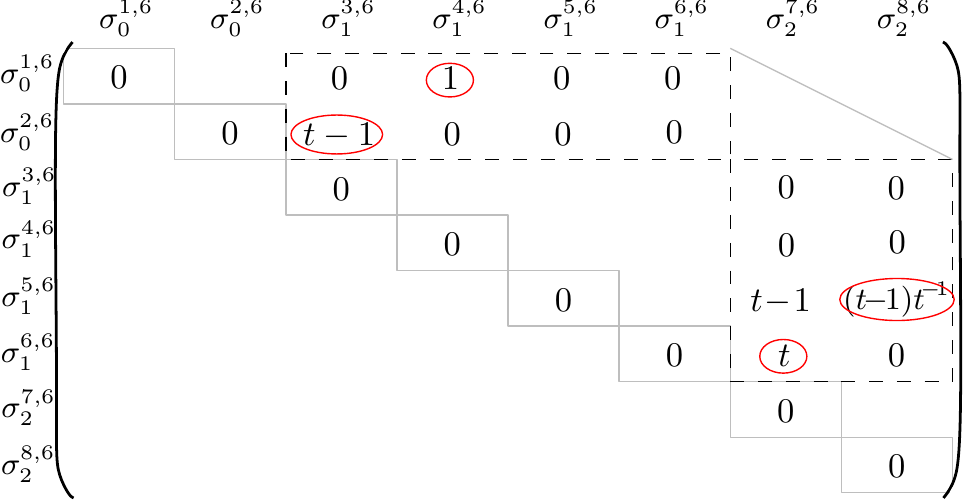}\\
       {\footnotesize    $\sigma^{6,6}_{1}\!=\!h^{6}_{1} - \ h^{4}_{1}$;       \\
  $\sigma^{j,6}_{k} = \sigma^{j,5}_{k}$ for all the  remaining  $\sigma$'s.}\vspace{-0.3cm}
\caption{$\Delta^{6}$ for Example \ref{exemplo}.}
\label{fig:delta6ex}
\end{minipage}
\end{figure}

\end{ex}

\section{Characterization of the Novikov Matrices}\label{caracterizacao}

The primordial aim in this section is to show  that the SSSA is proved to be correct for all Novikov differentials of a 2-dimensional Novikov complex. This is done by showing that, given a Novikov differential $\Delta$, all primary pivots determined by the SSSA are invertible polynomials in the ring $\mathbb{Z}((t))$. In fact, they are monomials with coefficient $\pm 1$ or binomials with coefficients $\pm 1$.  Throughout this section, the term monomial (respectively, binomial) will be used to refer to polynomials in $\mathbb{Z}((t))$ of the form $\pm t^{\ell}$ (respectively,    $t^{\ell_{1}} - t^{{\ell_{2}}}$), where $\ell \in \mathbb{Z}$ (respectively, $\ell_{1}, \ell_{2} \in \mathbb{Z}$).

\begin{teo}\label{teo:primarypivots}
Given a Novikov differential $\Delta$, the primary pivots and the change of basis pivots in the SSSA are polynomials of the form $t^{\ell}$ or $ t^{\ell_1} - t^{\ell_2}$, where $\ell,\ell_{1},\ell_{2}\in \mathbb{Z}$. 
\end{teo}

The proof  of Theorem \ref{teo:primarypivots} is an immediate consequence of Theorem \ref{bloco1}, Lemma \ref{Tinverso} and Theorem \ref{bloco2}. Theorem \ref{bloco1} gives a characterization of the columns of the first block $J_{0}\times J_{1}$ of $\Delta$ and Theorem \ref{bloco2} gives a characterization of the rows of the second block $J_{0}\times J_{1}$ of $\Delta$.

Lemma \ref{lemapivolinha} asserts that we cannot have more than one primary pivot in a fixed row
or column. Moreover, if there is a primary pivot in row $i$, then there is no primary pivot in
column $i$.

\begin{lem}\label{lemapivolinha} 
Let $\Delta$ be a Novikov differential for which the SSSA is proved correct up to  step $R$. Let $\Delta^{1},\cdots,\Delta^{R}$ be the family of Novikov matrices produced by the SSSA until step $R$. Given two primary pivots, the $ij$- th entry $\Delta^{r}_{{i,j}}$ and  the $ml$-th entry $\Delta^{r}_{{m,l}}$, then $\{i,j\}\cap \{m,l\}=\emptyset$.
\end{lem}

We omit the proof of Lemma \ref{lemapivolinha}, since it is similar in nature to  proof of Proposition 3.2 in \cite{CdRS}, where the SSSA was designed for a Morse chain complex over $\mathbb{Z}$. The next lemma implies that, in order to know the pivots which will appear during the execution of the SSSA, one  can apply this algorithm separately in block $J_{0}\times J_{1}$ and $J_{1}\times J_{2}$.

\begin{lem}\label{Tinverso} 
Let $\Delta$ be a Novikov differential for which the SSSA is proved correct up to  step $R$. Then,
the change of basis caused by change-of-basis pivots in block $J_{0} \times J_{1}$ do not affect the pivots in block $J_{1}\times J_{2}$. In other words, 
 multiplication by $(T^r)^{-1}$ does not change the primary and change-of-basis  pivots in block  $J_{1}\times J_{2}$ .
\end{lem}
\dem
 Without loss of generality, suppose that there is only one change-of-basis pivot $\Delta^{r}_{i,j}$ in $\Delta^{r}$ with $j\in J_{1}$, and let $\Delta^{r}_{i,u}$ be the primary pivot in row i. The change of basis matrix $T^{r}$ has unit diagonal and the only non zero entry off the diagonal is $T^{r}_{u,j} = -\Delta^{r}_{i,j}(\Delta^{r}_{i,u})^{-1}$. Hence, $(T^{r})^{-1}$ has unit diagonal and the only non zero entry off the diagonal is $(T^{r})^{-1}_{u,j} = - T^{r}_{u,j} =\Delta^{r}_{i,j}(\Delta^{r}_{i,u})^{-1}$. Therefore, multiplication by $(T^{r})^{-1}$ will only affect  row $u$ of $\Delta^{r}$.
By Lemma \ref{lemapivolinha},  there are no  primary pivots in row $u$ and hence there are no change-of-basis pivot in row $u$ as well. 
\cqd

Before proving  Theorem \ref{teo:primarypivots}, we introduce the notation and terminology  that will be used in the proof. From now on, we consider the SSSA without realising the pre-multiplication by  $(T^{r})^{-1}$, unless mention otherwise.

\begin{figure}[!ht]
\centering
\includegraphics[scale=1]{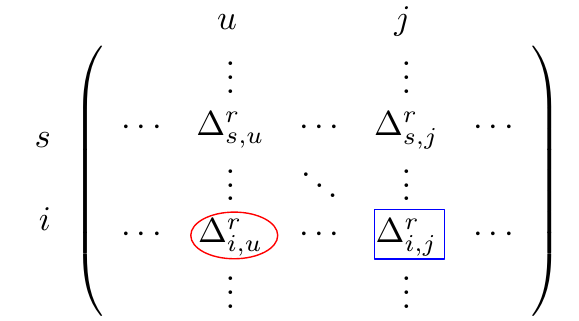} \\ \vspace{-0.5cm}
\caption{$\Delta^{r}_{J_{k-1}J_{k}}$; marking change-of-basis pivot.}\label{fig:generalmatrix}
\end{figure}

Let $\Delta^{r}_{i,j}$ be a change-of-basis pivot caused by a primary pivot $\Delta^{r}_{i,u}$.
 Suppose a change of basis determined by  $\Delta^{r}_{i,j}$ is performed by the SSSA in the matrix $\Delta^{r}$, i.e.,  in the step $(r+1)$, one has 
$$\sigma^{j,{r+1}}_{k}=\sigma^{j,r}_{k} - \Delta^{r}_{{i,j}} (\Delta^{r}_{{i,u}})^{-1}  \sigma^{u,r}_{k}.$$
Hence, whenever this change of basis  occurs, only column $j$ of the matrix $\Delta^{r}$ is modified, in fact, for each $s=0,\dots,i$, the entry $\Delta^{r}_{s,j}$  is  added to a  multiple of the entry $\Delta^{r}_{s,u}$,  where $u$  is the column of the primary pivot  in row $i$. In other words, $\Delta^{r+1}_{s,j} = \Delta^{r}_{s,j} -       \Delta^{r}_{{i,j}} (\Delta^{r}_{{i,u}})^{-1} \Delta^{r}_{s,u}$. See the matrix in Figure \ref{fig:generalmatrix} (this figure shows part of the block associated with index $k$, as the $r$-th diagonal is swept).

\begin{defin}\label{def:creation} In the situation described above and represented in Figure \ref{fig:generalmatrix}, we assert that the entry $\Delta^{r}_{s,u}$ in column $u$ {\it generates}\footnote{If an entry $\Delta^{\xi}_{s,u}$  with $\xi <r$ does not change until step $r$, i.e. $\Delta^{\xi}_{s,u}= \Delta^{\xi +1}_{s,u} = \cdots = \Delta^{r}_{s,u}$, and $\Delta^{r}_{s,u}$ generates $\Delta^{r+1}_{s,j}$, we say that $\Delta^{\xi}_{s,u}$ generates $\Delta^{r+1}_{s,j}$.}  the entry $\Delta^{r+1}_{s,j}$ in column $j$ in  $\Delta^{r+1}$, whenever $\Delta^{r}_{s,u}\neq 0$.
\end{defin}

Note that if an entry in a column $j$  generates another entry in a column $t$ then $t>j$, i.e, $\Delta^{r}_{s,j}$ generates an entry in a column on the right of  the column $j$.

\begin{figure}[!ht]
\centering
\includegraphics[scale=1]{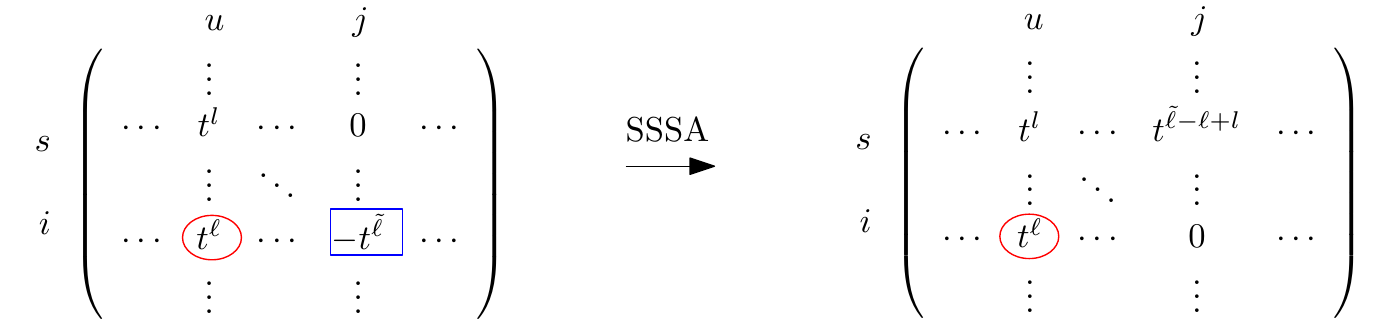}
\vspace{-0.4cm}
\caption{$\Delta^{r}_{J_{k-1}J_{k}}$ and $\Delta^{r+1}_{J_{k-1}J_{k}}$, respectively.}\label{fig:generalmatrix1}
\end{figure}

\begin{figure}[!ht]
\centering
\includegraphics[scale=1]{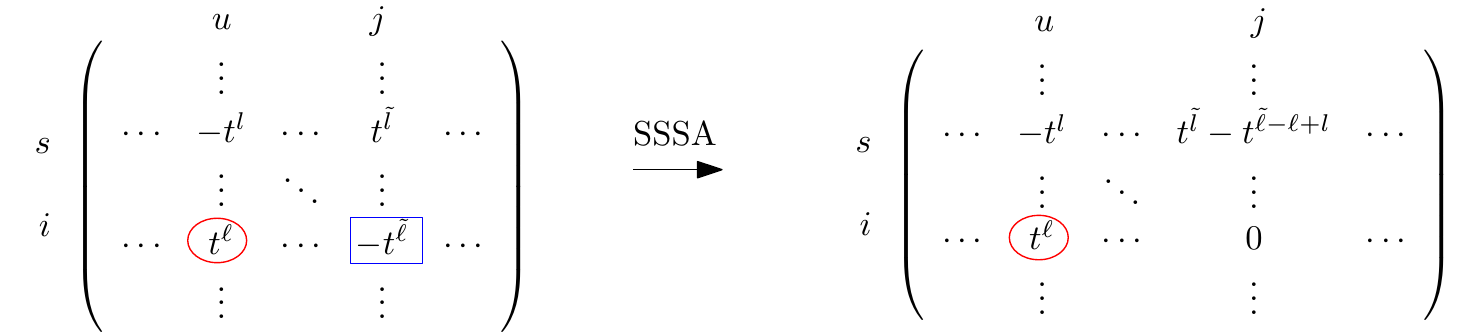}
\vspace{-0.4cm}
\caption{$\Delta^{r}_{J_{k-1}J_{k}}$ and $\Delta^{r+1}_{J_{k-1}J_{k}}$, respectively.}\label{fig:generalmatrix2}
\end{figure}

\begin{figure}[!ht]
\centering
\includegraphics[scale=1]{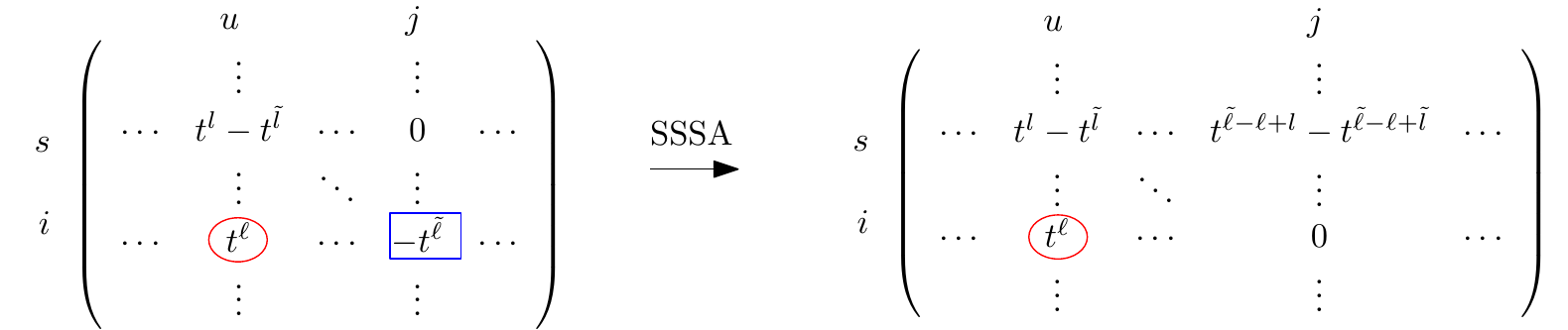}
\vspace{-0.4cm}
\caption{$\Delta^{r}_{J_{k-1}J_{k}}$ and $\Delta^{r+1}_{J_{k-1}J_{k}}$, respectively.}\label{fig:generalmatrix3}
\end{figure}

\begin{ex}
It is helpful to keep in mind  some  configurations that allow an entry $\Delta^{r}_{s,u}\neq 0$ to generate another entry $\Delta^{r+1}_{s,j}$. 
Consider for instance that $\Delta^{r}_{i,u}=t^{\ell}$ and $\Delta^{r}_{i,j}=-t^{\tilde{\ell}}$. We list some of the possibilities for the entries in positions $(s,u)$ and $(s,j)$ of $\Delta^{r}$:
\begin{enumerate}
\item $\Delta^{r}_{s,u} = t^{l}$ and $\Delta^{r}_{s,j} = 0$.
In this case, $\Delta^{r+1}_{s,j} = \Delta^{r}_{s,j} -       \Delta^{r}_{{i,j}} (\Delta^{r}_{{i,u}})^{-1} \Delta^{r}_{s,u} =  t^{\tilde{\ell}}t^{-\ell}t^{l}$, see Figure \ref{fig:generalmatrix1}.
\item $\Delta^{r}_{s,u} = -t^{l}$ and $\Delta^{r}_{s,j} = t^{\tilde{l}}$.
In this case, $\Delta^{r+1}_{s,j} = \Delta^{r}_{s,j} -       \Delta^{r}_{{i,j}} (\Delta^{r}_{{i,u}})^{-1} \Delta^{r}_{s,u} = t^{\tilde{l}} - t^{\tilde{\ell}}t^{-\ell}t^{l}$, see Figure \ref{fig:generalmatrix2}.
\item $\Delta^{r}_{s,u} = t^{l}-t^{\tilde{l}}$ and $\Delta^{r}_{s,j} = 0$.
In this case, $\Delta^{r+1}_{s,j} = \Delta^{r}_{s,j} -       \Delta^{r}_{{i,j}} (\Delta^{r}_{{i,u}})^{-1} \Delta^{r}_{s,u} =  t^{\tilde{\ell}}t^{-\ell}(t^{l} - t^{\tilde{l}})$, see Figure \ref{fig:generalmatrix3}.
\end{enumerate} 
\noindent Figures \ref{fig:generalmatrix1}, \ref{fig:generalmatrix2} and \ref{fig:generalmatrix3} show part of the block associated with index $k$, as the $r$-th diagonal is swept.
\end{ex}


In fact, we will see that,  as a consequence of Theorem \ref{bloco2}, the cases shown in the previous example are the only possibilities up to sign of generating entries by  a change of basis in block $J_{1}\times J_{2}$. Of course, in block $J_{0}\times J_{1}$ there are more possibilities.

\begin{lem}\label{creation}
An entry which is or will be marked as a primary or a change-of-basis  pivot never generates  entries.
\end{lem}

\dem
If an entry $\Delta^{r}_{s,u}$ generates an entry in  $\Delta^{r+1}$, then there must be a primary pivot in column $u$ and row $i>s$, which was marked in step $\xi < r$. Hence, $\Delta^{r}_{s,u}$ can not be marked as a pivot in any given step $r$. 
\cqd

\begin{defin} Let $\Delta$ be a Novikov differential.
\begin{enumerate}
\item[(a)]
When an entry $\Delta^{r}_{s,u}$ generates another one $\Delta^{r+1}_{s,j}$, we say that $\Delta^{r+1}_{s,j}$ is an {\it immediate successor} of $\Delta^{r+1}_{s,u}$.
\item[(b)] A sequence of entries $\{ \Delta^{\xi_{0}}_{s,j_{0}}, \Delta^{\xi_{2}}_{s,j_{2}},\cdots, \Delta^{\xi_{f}}_{s,j_{f}} \}$  such that  each entry is an immediate successor of the previous one is called a {\it generation sequence}.
\item[(c)] Given an entry $\Delta_{s,u}^{\xi}$ of $\Delta^{\xi}$, the {\it $\Delta_{s,u}^{\xi}$-lineage} is defined to be the set of all   generation sequences  whose first element is $\Delta_{s,u}^{\xi}$.
\end{enumerate}
\end{defin}

We will say that all the  elements in these sequences are in the same lineage or in  $\Delta^{\xi}_{s,u}$-lineage. Also, an  element of   a generation sequence is said to be successor of every element   of this sequence which is to its left.

\begin{lem}
Let $\Delta$ be a Novikov differential for which the SSSA is proved correct up to step $R$.  If $\Delta$ has the property that  at most one change-of-basis pivot is marked in a row during the SSSA until step $R$, then every lineage is formed by  a unique  generation sequence. 
\end{lem}
\dem
By  hypothesis,  one has that in each row $i$ at most one change-of-basis pivot is marked through out the algorithm and if so the mark up is done in  step  $2\leq \xi \leq m-1$, where $m$ is the order of $\Delta$. Then an entry $\Delta^{\xi_{0}}_{s,j_{0}}$ in row $s$ generates at most one entry $\Delta^{\xi_{1}}_{s,j_{1}}$ through out the algorithm,   and this entry  will  necessarily be in a column $j_{1} > j_{0}$ and diagonal $\xi_{1}$ with $\xi_{0}< \xi_{1}\leq m-1$.  In fact, if $\Delta^{\xi_{0}}_{s,j_{0}}$ generates two entries,  then either there would be  two change-of-basis pivots in row $i$, which contradicts our initial  hypothesis,   or two primary pivots in column $j_{0}$, which can not occur by the definition of primary pivots. 

Now, if $\xi_{1} \leq  m-1$, then $\Delta^{\xi_{1}}_{s,j_{1}}$ can  generate  at most a unique entry $\Delta^{\xi_{2}}_{s,j_{2}}$ where $ \xi_{1}<\xi_{2} \leq m-1$ and $j_{2}> j_{1}$ and this can be done successively. More specifically,
the entry $\Delta^{\xi}_{s,u}$ is responsible for  generating a unique immediate successor and this successor can in turn generate a unique immediate successor and  thus it determines a full lineage of entries represented in one finite sequence  
 $\{  \Delta^{\xi_{0}}_{s,j_{0}}, \Delta^{\xi_{1}}_{s,j_{1}}, \Delta^{\xi_{2}}_{s,j_{2}}, \cdots \}$, where $ 2 \leq  \xi_{0} < \xi_{1} < \xi_{2} < \cdots \leq m -1 $, and $j_{0} < j_{1} < j_{2} < \cdots \leq m $. 
 \cqd

\begin{cor}\label{uniqueseq}
Let $\Delta$ be a matrix for which the SSSA is proved correct up to step $R$. Suppose that  the second block of $\Delta$ has the property that at most one change-of-basis pivot is marked in each row from the beginning until the end of the SSSA. Therefore,  each $\Delta_{s,u}^{\xi}$-lineage with $s \in J_{1}$ is formed by  a unique  generation sequence.
\end{cor}
\dem
The SSSA  applied to the first block $J_{0}\times J_{1}$ of $\Delta$ does not interfere  in the number of change-of-basis pivots identified in block $J_{1}\times J_{2}$.
\cqd

\begin{figure}[!ht]
\centering
\includegraphics[scale=1]{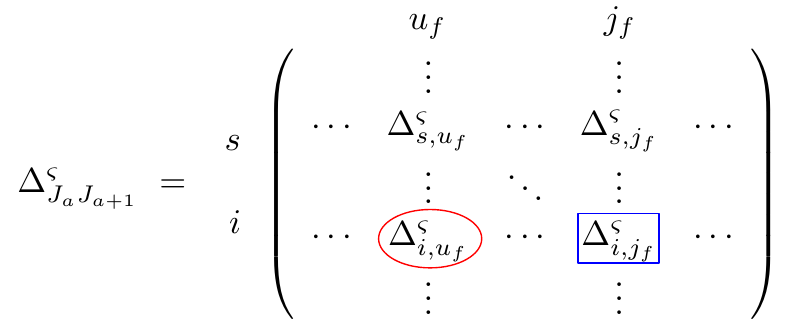} \\ \vspace{-0.5cm}
\caption{Generating a binomial from two monomial lineages.}\label{fig:binomiallineage}
\end{figure}

Consider a  $\Delta_{s,u}^{\xi}$-lineage which is formed by  a unique  generation sequence. If this generation sequence contains only monomials (binomials, resp.) then one says that $\Delta_{s,u}^{\xi}$-lineage is a {\it monomial} ({\it binomial}, resp.) lineage. However, if $\Delta_{s,u}^{\xi}$ is a monomial, then the  $\Delta_{s,u}^{\xi}$-lineage  could eventually contain binomials. One way that this can occur is when row $s$ is of type $3$ in $\Delta$, $\Delta_{s,u}^{\xi}$ and $\Delta_{s,j}^{\zeta} $ are monomials and the lineage determined by these entries merge giving rise to a binomial. More specifically, 
suppose that the first binomial in row $s$ appears in $\Delta^{\varsigma+1}$    then one has two monomial lineages 
 $\{  \Delta_{s,u}^{\xi}, \Delta^{\xi_{1}}_{s,u_{1}}, \Delta^{\xi_{2}}_{s,u_{2}}, \cdots, \Delta^{\xi_{f}}_{s,u_{f}} \}$
 and  $\{  \Delta_{s,j}^{\zeta}, \Delta^{\zeta_{1}}_{s,j_{1}}, \Delta^{\zeta_{2}}_{s,j_{2}}, \cdots, \Delta^{\zeta_{f}}_{s,j_{f}} \}$, where $\xi_{f}, \zeta_{f} \leq \varsigma$; observe  that $\Delta^{\xi_{f}}_{s,u_{f}} = \Delta^{\varsigma}_{s,u_{f}} $  and $\Delta^{\zeta_{f}}_{s,j_{f}} = \Delta^{\varsigma}_{s,j_{f}} $.
  The binomial will appear in $\Delta^{\varsigma+1}$ as a consequence of a change of basis caused by a change-of-basis pivot $\Delta_{i,j_{f}}^{\varsigma}$ and a primary pivot $\Delta_{i,u_{f}}^{\varsigma}$ in a row  $i > s$, as in Figure \ref{fig:binomiallineage}.
In this case, $\Delta^{\varsigma}_{s,u_{f}}$ is the generator  of the binomial $\Delta^{\varsigma + 1}_{s,j_{f}}$. Hence, we say that the 
 $\Delta_{s,j}^{\zeta}$-lineage {\it ceases}, i.e, this lineage remains the same until  $\Delta^{r+1}$  and the  $\Delta_{s,u}^{\xi}$-lineage is an {\it eventual binomial lineage}. From this point on,  this lineage contains only binomials.  

Once an element of a lineage  is marked as a pivot, this lineage ceases, since pivots do not generate entries, by Lemma \ref{creation}.

The next theorem provides a characterization of columns in the first block $J_{0}\times J_{1}$ of the  Novikov matrices $\Delta^{r}$ as the diagonals are swept.

\begin{teo}[First Block Characterization]\label{bloco1} Let $\Delta$ be a Novikov differential for which the SSSA is proved correct up to step $R$. Then we have the following possibilities for a column $j$ of the   matrix $\Delta^{r}$   with $j \in J_{1}$ and  $r \leq R$:
\begin{enumerate}
\item all the entries in column $j$ are equal to zero, .i.e, $\Delta^{r}_{{\bullet, j}}=0$;
\item there is only one non-zero entry in column $j$ and it is a binomial $ t^{\ell} - t^{\tilde{\ell}}$, where $\ell,\tilde{\ell} \in \mathbb{Z}$;
\item there are exactly two non-zero entries in column $j$ and they are monomials $ t^{\ell}$  and $ -t^{\tilde{\ell}}$, where $\ell, \tilde{\ell} \in \mathbb{Z}$; 
\item there is only one non-zero entry in column $j$ and it is a monomial $ t^{\ell}$, where $\ell\in \mathbb{Z}$.
\end{enumerate}
\end{teo}

\dem
The proof is done by induction. Note that the result is trivial for $\Delta^{1}$ and $\Delta^{2}$. In fact, the first change-of-basis pivot can only be detected from the second diagonal of $\Delta$, which implies that the entries of $\Delta$ may change from $r=3$ onwards. Because of that the base of the induction is $r=3$.

 \noindent {\bf r=3:} \\
To prove that the rows of  $\Delta^{3}$ are of type 1-4, we will analyze the effect a change-of-basis pivot marked in $\Delta^{2}$ has on $\Delta^{3}$.
Suppose, without loss of generality, that there is a change-of-basis pivot $\Delta^{2}_{i,i+2}$ on the second diagonal. Consequently, $\Delta^{2}_{i,i+1}$ is a primary pivot marked in $\Delta^{1}$. Recall that  columns in $\Delta^{2}$  are characterized  by  Corollary \ref{cor:charac}.  
Hence, one has the following possibilities:
\begin{enumerate}
\item  Column $i+1$ is of type $2$. Column $i+2$ can be of type $2,3$ or $4$. In each case, $\Delta^{3}_{i,i+2} = 0 $ and all the other entries in column $i+2$ remain the same, i.e, $\Delta^{3}_{s,i+2} = \Delta^{2}_{s,i+2}$ for all $s \neq i$. If column $i+2$ is of type $2$ or $4$ (respectively, $3$) it turns into a column of type $1$ (respectively, $4$).
\item Column $i+1$ is of type $3$. Then $\Delta^{2}_{i,i+1}$ is  a monomial and there is $s< i$ such that $\Delta^{2}_{s,i+1}$ is also a monomial.
\begin{enumerate}
\item  If the column $i+2$ is of type $2$ (resp., type $4$), then  $\Delta^{3}_{i,i+2} = 0$ and $\Delta^{3}_{s,i+2}  = - \Delta^{2}_{s,i+1}(\Delta^{2}_{i,i+1})^{-1} \Delta^{2}_{i,i+2}$, which is a binomial (resp., monomial). Hence, the column $i+2$  remains of  type $2$ (resp., type $4$).
\item  If the column $i+2$ is of type $3$, then change-of-basis pivot $\Delta^{2}_{i,i+2} $ is a monomial and there is $\bar{s} < i$ such that $\Delta^{2}_{\bar{s},i+2} $ is also a monomial. If $s=\bar{s}$, then $\Delta^{3}_{i,i+2} = 0 $ and $\Delta^{3}_{s,i+2} = \Delta^{2}_{s,i+2} - \Delta^{2}_{s,i+1}(\Delta^{2}_{i,i+1})^{-1} \Delta^{2}_{i,i+2} $, which is either a binomial or zero. Hence, column $i+2$ turns into a column of type $2$ or $1$, respectively. On the other hand, if $s\neq\bar{s}$, then $\Delta^{3}_{i,i+2} = 0 $, $\Delta^{3}_{s,i+2} = - \Delta^{2}_{s,i+1}(\Delta^{2}_{i,i+1})^{-1} \Delta^{2}_{i,i+2} $, which is a monomial, and the other entries of column $i+2$ remain the same. Hence, column $i+2$ remains of type $3$.
\end{enumerate}
\item Column $i+1$ is of type $4$. It is analogous to Case $1$.
\end{enumerate}

{\bf Induction hypothesis:} Suppose that the conclusion of the Theorem holds for $3 \leq r < R$. We will show, that it also holds for $r+1$.

Suppose that $\Delta^{r}_{i,j}$ is a change-of-basis pivot in the $r$-th diagonal. Then there is a primary pivot $\Delta^{r}_{i,u}$ in a column  $u < j$. we have one of the possibilities:

\begin{enumerate}
\item Column  $u$ is of type $2$. 
Column $j$ can be of type $2,3$ or $4$. In all cases,  $\Delta^{r+1}_{i,j} = 0 $ and all the other entries in column $r+1$ remain the same. If column $j$ is of type $2$ or $4$ (resp., type $3$) it turns into a column of type 1 (resp., type $4$).
\item Column $u$ is of type $3$. Then $\Delta^{r}_{i,u}$ is  a monomial and there is $s< i$ such that $\Delta^{r}_{s,u}$ is also a monomial.
\begin{enumerate}
\item  If the column $j$ is of type $2$ (resp., type $4$), then the change-of-basis pivot $\Delta^{r}_{i,j} $ is a binomial (resp., monomial). In this case, $\Delta^{r+1}_{i,j} = 0$ and $\Delta^{r+1}_{s,j}  = - \Delta^{r}_{s,u}(\Delta^{r}_{i,u})^{-1} \Delta^{r}_{i,j}$, which is a binomial (resp., monomial). Hence, the column $j$  remains of  type $2$, (resp., type $4$).
\item  If the column $j$ is of type $3$, then the change-of-basis pivot $\Delta^{r}_{i,j} $ is a monomial and there exists $\bar{s} < i$ such that $\Delta^{r}_{\bar{s},j} $ is also a monomial. If $s=\bar{s}$, then $\Delta^{r+1}_{i,j} = 0 $ and $\Delta^{r+1}_{s,j} = \Delta^{r}_{s,j} - \Delta^{r}_{s,u}(\Delta^{r}_{i,u})^{-1} \Delta^{r}_{i,j} $, which is either a binomial or zero. Hence, column $j$ turns into a column of type $2$ or $1$, respectively. On the other hand, if $s\neq\bar{s}$, then $\Delta^{r+1}_{i,j} = 0 $, $\Delta^{r+1}_{s,j} = - \Delta^{r}_{s,u}(\Delta^{r}_{i,u})^{-1} \Delta^{r}_{i,j} $, which is a monomial, and the other entries of column $j$ remain the same. Hence, column $j$ remains of type $3$.
\end{enumerate}
\item  Column $u$ is of type $4$. The only non-zero entry in column $u$ is the primary pivot $\Delta^{r}_{i,u}$ which  is a monomial.  Hence, in column $j$ all the entries remain the same besides the change-of-basis $\Delta^{r+1}_{i,j} =0 $. If $j$ is a column  of type $2$ or $4$ it turns into a column  of type $1$; if $j$ is of type $3$ it turns into a column of type $4$.

\end{enumerate}
\cqd

\begin{teo}[Second Block Characterization]\label{bloco2} Let $\Delta$ be a Novikov differential for which the SSSA is proved correct up to step $R$.  Then we have the following possibilities for a  non-zero row $s \in J_{1}$ of the matrix $\Delta^{r}$ produced by the SSSA  in step $r\leq R$ without realizing  the pre-multiplication by $(T^{r})^{-1}$:
\begin{enumerate}
\item[(A)]  all non null entries are binomials of the form $t^{\ell}- t^{\tilde{\ell}}$, where $\ell,\tilde{\ell} \in \mathbb{Z}$;
\item[(B)] all non null entries are monomials of the form $t^{\ell}$, where $\ell \in \mathbb{Z}$;
\item[(C)] all non null entries are either monomials $t^{\ell}$ or binomials $t^{\ell}- t^{\tilde{\ell}}$. Moreover, if a column $j\in J_{2}$ contains a binomial, then there are no monomials in columns $j'\in J_{2}$ with $j'> j$.
\end{enumerate}
\end{teo}

\dem
We will  prove this theorem by induction  in $r \leq R$.
In the course of the proof, we will also prove the following set of statements: 
\begin{enumerate}
\item[$(i)$] If an entry  $t^{\ell}$ is a primary pivot in row $i$ then  at most one entry will be marked as a change-of-basis pivot in row $i$.
\item[$(ii)$] An entry $t^{\ell} - t^{\tilde{\ell}} $ is never marked as a change-of-basis pivot, i.e. all  change-of-basis pivots are monomials $t^{\ell}$, for some $\ell \in \mathbb{Z}$.  
\item[$(iii)$]  If  $\Delta^{r}_{s,j} = t^{l}- t^{\tilde{l}} $,  $\Delta^{r}_{i,j}$ is a change-of-basis pivot in  row $i > s$ and $\Delta^{r}_{i,u}$  is the primary pivot in  row $i$ with $u<j$, then $\Delta^{r}_{s,u}$ is zero. 
\item[$(iv)$] A primary pivot $\Delta^{r}_{i,u}$ and   a change-of-basis pivot $\Delta^{r}_{i,j}$ in row $i$ are always  monomials  with opposite signs, i.e., $\Delta^{r}_{i,j} = \pm t^{\ell}$ and $\Delta^{r}_{i,u}= \mp t^{\tilde{\ell}}$, for $\ell, \tilde{\ell}\in \mathbb{Z}$.
\item[$(v)$]  A monomial  $\Delta^{r}_{s,u}$  above  a primary pivot $\Delta^{r}_{i,u}$ and  a monomial $\Delta^{r}_{s,j}$ above a change-of-basis pivot $\Delta^{r}_{i,j}$ always  have opposite signs, i.e.,  $\Delta^{r}_{s,j} = \pm t^{l}$ and $\Delta^{r}_{s,u}= \mp t^{\tilde{l}}$, for $l, \tilde{l}\in \mathbb{Z}$.
\end{enumerate}

Observe that the matrices $\Delta^{1}$ and $\Delta^{2}$ differ from the initial matrix $\Delta$ only in the mark-ups of primary and change-of-basis pivots, since the entries can only change as of the $3$-rd step of the SSSA. 

 \noindent{\bf Base case r=3:} 
 In order to prove that the rows of $\Delta^{3}$ satisfies conditions $(A), (B)$ and $(C)$ of the theorem, we must  analyze the effect on a row of $\Delta^{3}$ caused by  a change-of-basis pivot marked in the second step $r=2$ of the SSSA. 

\begin{figure}[!ht]
\centering
\includegraphics[scale=1]{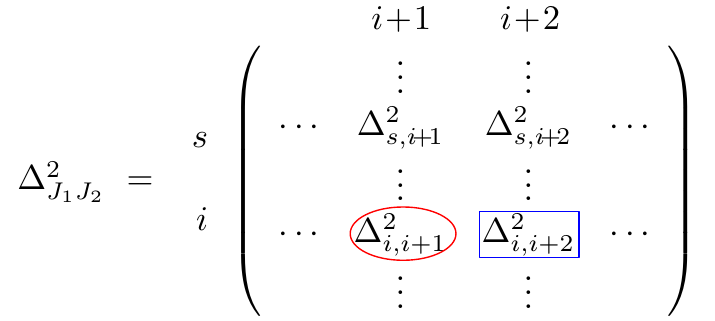}
\vspace{-0.4cm}
\caption{Primary  and  change-of-basis pivots in the first and second diagonals of $\Delta^{2}$, respectively.}\label{fig:2dimmatrix}
\end{figure}

\begin{figure}[!ht]
\centering
\includegraphics[scale=1]{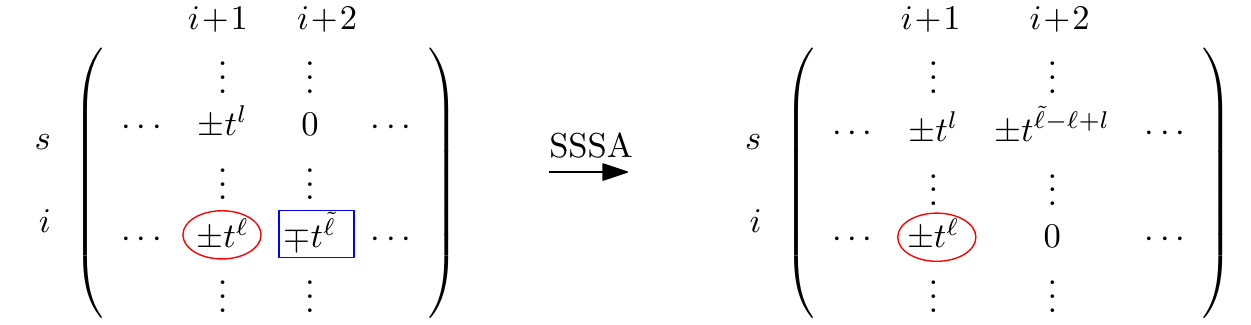}
\vspace{-0.4cm}
\caption{$\Delta^{2}_{J_{1}J_{2}}$ and $\Delta^{3}_{J_{1}J_{2}}$, respectively.  } \label{fig:dem-delta21}
\end{figure}

\begin{figure}[!ht]
\centering
\includegraphics[scale=1]{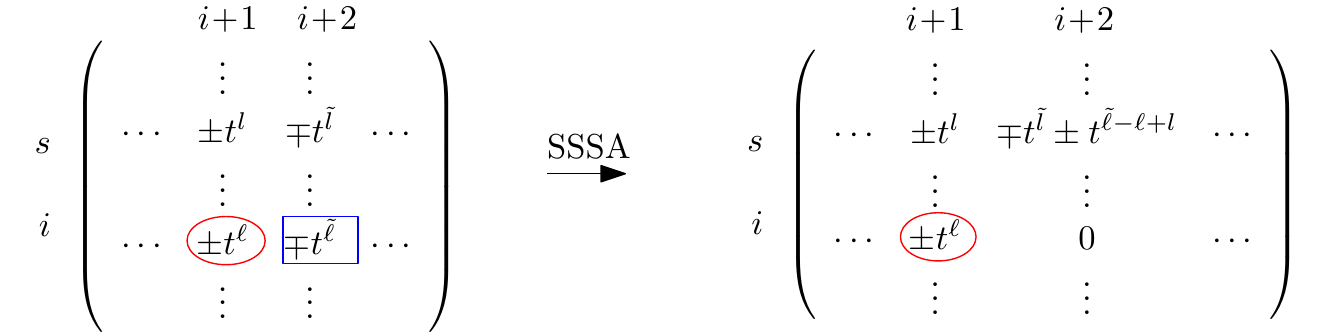}
\vspace{-0.4cm}
\caption{$\Delta^{2}_{J_{1}J_{2}}$ and $\Delta^{3}_{J_{1}J_{2}}$, respectively. } \label{fig:dem-delta22}
\end{figure}

Suppose, without loss of generality, that there is a change-of-basis pivot $\Delta^{2}_{i,i+2}$ on the second diagonal. Consequently, $\Delta^{2}_{i,i+1}$ is a primary pivot marked in the first step of the SSSA, see Figure \ref{fig:2dimmatrix}. Recall that the rows in $\Delta^{2}$  are characterized by  Corollary \ref{cor:charac}. A change-of-basis pivot only occurs in row $i$  and column $i+2$  if this row is of type $3$. In what follows,  we analyze  the effect of  this  change of basis on a row $s$  with $s<i$:
\begin{enumerate}
\item  If  row $s$ is null (i.e., of type $1$) then only row $i$ is altered  and becomes a row of type $B$.
\item Suppose that row $s$ is of type $2$. If the only non-zero entry in row $s$ is in a column different from $i+1$, then row  $s$ remains unaltered and row  $i$ turns  into a row of type  $B$. On the other hand, if this non zero entry is in the column  $i+1$, which is the same column as that of the primary pivot, then row $i$ turns into a row of type $B$ and row $s$ turns into a row of type $A$, as one can see in Figure \ref{fig:dem-delta2}.

\begin{figure}[!ht]
\centering
\includegraphics[scale=1]{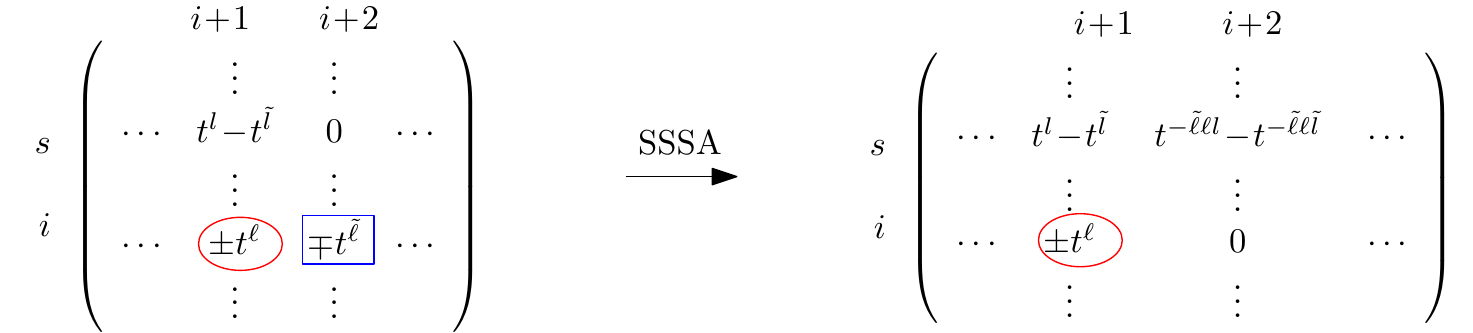}
\vspace{-0.4cm}
\caption{$\Delta^{2}_{J_{1}J_{2}}$ and $\Delta^{3}_{J_{1}J_{2}}$, respectively.} \label{fig:dem-delta2}
\end{figure}

\item Suppose that row $s$ is of type $3$. If $\Delta_{s,i+1}$ is zero, then row $s$ remains unaltered. If $\Delta^{2}_{{s,i+1}}=t^{l}$, one has two possibilities for $\Delta^{2}_{s,i+2}$, namely, $0$ or $t^{\tilde{l}}$.  In the first case, after performing the change of basis,  row $s$ turns into a row of type $B$ (see Figure \ref{fig:dem-delta21}),  and in the second case it turns into a row of type $C$ (see Figure \ref{fig:dem-delta22}). 
\end{enumerate}


In the base case, it is ease to see that $(i)$ through $(v)$ hold.

In order to prove $(i)$ note that, the only case that needs to be  analyzed  is when  $\Delta^{2}_{i,i+1}$ is a primary pivot and  $\Delta^{2}_{i,i+2}$ is a change-of-basis pivot. Observe that $\Delta^{2}_{i,i+3} = 0$, since rows of $\Delta^{2}$ have at most two non zero entries.  As pivots do not generate entries, by Lemma \ref{creation},  the entry $\Delta^{2}_{i,i+3}$ is not altered by change-of-basis pivots marked in step $2$. Hence, $\Delta^{3}_{i,i+3} = \Delta^{2}_{i,i+3}= 0$ and it is not a change-of-basis pivot.

In order to prove $(ii)$,  we must consider each row $s$ where the entry $\Delta^{3}_{s,s+3}=t^{l} -t^{\tilde{l}}$ was generated in $\Delta^{3}$, otherwise this entry would be in a row of type $2$ and hence,  could not be a change-of-basis pivot. There are exactly two  ways that this entry can be generated:  when row $s$ was of type $2$ and of type $3$ in $\Delta^{2}$.  Suppose by contradiction that $\Delta^{3}_{s,s+3}$ is a change-of-basis pivot. Since primary pivots do not generate entries  and there exists at most one change-of-basis pivots in a row then  row $s$ could not be of type $2$ in $\Delta^{2}$, see Figure \ref{fig:generalmatrix2}. If row $s$ was of type $3$ in $\Delta^{2}$ then the primary pivot would generate $\Delta^{3}_{s,s+3}$ which also contradicts   Lemma \ref{creation}, see Figure \ref{fig:generalmatrix3}.  Consequently,  $ \Delta^{3}_{s,s+3}$ is not a change-of-basis pivot.

In order to prove $(iii)$,
 suppose by contradiction that $\Delta^{3}_{s,u} \neq 0$. Since  rows of the initial matrix $\Delta$ do not admit two non zero entries with one of them being a binomial, see Corollary \ref{cor:charac}, then  $\Delta^{2}_{s,u} = \Delta^{3}_{s,u}$ generates $\Delta^{3}_{s,i+3}$ in $\Delta^{3}$. By Definition \ref{def:creation}, there must exist a change-of-basis pivot in the second diagonal  in a row $\tilde{i} >i$ and a primary pivot in row $\tilde{i}$ and column $u$, which contradicts the fact that  each column has  at most one primary pivot.

Items $(iv)$ and $(v)$ are trivially true.
 
\vspace{0.3cm}

 {\bf Induction hypothesis:} Suppose that Theorem \ref{bloco2}  and item $(i)$ through $(v)$ hold for $\xi \leq r < R$. We will show that they also hold for $r+1$. 
First, note that by the induction hypothesis, at most one entry in a fixed row in $J_{1}$ is marked as  change-of-basis pivot up to  step $r$. By Corollary \ref{uniqueseq}, given an entry $\Delta_{s,u}$, the $\Delta_{s,u}$-lineage is formed by a unique generation sequence until  $\Delta^{r+1}$, which is either a binomial lineage, if $\Delta_{s,u}$ is a binomial; or  a monomial lineage or an eventual binomial lineage, if $\Delta_{s,u}$ is a monomial{\footnote{Note that in this proof, whenever we assume the induction hypothesis for $r$,   the index $r$  is  shifted by one, i.e. $r+1$, when referring  to the lineages.}}. More specifically, if row $s$ is of type $2$ in $\Delta$, where $\Delta_{s,u}\neq 0$ is a binomial, then the $\Delta_{s,u}$-lineage is a binomial lineage. If row $s$ is of type $3$, where $\Delta_{s,u}$ and $\Delta_{s,j}$ are monomials, then each monomial determines a lineage, which are either both monomial lineages or one monomial lineage which ceases and merges with the other to create an eventual binomial lineage. If row $s$ is of type 4 in $\Delta$, where $\Delta_{s,u}$ is a monomial, then $\Delta_{s,u}$-lineage is monomial lineage.
 It is important to keep in mind that, if there is at most one change-of basis per row  up to step $r$, then the  $\Delta_{s,u}$-lineage is formed by a unique  generation sequence until $\Delta^{r+1}$. This follows since entries in $\Delta^{r+1}$ can only be generated  by change of basis determined in step $r$.

By the induction hypothesis, it follows that if  $\Delta^{r}_{i,j}$ is a change-of-basis pivot on the $r$-th diagonal and  $\Delta^{r}_{i,u}$ is the primary pivot of row $i$, then these entries must be monomials. Moreover,
  $\Delta^{r}_{i,j} = \pm t^{\ell}$ and $\Delta^{r}_{i,u}= \mp t^{\tilde{\ell}}$, for $\ell, \tilde{\ell}\in \mathbb{Z}$.

\vspace{0.3cm}

Now we will prove that the statement of the Theorem \ref{bloco2} and item $(i)$ through $(v)$ hold  for $r+1$.

$\bullet$ We will first show that the non-zero rows in $\Delta^{r+1}$ are of type $A,B$ or $C$.

In order to prove this, we will perform all the possible  changes of basis that could occur due to a change-of-basis pivot in the $r$-th diagonal. Let  the monomial $\Delta^{r}_{i,u}$ be the primary pivot in row $i$  and  the monomial $\Delta^{r}_{i,j}$ be a change-of-basis pivot in diagonal $r$.  Observe that after the change of basis, row $i$ will remain of the same type. Each row $s < i$ in $\Delta^{r}$ is of type $A,B$ or $C$, and thus one has the following cases to analyze the effect a change of basis causes in  row $s$ (bear in mind the configuration of the matrix in Figure \ref{fig:generalmatrix}):

\begin{enumerate}
\item If $\Delta^{r}_{s,u} =0$, then row $s$ remains unaltered after performing the change of basis.
\item If $\Delta^{r}_{s,u} \neq 0$ and $\Delta^{r}_{s,j}=0$, then $\Delta^{r+1}_{s,j}= - \Delta^{r}_{i,j} (\Delta^{r}_{i,u})^{-1} \Delta^{r}_{s,u}$, which is a monomial if $\Delta^{r}_{s,u}$ is a monomial, or a binomial if $\Delta^{r}_{s,u}$ is a binomial. Hence, row $s$ remains of the same type.
\item If $\Delta^{r}_{s,u} \neq 0$ and $\Delta^{r}_{s,j}\neq 0$, by the induction hypothesis the only possible case is $\Delta^{r}_{s,u} = \pm t^{l}$ and $\Delta^{r}_{s,j}= \mp t^{\tilde{l}}$. Then $\Delta^{r+1}_{s,j}= \Delta^{r}_{s,j} - \Delta^{r}_{i,j} (\Delta^{r}_{i,u})^{-1} \Delta^{r}_{s,u}$, which  ensures that $\Delta^{r+1}_{s,j}$ is a zero entry or a binomial with coefficients equal to $\pm 1$.  Hence,  row $s$ turns into a row of type $B$ or $C$, respectively. 

Hence,   every row $s \in J_{1}$ of $\Delta^{r+1}$ is also of type $A$, $B$ or  $C$.
\end{enumerate}

$\bullet$ We will now show that  item $ (i)$ holds for $\Delta^{r+1}$.

Let  $\Delta^{r+1}_{i,j }$ be a change-of-basis pivot marked in  diagonal $r+1$.  Suppose by contradiction that an entry $\Delta^{\xi}_{i,t}$,  where $t< j$, was marked as a change-of-basis pivot in an earlier step $\xi < r+1 $. Consequently,  the  primary pivot in row $i$, which is in a column $u < t$, was marked in a previous step $ < \xi$. By   item $(ii)$ of the induction hypothesis, 
this  primary pivot is a monomial  and it can not generate entries, which implies that $\Delta^{\xi}_{i,t }$ and $\Delta^{r+1}_{i,j}$ are not in the lineage of this primary pivot.  Therefore,   
$\Delta^{\xi}_{i,t }$ and $\Delta^{r+1}_{i,j}$ must be in the same lineage,  since these entries were generated up to  step $r+1$ by the change-of-basis pivots up to  step $r$, 
and by the induction hypothesis, there is only one change-of-basis pivots  per row up to  step $r$. This is a contradiction, since by Lemma \ref{creation} the change-of-basis pivot $\Delta^{\xi}_{i,t }$ can not generate entries.

$\bullet$  We will prove  item $(ii)$ for $r+1$:  

Let $\Delta^{r+1}_{i,j}$ be a binomial in  diagonal $r+1$.  Suppose by contradiction that this entry  is a change-of-basis pivot marked in step $r+1$. Let $u$ be the column  of the primary pivot in row $i$, hence $u < j$.

If  row $i$  was originally of type $2$ in $\Delta^{2}$, then there is only one lineage in row $i$  until $\Delta^{r+1}$, which is a binomial lineage. Hence, $\Delta^{r+1}_{i,j}$ is a successor of the primary pivot in row $i$, which must have generated an entry,  contradicting Lemma \ref{creation}.

If row $i$ was of type $3$ in $\Delta^{2}$, then originally there were two lineages that merged in order to create a binomial. Note that the primary pivot in  row $i$ can not be a monomial, i.e., it can not be marked before the two sequences have merged, since pivots do not generate entries. Hence the primary pivot must be a binomial. As we have seen in the previous paragraph, this contradicts Lemma \ref{creation}.  

Note that row $i$ could not be originally of type $4$ in $\Delta^{2}$, since  this row would not contain binomials.

$\bullet$ We will prove  item $(iii)$ for $r+1$:

Let $\Delta^{r+1}_{i,j}$ be a change-of-basis pivot marked in  diagonal $r+1$, $\Delta^{r+1}_{s,j}$ be a binomial and let $\Delta^{r+1}_{i,u}$ be the primary pivot  in row $i$. Suppose by contradiction that the entry $\Delta^{r+1}_{s,u}$ is non-zero.

 If row $s$ is of type $2$ in $\Delta^{2}$, then $\Delta^{r+1}_{s,u}$ and $\Delta^{r+1}_{s,j}$ are in the same lineage, i.e., $\Delta^{r+1}_{s,j}$ must be a successor of $\Delta^{r+1}_{s,u}$, which is a contradiction. In fact, $\Delta^{r+1}_{s,j}$ is not an immediate successor of $\Delta^{r+1}_{s,u}$, since in this case it would imply the existence of two primary pivots in column $u$. Moreover, $\Delta^{r+1}_{s,j}$ is not an eventual successor of $\Delta^{r+1}_{s,u}$, since it would imply the existence of two change-of-basis pivots in row $i$ marked up to step $r$.

If row $s$ is of type $3$ in $\Delta^{2}$ and if $\Delta^{r+1}_{s,u}$ is a binomial  then the argument is the same as the one above. However, if $\Delta^{r+1}_{s,u}$ is a monomial, one has two cases to consider: $\Delta^{r+1}_{s,u}$ and $\Delta^{r+1}_{s,j}$ are in the same lineage or in different lineages.  In both cases there exists a $\xi<r+1$ such that the entry $\Delta^{\xi}_{s,u}$  generated an entry in $\Delta^{\xi+1}$. Indeed, if they are in the same lineage, which is an eventual binomial lineage, then the $(s,j)$-entry is a successor of the $(s,u)$-entry. Now, if they are in different lineages, the lineage containing  the $(s,u)$-entry
ceases due to the appearing of the first binomial in a step $\xi^{\ast}<r +1$, which appears in the entry corresponding to the last element of this lineage. Hence, the last element in this lineage can not be $\Delta^{\xi^{\ast}-1}_{s,u}$ since the entry $\Delta^{r+1}_{s,u}$ is still a monomial. Therefore there exists a $\xi<\xi^{\ast}<r +1$ such that $\Delta^{\xi}_{s,u}$  generates an entry.

Note that  row $s$ could not be originally of  type $4$ in $\Delta^{2}$, 
since  this row would not contain binomials.

$\bullet$ 
We will prove item for $(iv)$ for $r+1$, i.e. that 
a primary pivot $\Delta^{r+1}_{i,u}$ and   a change-of-basis pivot $\Delta^{r+1}_{i,j}$ in row $i$ are always  monomials  with opposite signs.

As consequence of the induction hypothesis $(iv)$ and $(v)$ for $r$,
 all elements in a monomial lineage have the same sign, up to step $\xi \leq r+1$. Moreover, if row $i$ was of type 3 in $\Delta$, then the two monomial lineages of this row have opposite signs.
 Now, suppose that there is a change-of-basis pivot $\Delta^{r +1}_{i,j}$ in  diagonal $r+1$ and let $\Delta^{r +1}_{i,u}$ be the primary pivot in row $i$ with $u<j$. Since pivots do not generate entries by Lemma \ref{creation},  the entries $\Delta^{r +1}_{i,j}$ and $\Delta^{r +1}_{i,u}$ are clearly in different lineages, therefore  they have opposite signs.

$\bullet$ We will prove  item $(v)$ for $r+1$, i.e, that a monomial  $\Delta^{r+1}_{s,u}$  above  a primary pivot $\Delta^{r+1}_{i,u}$ and  a monomial $\Delta^{r+1}_{s,j}$ above a change-of-basis pivot $\Delta^{r+1}_{i,j}$ always  have opposite signs.

As a consequence of the induction hypothesis $(iv)$ and $(v)$,
 until $\Delta^{r+1}$, all elements in a monomial lineage have the same coefficient, which is either $+1$ or $-1$. Moreover, if row $i$ was of type 2 in $\Delta$, then the two monomial lineages of this row have opposite signs.
Let $\Delta^{r+1}_{s,u}$ be a monomial    above  a primary pivot $\Delta^{r+1}_{i,u}$ and $\Delta^{r+1}_{s,j}$ be a monomial  above a change-of-basis pivot $\Delta^{r+1}_{i,j}$. Observe that, by the induction hypothesis,  $\Delta^{r+1}_{s,u}$ and $\Delta^{r+1}_{s,j}$ are not in the same lineage.  Hence, they have opposite signs. 

\cqd

We now proceed with the proof Theorem \ref{teo:primarypivots}.

\noindent {\bf Proof of Theorem \ref{teo:primarypivots}:}
Let $\Delta$ be a Novikov differential. By the characterization of the initial matrix, see  Corollary \ref{cor:charac}, the entries of $\Delta$ are invertible in $\mathbb{Z}((t))$; hence, one can apply the SSSA in $\Delta$. Since the first change of basis can only occur from step $2$ to step $3$, then the SSSA is correct until step $2$, and the entries of  $\Delta^{1}$ are equal to the entries of  $  \Delta^{2}$. Now, using Lemma \ref{Tinverso}, one can apply the SSSA to each block of $\Delta$. Theorem \ref{bloco1} and \ref{bloco2} imply that the pivots in $\Delta^{3}$ are invertible.  Hence, the SSSA is also correct for $\Delta^{3}$. By an induction argument, one can suppose  that the SSSA is proved correct until step $r$. Theorems \ref{bloco1} and \ref{bloco2} also imply that the SSSA
 is proved correct for $\Delta^{r+1}$. Therefore, Theorem \ref{teo:primarypivots} follows.
\cqd

Observe that, if $\Delta^{L}$ is the last matrix produced by the SSSA, then  the non null columns of $\Delta^{L}$ are the columns containing primary pivots. The primary pivots are non-zero and are unalterable after being identified. Moreover, $\Delta^{L}\circ \Delta^{L}=0$.

The next  example  shows that during the application of the SSSA infinite series can appear as entries of an intermediate matrix due to multiplication by $(T^{r})^{-1}$.  However, the last matrix produced by the SSSA   does not contain entries which are infinite series.

\begin{ex}\label{exemplo2} 
Applying the SSSA to the Novikov matrix $\Delta$ in Figure \ref{fig:delta0ex2}, one  obtains the sequence of Novikov matrices $\Delta^{1},\cdots,\Delta^{6}$ presented in Figures \ref{fig:delta1ex}, $\cdots$,  \ref{fig:delta6ex}, respectively.  Observe that the 
entry $\Delta^{3}_{3,7}$ is an infinite Laurent series.

\end{ex}

The next two results imply that the last matrix $\Delta^{L}$ produced by the SSSA  never contains an entry which is an infinite series, rather all entries are Laurent polynomials in $\mathbb{Z}((t))$.

The proof of  Theorem \ref{lemalinhanula} follows the same steps of  the proof of its analogous  version  in \cite{MdRS}, where the SSSA is done over a field $\mathbb{F}$.

\begin{teo}\label{lemalinhanula}
Given a Novikov complex $(\mathcal{N}_{\ast}(f),\Delta)$ on surfaces, let $\Delta^{L}$ be the last matrix produced by the SSSA over $\mathbb{Z}((t))$.
If column $j$ of $\Delta^{L}$ is non null then  row $j$ is null.
\end{teo}

\dem
The statement of the lemma is equivalent to say that $\Delta^{L}_{j \bullet}\Delta^{L}_{\bullet j} = 0$ for all $j$. If $\Delta^{L}_{\bullet j }=0 $ then it is trivial that $\Delta^{L}_{j \bullet}\Delta^{L}_{\bullet j} = 0$. Suppose that $\Delta^{L}_{\bullet j } \neq 0 $. Let $s$ be an integer such that $j \in J_{s}$. Labelling the  primary pivots in block $J_{s}$  such that, if $\Delta^{L}_{i_{1},j_{1}}, \cdots, \Delta^{L}_{i_{a},j_{a}}$ are the primary pivots in block $J_{s}$, then $i_{1}<i_{2}<\cdots <i_{a}$, one has that $j_{1}, \cdots, j_{a}$ are the non null columns of $J_{s}$. Moreover, $\Delta^{L}_{i_{a},j_{a}}$ is the  unique non zero entry in row $i_{a}$. Row $i_{a-1}$ has non zero entry in column $j_{a-1}$ and may have another one non zero entry in column $j_{a}$, and so on. 
Since $\Delta^{L}\circ \Delta^{L}=0$, one has 
$$ 0 = \Delta^{L}_{i_{a}\bullet} \Delta^{L}_{\bullet j'} = \Delta^{L}_{i_{a}j_{a}} \Delta^{L}_{j_{a}j'}, $$
for all $j'$. Since $\Delta^{L}_{i_{a}j_{a}}$ is a primary pivot, hence non null, then $\Delta^{L}_{j_{a}j'}=0$ for all $j'$, i.e., $\Delta^{L}_{j_{a}\bullet }=0$.
Analogously, one has
$$ 0 = \Delta^{L}_{i_{a-1}\bullet} \Delta^{L}_{\bullet j'} = \Delta^{L}_{i_{a-1}j_{a-1}} \Delta^{L}_{j_{a-1}j'}+ \Delta^{L}_{i_{a-1}j_{a}} \Delta^{L}_{j_{a}j'}, $$
for all $j'$. Since $\Delta^{L}_{i_{a-1}j_{a-1}}\neq 0$ and $\Delta^{L}_{j_{a}j'}=0$, it follows that $\Delta^{L}_{j_{a-1}j'}=0$ for all $j'$, i.e., $\Delta^{L}_{j_{a-1}\bullet }=0$. Proceeding in this way, one can show  the nullity of rows $j_{a-2}, \cdots, j_{1}$.
\cqd

\begin{cor}
Given a Novikov complex $(\mathcal{N}_{\ast}(f),\Delta)$ on a surface, let $\Delta^{L}$ be the last matrix produced by the SSSA over $\mathbb{Z}((t))$. Then, the entries of $\Delta^{L}$ are monomials $t^{\ell}$ or binomials $t^{\ell}-t^{\tilde{\ell}}$.
\end{cor}

\dem
By Theorems  \ref{bloco1} and \ref{bloco2}, without performing the pre-multiplication by $(T^{r})^{-1}$, the entries of $\Delta^{r}$ are  monomials $t^{\ell}$ or binomials $t^{\ell_{1}}-t^{\ell_{2}}$. Moreover, the pre-multiplication  by $(T^{r})^{-1}$ only affects  row $j$
if column $j$ contains a primary pivot. However, by Lemma \ref{lemalinhanula}, these rows will be zeroed  out by the time SSSA reaches the last matrix $\Delta^{L}$.
\cqd

\section{Spectral Sequence $(E^{r},d^{r})$ associated to SSSA}\label{ultimo}

Let $M$ be a smooth closed orientable 2-dimensional manifold, $f:M\to S^{1}$ be a circle-valued Morse function.
Let  $(\mathcal{N}_{\ast}(f),\Delta)$ be a  Novikov chain complex with the finest filtration ${F}$ and $(E^{r},d^{r})$ be  the   associated spectral sequence, as in Section  \ref{spectralsequence}.  In this section, we show that the SSSA provides a mechanism to recover the modules $E^{r}$ and differentials $d^{r}$ of the spectral sequence. More specifically, the SSSA provides a system to detect the generators of  $E^{r}$ in terms of the original basis of $\mathcal{N}_{\ast}(f)$ and to identify the differentials $d^{r}$ with the primary pivots in the  $r$-th diagonal.

For instance, consider the Novikov complex  $(\mathcal{N}_{\ast}(f),\Delta)$ presented in Example \ref{exemplo}. 
The matrices produced by the SSSA are illustrated in Figures 
\ref{fig:delta1ex}, $\cdots$,  \ref{fig:delta6ex}. 
The Novikov homology of the complex    $(\mathcal{N}_{\ast}(f),\Delta)$ is given by
$$ H^{Nov}_{0}(M,f) = 0 , \ \ H^{Nov}_{1}(M,f) = 0 , \ \ H^{Nov}_{2}(M,f) = 0   .$$

\begin{figure}[!ht]
\centering
\includegraphics[scale=0.99]{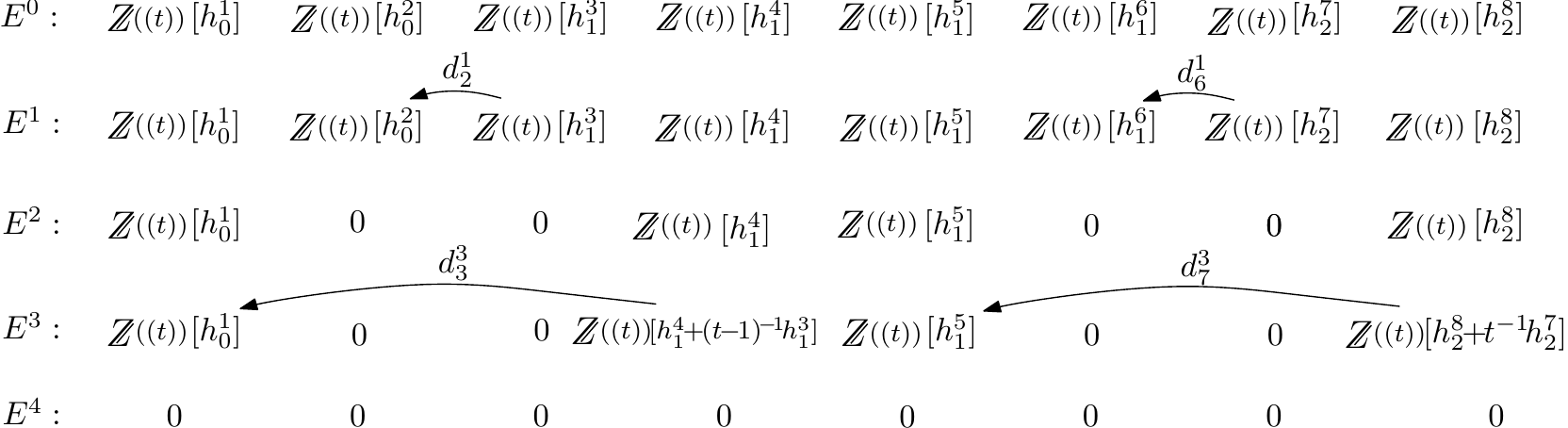}
\caption{A spectral sequence associated to the Novikov complex $(\mathcal{N}_{\ast},\Delta)$  presented in Example \ref{exemplo}.  }\label{fig:nov_complex_ss2}
\end{figure}

Consider the filtration ${F}$  on $(\mathcal{N}_{\ast}(f),\Delta)$ defined by 
$$ F_{p}\mathcal{N}_{k} = \displaystyle\bigoplus_{h^{\ell}_{k}, \ \ell \leq p+1} \mathbb{Z}((t)) \langle h^{\ell}_{k} \rangle . $$
The spectral sequence associated to  $(\mathcal{N}_{\ast}(f),\Delta)$  endowed with this filtration $F$ is presented in Figure \ref{fig:nov_complex_ss2}.
Note that: 
\begin{enumerate}
\item[(1)] Each $\mathbb{Z}((t))$-module $E^{r}_{p}$ is generated by a $k$-chain $\sigma^{p+1,r}_{k}$ determined in the $r$-th step of the SSSA.
\item[(2)] The differentials $d^{1}_{2}$, $d^{1}_{6}$, $d^{3}_{3}$ and $d^{3}_{7}$ are isomorphisms which could be interpreted as the multiplication by the invertible polynomials  $\Delta^{1}_{2,3}$, $\Delta^{1}_{6,7}$, 
$\Delta^{3}_{1,4}$ and $\Delta^{3}_{5,8}$, respectively. 
\item[(3)]  The spectral sequence $(E^{r},d^{r})$ converges to the Novikov homology of $(\mathcal{N}_{\ast}(f),\Delta)$.
\end{enumerate}

In the $2$-dimensional setting, we will prove that the remarks in $(1)$ and $(2)$ hold for any Novikov chain complex  endowed with a finest filtration.
The first step in this direction is the next proposition, which  establishes a formula for the modules $Z^{r}_{p,k-p}$ via the chains $\sigma_{k}^{i,j}$'s determined by the SSSA  applied to $\Delta$.

\begin{prop}\label{prop:formuladoz} Let
$\kappa$ be the first column in $\Delta$ associated to a
$k$-chain and consider  $\mu^{j,\zeta}=0$ whenever the primary pivot of  column $j$ is below row $(p-r+1)$ and $\mu^{j,\zeta}=1$
otherwise. Then 
$$Z^{r}_{p}=
\mathbb{Z}((t))[\mu^{p+1,{r}}\sigma_{k}^{p+1,{r}},\mu^{p,{r-1}}\sigma_{k}^{p,{r-1}},\ldots
,\mu^{\kappa,{r-p-1+\kappa}}\sigma_{k}^{\kappa,{r-p-1+\kappa}}].$$
\end{prop}
\dem By definition, $\sigma_{k}^{p+1-\xi,{r-\xi}}$ is associated to column $(p+1-\xi)$ of the matrix $\Delta^{r-\xi}$, for $\xi\in\{0,\ldots, p+1-\kappa\}$, and $\mu^{p+1-\xi,{r-\xi}}=1$ if and only if the primary
pivot on column $(p+1-\xi)$ is in or above  row
$(p+1-\xi)-(r-\xi)=p-r+1$ or if this column does not have a primary pivot. If $\sigma_{k}^{p+1-\xi,{r-\xi}}$ is such that  $\mu^{p+1-\xi,{r-\xi}}=1$, we show that $\sigma_{k}^{p+1-\xi,{r-\xi}}$ is a
$k$-chain which corresponds to a generator of $Z^r_p$. In fact,  $\sigma_{k}^{p+1-\xi,{r-\xi}}$ is in $F_p\mathcal{N}_k$ for
$\xi\geq 0$. Furthermore, all nonzero entries
of  column $(p+1-\xi)$ of $\Delta^{r-\xi}$ are in or above row
$(p-r+1)$, since step $(r-\xi)$ of SSSA has zeroed out all entries below 
diagonal $(r-\xi)$. Hence, the boundary of
$\sigma_{k}^{p+1-\xi,{r-\xi}}$ is in $F_{p-r}\mathcal{N}_{k-1}$.

On the other hand, any element in $Z_{p}^r$ is a linear
combination over $\mathbb{Z}((t))$ of
$\mu^{p+1-\xi,{r-\xi}}\sigma_{k}^{p+1-\xi,{r-\xi}}$ for
$\xi=0,\ldots ,p+1-\kappa$, as we prove below  by multiple induction in
$p$ and $r$.

\noindent\begin{bf}Base case:\end{bf}
\begin{itemize}
\item Denote by $\kappa$  the first column of $\Delta$
associated to a $k$-chain and  let  $\xi$ be such that the boundary of
$h_k^{\kappa}$ is in $F_{\kappa-1-\xi}\mathcal{N}_k$ but it is not in
$F_{\kappa-1-\xi-1}\mathcal{N}_k$. We will show that $Z_{\kappa-1}^{r}=\mathbb{Z}((t))[\mu^{\kappa,{r}}\sigma_{k}^{\kappa,{r}}]$.

 Since $Z_{\kappa-1}^{r}$ is generated by  $k$-chain
 in $F_{\kappa-1}\mathcal{N}_k$ with boundaries in $F_{\kappa-1-r}\mathcal{N}_{k-1}$ and
  there is only one chain $h_k^{\kappa}$ in $F_{\kappa-1}\mathcal{N}_k$ then:
  \begin{enumerate}
  \item[(a)] If $\xi <r$ then $\partial h_k^{\kappa}\notin
  F_{\kappa-1-r}\mathcal{N}_{k-1}$. Thus, $Z_{\kappa-1}^{r}=0$.
  \item[(b)] If $\xi \geq r$ than $\partial h_k^{\kappa}\in F_{\kappa-1-r}\mathcal{N}_{k-1}$. Thus,
  $Z_{\kappa-1}^{r}=\mathbb{Z}((t))[ h_k^{\kappa}]$.
  \end{enumerate}
 On the other hand, since there
 is no change of basis caused by the SSSA that affects the first column of
 $\Delta_k$, $\sigma_{k}^{\kappa,{r}}=h_k^{\kappa}$, where $\sigma_{k}^{\kappa,{r}}$ is a $k$-chain
 associated to the column $\kappa$  of $\Delta^{r}$. Furthermore,
 $\mu^{\kappa,{r}}=1$ if and only if the boundary of
 $h_k^{\kappa}=\sigma_{k}^{\kappa,{r}}$ is in or above the $r$-th diagonal.
 Hence
  \begin{enumerate}
  \item[(a)] If $\xi <r$ then $\mu^{\kappa,{r}}=0$. Thus  $\mathbb{Z}((t))[\mu^{\kappa,{r}}\sigma_{k}^{\kappa,{r}}]=0$
  \item[(b)] If $\xi\geq r$ then $\mu^{\kappa,{r}}=1$. Thus
  $\mathbb{Z}((t))[\mu^{\kappa,{r}}\sigma_{k}^{\kappa,{r}}]=\mathbb{Z}((t))[\sigma_{k}^{\kappa,{r}}]=\mathbb{Z}((t))[
  h_k^{\kappa}]$.
 \end{enumerate}
It follows that in both cases
$Z_{\kappa-1}^{r}=\mathbb{Z}((t))[\mu^{\kappa,{r}}\sigma_{k}^{\kappa,{r}}]$.

\item Denote by $\xi_1$ 
the first diagonal in $\Delta$ which contains a nonzero entry in $\Delta_k$. Note that the nonzero entries of the columns of
$\Delta$ corresponding to the chains $h_{k}^{\kappa}, \ldots, h_{k}^{p+1}
$ are in or  above  the row $(p-\xi_1+1)$. We will show that $Z_{p}^{\xi_1}=\mathbb{Z}((t))[\mu^{p+1,{\xi_1}}\sigma_{k}^{p+1,{r}},\ldots,\mu^{\kappa,{\kappa-p+1+\xi_1}}
\sigma_{k}^{\kappa,{\kappa-p+1+\xi_1}}]$.

By the definition of  $Z_{p}^{\xi_1}$ and the remark above we have that 
$Z_{p}^{\xi_1}=\mathbb{Z}((t))[h_{k}^{p+1},\ldots , h_{k}^{\kappa}].$
On the other hand, it is easy to see that  $\sigma_{k}^{j,{\xi_1}}=h_k^{j}$, $j=\kappa,\ldots p+1$ and $\mu ^{j,{\xi_1}}=1$, $j=\kappa,\ldots p+1$. Hence,
$$\mathbb{Z}((t))[\mu^{p+1,{\xi_1}}\sigma_{k}^{p+1,{r}},\ldots,\mu^{\kappa,{\kappa-p+1+\xi_1}}
\sigma_{k}^{\kappa,{\kappa-p+1+\xi_1}}]=\mathbb{Z}((t))[h_{k}^{p+1},\ldots ,
h_{k}^{\kappa}].$$
\end{itemize}

\noindent\begin{bf}Induction Hypothesis: \end{bf} Suppose
$$Z^{r-1}_{p-1}=\mathbb{Z}((t))[\mu^{p,{r-1}}\sigma_{k}^{p,{r-1}},\ldots
,\mu^{\kappa,{r-p-1+\kappa}}\sigma_{k}^{\kappa,{r-p-1+\kappa}}].$$

Note that if $\mu^{p+1,r}=0$ then $Z^{r}_p=Z^{r-1}_{p-1}$ and result in this case follows by the induction hypothesis. 
Suppose that  $\mu^{p+1,r}=1$ and let
$\mathfrak{h}_k=b^{p+1}h_k^{p+1}+\cdots +b^{\kappa}h_k^{\kappa}\in Z_{p,k-p}^r$.  If $b^{p+1}=0$ then $\mathfrak{h}_k\in
Z^{r-1}_{p-1}$ and the result follows. If $b^{p+1}\neq 0$, then 
$$\mathfrak{h}_k=b^{p+1}\sigma^{p+1,r}_{k}+(b^{p}-b^{p+1} c_{p}^{p+1,r})h_{k}^{p}+\cdots
+(b^{\kappa}-b^{p+1}c_{\kappa}^{p+1,r})h_{k}^{\kappa}.$$
and hence
$\mathfrak{h}_k-b^{p+1}\sigma^{p+1,r}_{k}\in F_{p-1}\mathcal{N}_k$. Since
$\mathfrak{h}_k$ and $\sigma^{p+1,r}_{k}$ have their boundaries on and
above row $p-r+1$ then the boundary of
$\mathfrak{h}_k-b^{p+1}\sigma^{p+1,r}_{k}$ is on and above 
row $p-r+1$. Hence
$\mathfrak{h}_k-b^{p+1}\sigma^{p+1,r}_{k}\in Z^{r-1}_{p-1}$. By
the induction hypotheses, we have that
$\mathfrak{h}_k-b^{p+1}\sigma^{p+1,r}_{k}=\alpha_p\mu^{p,r-1}\sigma_{k}^{p,{r-1}}+\cdots
+
\alpha_{\kappa}\mu^{\kappa,{r-p-1+\kappa}}\sigma_{k}^{\kappa,{r-p-1+\kappa}}$,
i.e.
$$\mathfrak{h}_k=b^{p+1}\sigma^{p+1,r}_{k}+\alpha_p\mu^{p,r-1}\sigma_{k}^{p,{r-1}}+\cdots
+
\alpha_{\kappa}\mu^{\kappa,{r-p-1+\kappa}}\sigma_{k}^{\kappa,{r-p-1+\kappa}}.$$
\cqd

\begin{lem}\label{prop:formulabordo}If $\partial Z^{r-1}_{p+r-1}\nsubseteq
Z^{r-1}_{p-1}$, then
$Z^{r-1}_{p-1}+\partial
Z^{r-1}_{p+r-1}=\mathbb{Z}^{r}_{p}.$
\end{lem}
\dem Since $\partial Z^{r-1}_{p+r-1}\nsubseteq
Z^{r-1}_{p-1}$ then $Z^{r-1}_{p-1}+\partial
Z^{r-1}_{p+r-1}$ is a submodule of
$$Z^{r}_{p}=\mathbb{Z}((t))[\mu^{p+1,r}\sigma_k^{p+1,r},\mu^{p,{r-1}}\sigma_{k}^{p,{r-1}},\ldots
,\mu^{\kappa,{r-p-1+\kappa}}\sigma_{k}^{\kappa,{r-p-1+\kappa}}]$$
where
$\kappa$ is the first column associated to a $k$-chain, 
but it is not a submodule of
$$Z^{r-1}_{p-1}=\mathbb{Z}((t))[\mu^{p,{r-1}}\sigma_{k}^{p,{r-1}},\mu^{p-1,{r-2}}\sigma_{k}^{p-1,{r-2}},\ldots
,\mu^{\kappa,{r-p-1+\kappa}}\sigma_{k}^{\kappa,{r-p-1+\kappa}}].$$
Then $\mu^{(p+1),r}=1$ and $Z^{r-1}_{p-1}+\partial
Z^{r-1}_{p+r-1}$ contains a multiple 
of $\sigma_k^{p+1,r}$ over $\mathbb{Z}((t))$. 
We will show that  $\sigma_k^{p+1,r}\in Z^{r-1}_{p-1}+\partial
Z^{r-1}_{p+r-1}$. Note that
\begin{equation}\label{eq:1}\partial Z^{r-1}_{p+r-1}=\mathbb{Z}((t))[\mu^{p+r,{r-1}}\partial \sigma_{k+1}^{p+r,{r-1}},
\mu^{p+r-1,{r-2}}\partial\sigma_{k+1}^{p+r-1,{r-2}},\ldots,\mu^{\overline{\kappa},\overline{\kappa}-p-1}
\partial \sigma_{k+1}^{\overline{\kappa},\overline{\kappa}-p-1}].\end{equation}
where $\overline{\kappa}$ is the first column associated to a
$(k+1)$-chain.  For $\xi=0,\ldots
,p+r-\overline{\kappa}$ with $\mu^{p+r-\xi,r-1-\xi}=1$ we have
$\Delta^{r-1-\xi}_{i,p+r-\xi}=0$ for all $i>p+1$ and hence
$$\partial
\sigma_{k+1}^{p+r-\xi,r-1-\xi}=\Delta^{r-1-\xi}_{p+1,p+r-\xi}\sigma_{k}^{p+1,{r-1-\xi}}+\cdots+\Delta^{r-1-\xi}_{\kappa,p+r-\xi}\sigma_{k}^{\kappa,{r-1-\xi}}.$$
 In fact, the boundaries $\partial
\sigma_{k+1}^{p+r-\xi,r-1-\xi}$ with
$\Delta^{r-1-\xi}_{i,p+r-\xi}\neq 0$ for some $i>p+1$ correspond
exactly to the columns which have the primary pivots below the
$(p+1)$-st row and therefore $\mu^{p+r-\xi,r-1-\xi}=0$.

Hence, for $\xi=0,\ldots ,p+r-\overline{\kappa}$, when
$\mu^{p+r-\xi,r-1-\xi}=1$ we have
\begin{equation}\label{eq:a1}Z^{r-1}_{p-1}+[\partial
\sigma_{k+1}^{p+r-\xi,r-1-\xi}]=Z^{r-1}_{p-1}+[\Delta^{r-1-\xi}_{p+1,p+r-\xi}\sigma_{k}^{p+1,{r-1-\xi}}+\cdots+\Delta^{r-1-\xi}_{\kappa,p+r-\xi}\sigma_{k}^{\kappa,{r-1-\xi}}]
.\end{equation} On the other hand, since $Z^{r-1}_{p-1}+[\partial
\sigma_{k+1}^{p+r-\xi,r-1-\xi}]\subset Z^{r-1}_{p-1}+\partial
Z^{r-1}_{p+r-1}$ then
\begin{equation}\label{eq:a2}Z^{r-1}_{p-1}+[\partial
\sigma_{k+1}^{p+r-\xi,r-1-\xi}]=[\ell_\xi\sigma_k^{p+1,r},\mu^{p,{r-1}}\sigma_{k}^{p,{r-1}},\mu^{p-1,{r-2}}\sigma_{k}^{p-1,{r-2}},\ldots
,\mu^{\kappa,{r-p-1+\kappa}}\sigma_{k}^{\kappa,{r-p-1+\kappa}}].
\end{equation}
The coefficient of $h_{k}^{p+1}$ on the set of generators of the
$\mathbb{Z}((t))$-module in (\ref{eq:a1}) is
$\Delta^{r-1-\xi}_{p+1,p+r-\xi}$. On the other
hand, the coefficient of $h_{k}^{p+1}$ on the set of the
generators of the $\mathbb{Z}((t))$-module in (\ref{eq:a2}) is $\ell_\xi$. 
Since for each $\xi=0,\ldots, p+r-\overline{\kappa}$, $\Delta^{r-1-\xi}_{i,p+r-\xi}=0$ for all $i>p+1$ then $\Delta^{r-1-\xi}_{p+1,p+r-\xi}$ is either a pivot or a zero entry. Note that the entries $\Delta^{r-1-\xi}_{p+1,p+r-\xi}$ can not be all zeros, since it would contradict the hipothesis of $\partial Z^{r-1}_{p+r-1}\nsubseteq
Z^{r-1}_{p-1}$.  It follows from Lemma \ref{teo:primarypivots} that if $\Delta^{r-1-\xi}_{p+1,p+r-\xi}$ is nonzero then it is invertible in $\mathbb{Z}((t))$. Then $\sigma_k^{p+1,r}\in Z^{r-1}_{p-1}+\partial
Z^{r-1}_{p+r-1}$ and hence $Z^{r-1}_{p-1}+\partial
Z^{r-1}_{p+r-1}=Z^{r}_{p}$.
\cqd

\begin{teo}\label{teo:inducaodoe} The matrix $\Delta ^r$ obtained from the SSSA applied to $\Delta$ determines 
$$E^{r}_{p}=\frac{Z^r_{p}}{Z^{r-1}_{p-1}+\partial
Z^{r-1}_{p+r-1}}. $$
More specifically, $E^r_p$ is either zero or a finitely
generated $\mathbb{Z}((t))$-module whose generator corresponds to a $k$-chain associated to column $(p+1)$ of $\Delta^r$.
\end{teo}
\dem The entry $\Delta^{r}_{p-r+1,p+1}$ on the $r$-th diagonal plays a
crucial role in determining the generators of $E^{r}_{p}$. We will  analyze each possibility of  $\Delta^{r}_{p-r+1,p+1}$ and prove how to determine $E^r_p$ for each case.
\begin{enumerate}
\item The entry $\Delta^{r}_{p-r+1,p+1}$ is identified by the SSSA
as a primary pivot, a change-of-basis pivot or a zero entry with column of zeros below it.

In these cases $\Delta^{r}_{s,p+1}=0$ for all $s>p-r+1$ and hence the generator $\sigma_k^{p+1,r}$ corresponding to the $k$-chain associated to column $(p+1)$ in $\Delta^{r}$ is a generator of $Z^{r}_{p}$.
Thus we must analyze row $(p+1)$. We have the following possibilities:
\begin{enumerate}
\item $\partial Z^{r-1}_{p+r-1}\subseteq Z^{r-1}_{p-1}$,
i.e, all the boundaries of the elements in
$Z^{r-1}_{p+r-1}$ are above row $p$.

In this case, as before, by Proposition \ref{prop:formuladoz}
$E^{r}_{p}=\displaystyle\mathbb{Z}((t))[\sigma_k^{p+1,r}]$.

\item $\partial Z^{r-1}_{p+r-1}\nsubseteq
Z^{r-1}_{p-1}$, i.e, there exist elements in $
Z^{r-1}_{p+r-1}$ whose boundary has a nonzero entry in
row $(p+1)$. By Proposition \ref{prop:formuladoz} and Lemma
\ref{prop:formulabordo} $E^{r}_{p}=0$. 
\end{enumerate}

Note that if $\Delta^{r}_{p-r+1,p+1}$ has been identified by the
SSSA as a primary pivot then \\ $\partial
Z^{r-1}_{p+r-1}\subseteq Z^{r-1}_{p-1}$. For more details see \cite{CdRS}.

In fact, by Lemma \ref{lemapivolinha}, row $(p+1)$ cannot contain
a primary pivot since we have assumed that column $(p+1)$ has a
primary pivot. Therefore, the entries of these $h_{k+1}$ columns in
row $(p+1)$ must be zeroes. It follows that $\partial
Z^{r-1}_{p+r-1}$ does not contain in its set of
generators the generator $\sigma_k^{p+1,r}$.

\item The entry $\Delta^{r}_{p-r+1,p+1}$ is an entry above a primary
pivot, i.e. there exists $s>p-r+1$ such that
$\Delta^{r}_{s,p+1}$ is a primary pivot.  In this case,  the generator $\sigma^{p+1,r}_{k}$ corresponding to the $k$-chain associated to column $(p+1)$ is not a generator of $Z^{r}_{p}$  and hence $Z^{r-1}_{p-1}=Z^{r}_{p}$. It
follows that $E^{r}_{p}=0$.

\item The entry $\Delta^{r}_{p-r+1,p+1}$ is not in $\Delta^r_k$.
This includes the case where $p-r+1<0$, i.e,
$\Delta^{r}_{p-r+1,p+1}$ is not on the matrix $\Delta^r$.

The analyzes of $E^r_p$ is very similar to the previous one, i.e, we
have two possibilities:
\begin{enumerate}
\item There is a primary pivot in column $(p+1)$ in a
diagonal $\overline{r}<r$. In this case the generator corresponding
to the $k$-chain associated to column $(p+1)$, $\sigma^{p+1,r}_{k}$
is not a generator of $Z^{r}_{p}$. Hence
$Z^{r-1}_{p-1}=Z^{r}_{p}$ and $E^{r}_{p}=0$.

\item All the entries in $\Delta^r$ in column $(p+1)$
in diagonals lower than $r$ are zero, i.e, the generator
corresponding to the $k$-chain associated to column $(p+1)$,
$\sigma_k^{p+1,r}$ in $\Delta^{r}$ is a generator of
$Z^{r}_{p}$. Then we have to analyze row $(p+1)$.

\begin{enumerate}
\item If $\partial Z^{r-1}_{p+r-1}\subseteq Z^{r-1}_{p-1}$
then, by Proposition \ref{prop:formuladoz},
$E^{r}_{p}=\displaystyle\mathbb{Z}((t))[\sigma_k^{p+1,r}]$.

\item If $\partial Z^{r-1}_{p+r-1}\nsubseteq
Z^{r-1}_{p-1}$ then, by Proposition \ref{prop:formuladoz}
and Lemma \ref{prop:formulabordo}, $E^{r}_{p}=0$.
\end{enumerate}
\end{enumerate}

\end{enumerate}

\cqd

We will describe how the SSSA applied to $\Delta$
induces the differentials $d^r_p:E^{r}_p \to E^r_{p-r}$ of
spectral sequence.

\begin{teo}\label{teo:interpretacaodr} If $E^{r}_{p}$ and
$E^{r}_{p-r}$ are both nonzero, then the map $d^r_p:E^{r}_{p}\to
E^{r}_{p-r}$ is induced by the multiplication by
the entry $\Delta^r_{p-r+1,p+1}$.
\end{teo} 

\dem  Suppose that $E^{r}_{p}$ and $E^{r}_{p-r}$ are both nonzero. Let $\delta^{r}_{p}:\mathbb{Z}((t))[\sigma_{k}^{p+1,r}]\to\mathbb{Z}((t))[\sigma_{k-1}^{p-r+1,r}]$ be the multiplication by the entry $\Delta^r_{p-r+1,p+1}$ and ${\tilde{\delta}}^{r}_{p}$ the induced map in $E^{r}_{p}$. We
must show that
$$\displaystyle\frac{\Ker \tilde{\delta}^{r}_{p}}{\Ima
\tilde{\delta}^{r}_{p+r}}\cong E^{r+1}_p.$$ 
Since we are considering $E^{r}_{p}$ nonzero,
it follows from Theorem \ref{teo:inducaodoe}, that we must consider
three cases for the entry $\Delta^r_{p-r+1,p+1}$: primary pivot,
change-of-basis pivot and zero with a column of zeroes below it.
However, if $\Delta^r_{p-r+1,p+1}$ is a change-of-basis pivot then
there exists a primary pivot in row $(p-r+1)$ on a diagonal below
the $r$-th diagonal and hence $E^{r}_{p-r}=0$. Hence, whenever
$E^{r}_p$ and $E^r_{p-r}$ are both nonzero, the entry
$\Delta^r_{p-r+1,p+1}$ in $\Delta^r$ is either a primary pivot or a
zero with a column of zero entries below it.

In this case $E^{r}_{p}=\displaystyle\mathbb{Z}((t))[\sigma_k^{p+1,r}]$
and $E^{r}_{p-r}=\displaystyle\mathbb{Z}((t))[\sigma_{k-1}^{p-r+1,r}]$.

\begin{enumerate}
\item Suppose $\Delta^r_{p-r+1,p+1}$ is a primary pivot.

Since $\delta^{r}_{p}: \mathbb{Z}((t))[\sigma_{k}^{p+1,r}]\to
\mathbb{Z}((t))[\sigma_{k-1}^{p-r+1,r}]$ is multiplication by
$\Delta^r_{p-r+1,p+1}$, which is invertible in $\mathbb{Z}((t))$, then $\Ker \delta^{r}_{p}=0$. Since $\tilde{\delta}^{r}_{p}=\delta^{r}_{p}$, then
$\displaystyle\frac{\Ker\tilde{\delta}^{r}_{p}}{\Ima \tilde{\delta}^{r}_{p+r}}=0$.
On the other hand, since $\Delta^r_{p-r+1,p+1}$ is a primary pivot
then $\sigma_{k}^{p+1,r+1}=\sigma_{k}^{p+1,r}\notin Z^{r+1}_p$. Thus
$Z^{r+1}_p= Z^{r}_{p-1}$ and $E^{r+1}_p=0$.

\item Suppose $\Delta^{r}_{p-r+1,p+1}=0$ with a column of zeroes below it.
  In this case $\Ker \delta ^r_p \cong E^{r}_{p}=\Ker \tilde{\delta} ^r_p $ and
  $\sigma_k^{p+1,r}=\sigma_k^{p+1,r+1}$. There are three cases to consider:
  \begin{enumerate}
    \item If $\Delta^r_{p+1,p+r+1}$ is an entry above a primary
      pivot then we have $E^{r}_{p+r}=0$ and hence
      $\Ima\tilde{\delta}^{r}_{p+r}=0$. Thus,
      $$\displaystyle\frac{\Ker\tilde{\delta}^{r}_{p}}{\Ima\tilde{\delta}^{r}_{p+r}}= E^r_p.$$
      On the other hand, since $\mu^{p+r+1,r}=0$ then
      $E^{r+1}_{p}=E^{r}_{p}$.
      
       \item If $\Delta^{r}_{p+1,p+r+1}$ is a primary pivot then
      $E^{r}_{p+r}=\mathbb{Z}((t))[\sigma_k^{p+r+1,r}]$.
      Therefore $\delta^{r}_{p+r}$ is an isomorphism and hence
      \[\frac{\Ker\tilde{\delta}^{r}_{p}}{\Ima\tilde{\delta}^{r}_{p+r}}\cong
      \frac{\mathbb{Z}((t))[\sigma_k^{p+1,r}]}{\mathbb{Z}((t))[\sigma_k^{p+1,r}]}=0.\]
      On the other hand, since $\Delta^{r}_{p-r+1,p+1}$ is zero with a column of zero
entries below it then $\sigma_k^{p+1,r+1}\in Z^{r+1}_{p,k-p}$ and
hence $Z^{r}_{p-1}\varsubsetneq Z^{r+1}_{p,k-p}$. Moreover,
since $E^{r}_{p}=\displaystyle\mathbb{Z}((t))[\sigma_k^{p+1,r}]$ then
$\partial Z^{r-1}_{p+r-1}\subseteq
Z^{r-1}_{p-1}$. But the difference between $\partial
Z^{r-1}_{p+r-1}$ and $\partial
Z^{r}_{p+r}$ is that the last one includes the boundary
of column $(p+r+1)$. The element in column $(p+r+1)$ and row $(p+1)$
is $\Delta^r_{p+1,p+r+1}$.  Since $\Delta^r_{p+1,p+r+1}$ is a primary pivot then  $\partial Z^{r}_{p+r}\nsubseteq Z^{r}_{p-1}$ and
$E^{r+1}_{p}=0$.

    \item If $\Delta^{r}_{p+1,p+r+1}=0$ with a column of zero entries below
    it then $\Ima\delta^{r}_{p+r}=0$ and
      $$\displaystyle\frac{\Ker\tilde{\delta}^{r}_{p}}{\Ima\tilde{\delta}^{r}_{p+r}}= E^r_p.$$
      Analogously to the previous case, $\sigma_k^{p+1,r+1}\in Z^{r+1}_{p,k-p}$ and hence $Z^{r}_{p-1}\varsubsetneq Z^{r+1}_{p}$. Moreover,
$\partial Z^{r-1}_{p+r-1}\subseteq
Z^{r-1}_{p-1}$ and the only difference between $\partial
Z^{r-1}_{p+r-1}$ and $\partial
Z^{r}_{p+r}$ is that the last one includes the boundary
of column $(p+r+1)$. Since the element in column $(p+r+1)$ and row $(p+1)$
is $\Delta^r_{p+1,p+r+1}=0$  then $\partial
Z^{r}_{p+r}\subseteq Z^{r}_{p-1}$ and
$E^{r+1}_{p}=\displaystyle \mathbb{Z}((t))[\sigma_k^{p+1,r}]$.
      
\end{enumerate}
\end{enumerate}

Note that the case where $\Delta^{r}_{p+1,p+r+1}$ is a change of basis pivot does not have to be considered, since in this case $E^{r}_{p}=0$.

Therefore, in all cases
\[\frac{\Ker d^{r}_{p}}{\Ima d^{r}_{p+r}} \ = \ E^{r+1}_{p} \ = \ \frac{\Ker\tilde{\delta}^{r}_{p}}{\Ima\tilde{\delta}^{r}_{p+r}}.\]
\cqd

\begin{cor}\label{cor:novikov1} 
Each non zero differential $d^{r}_{p}$ of the spectral sequence $(E^{r},d^{r})$ is an isomorphism.
\end{cor}
\dem
In fact, by the proof of Theorem \ref{teo:interpretacaodr}, non zero differentials of $(E^{r},d^{r})$ are induced by primary pivots. Theorem
\ref{teo:primarypivots} states that each primary pivot produced by the SSSA is an invertible polynomial, hence each  induced non zero differential is an isomorphism.
\cqd

The next corollary states that the spectral sequence $(E^{r},d^{r})$ converges to the Novikov homology of $(\mathcal{N}_{\ast}(f),\Delta)$.

\begin{cor}\label{cor:novikov2}
 The modules $E^{\infty}_{p,q}$ of the spectral sequence are free for all $p$ and $q$. Moreover, the spectral sequence converges to the Novikov homology.
\end{cor}
\dem
By Corollary \ref{cor:novikov1},  the non zero differentials $d^{r}: E^{r}_{p} \rightarrow E^{r}_{p-r}$ of the spectral sequence are isomorphisms.
Since,  $E^{r+1}_{p} \cong \Ker \ d^{r}_{p} / \Ima \ d^{r}_{p-r}$, it follows that 
 the modules $E^{r}_{p}$'s are free for all $r\geq 0$ and $p \geq 0$.
Therefore, 
 $$ E^{\infty}_{p,q}\approx GH_*(\mathcal{N})_{p,q}=\frac{F_pH_{p+q}(\mathcal{N})}{F_{p-1}H_{p+q}(\mathcal{N})} \approx H^{Nov}_ {\ast}(M,f).$$
 by equation (\ref{eq:seqspec}).
\cqd

\section{Cancellation of Critical Points}\label{sec6}

The results obtained for a  spectral sequence of a filtered two dimensional Novikov chain complex $(\mathcal{N}_{\ast}(f),\Delta)$ associated to a circle-valued Morse function $f:M\rightarrow S^{1}$ lay the groundwork for dynamical exploration of connections in a negative gradient flow associated to $f$. 
In Sections \ref{caracterizacao} and \ref{ultimo} we proved that the SSSA produces from $(\mathcal{N}_{\ast}(f),\Delta)$ a sequence
of Novikov matrices from which the modules and differentials of the spectral sequence $(E^r,d^r)$ may be retrieved. By Corollary \ref{cor:novikov1}, if $d^{r}:E_{p}^{r}\rightarrow E^{r}_{p-r}$ is a non zero differential, then the modules $E^{r+1}_{p}$ and $E^{r+1}_{p-r}$ are zero. We refer to this situation as an {\it algebraic cancellation}. Hence, as one "turns the pages" of the spectral sequence, i.e.  considers progressively modules $E^r$, one  observes   algebraic cancellations occurring within the $E^r$'s.

We will now show how the dynamics follows the algebraic unfolding of the spectral sequence, that is, how
these algebraic cancellations can be associated to dynamical cancellations of singularities of the negative gradient flow on $M$. Let $f$ be a circle-valued Morse function and  $(\mathcal{N}_{\ast}(f),\Delta)$  the Novikov complex associated to $f$. Let the Morse function $F:\overline{M} \rightarrow \mathbb{R}$ be the lift of $f$ to the infinite cyclic covering space $\overline{M}$, as defined in Section \ref{novikovcomplex}.

Consider a finest filtration $\mathcal{F}= \{\mathcal{F}_{p}\mathcal{N}\}_{p\in P}$ on $(\mathcal{N}_{\ast}(f),\Delta)$, where 
$$F_{p}\mathcal{N} = \mathbb{Z}((t))\Big[h^1_{k_1},h^{2}_{k_2}\ldots , h^{p+1}_{k_{p+1}}\Big]$$
and 
$P=\{0,1,\ldots,m\}$ is an indexing set for the filtration $\mathcal{F}$, where $\#  Crit(f)=m+1$.
From the  filtration $\mathcal{F}$ one can induce an  infinite filtration $\{\mathcal{F}_{\lambda , p} \}_{\lambda \in \mathbb{Z}, p\in P}$ in the covering space $\overline{M}$ which corresponds to

$$
\mathcal{F}_{\lambda , p}\mathcal{N} = \mathbb{Z} \Bigg\langle
t^{\delta}h^{j}_{k} \ \  \Bigg|  \ 
\begin{array}{l}
\delta = \lambda \ \text{and} \ j \leq p+1 \  \ \text{or}\\
\delta > \lambda \ \text{and} \ j \in \{1,\ldots,m\}
\end{array}
 \Bigg\rangle
$$
where $t^{\lambda} h^{p+1}_{k}$ represents the lift of the point $h^{p+1}_{k}$ at level $\lambda$, for $\lambda\in \mathbb{Z}$.
Note that, for $p\in\{0,1,\ldots,m\}$ and $\lambda \in \mathbb{Z}$, we have 
$$ \mathcal{F}_{\lambda , p-1}\mathcal{N} \subset \mathcal{F}_{\lambda , p}\mathcal{N}  \ \ \text{and} \ \  \mathcal{F}_{\lambda , m}\mathcal{N} \subset \mathcal{F}_{\lambda -1 , 0}\mathcal{N},$$
and also,
$$  \mathcal{F}_{\lambda , p}\mathcal{N} \setminus \mathcal{F}_{\lambda , p-1}\mathcal{N} = t^{\lambda}h^{p+1}_{k} \ \ \text{and} \ \  \mathcal{F}_{\lambda -1 , 0}\mathcal{N} \setminus \mathcal{F}_{\lambda , m}\mathcal{N} = t^{\lambda-1}h^{1}_{k}.$$

Without loss of generality, one can assume  that $f$ has one critical point per critical level set. 
Let $c_{p}$ be a critical value  of $f$ such that $f(h_{k}^{p}) =c_p$.
Hence, $F$ also has one critical point per critical level set and $c_{\lambda,p}:=c_{p} -\lambda$ is a critical value  of $F$ such that $F(t^\lambda h_{k}^{p}) =  c_{\lambda,p}$, for all $\lambda \in \mathbb{Z}$.

\begin{ex}\label{exemplo-canc}
Consider the negative gradient flow on a torus  as  illustrated in Figure \ref{fig:exemplo-art}  which is associated to a  circle-valued Morse function $f$. 
\begin{figure}[!htb]
\centering
\begin{minipage}[b]{0.5\linewidth}
\includegraphics[width=7.5cm]{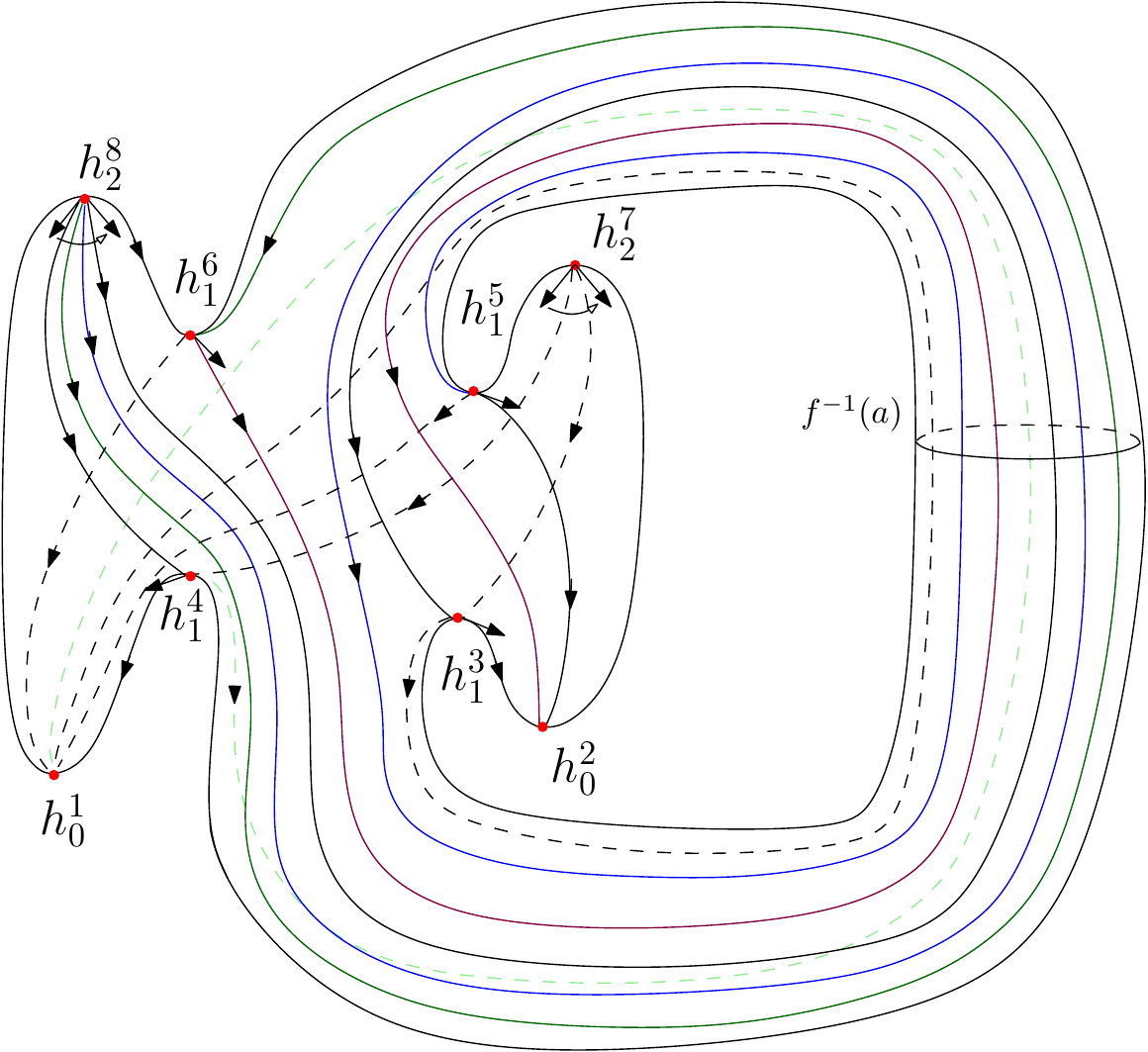}
\caption{Circle-valued Morse function on the torus.}
\label{fig:exemplo-art}
\vspace{1.9cm}
\includegraphics[width=7cm]{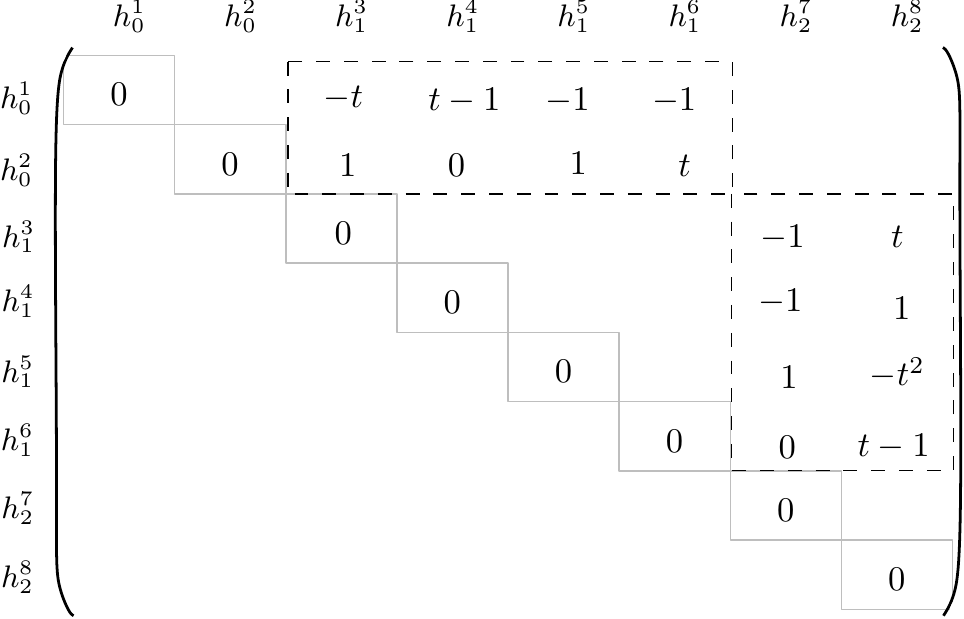}
\caption{Novikov differential $\Delta$.}
\label{fig:ddelta0ex2}
\end{minipage} 
\hspace{0.5cm}
\begin{minipage}[b]{0.45\linewidth}
\centering
\includegraphics[scale=0.6]{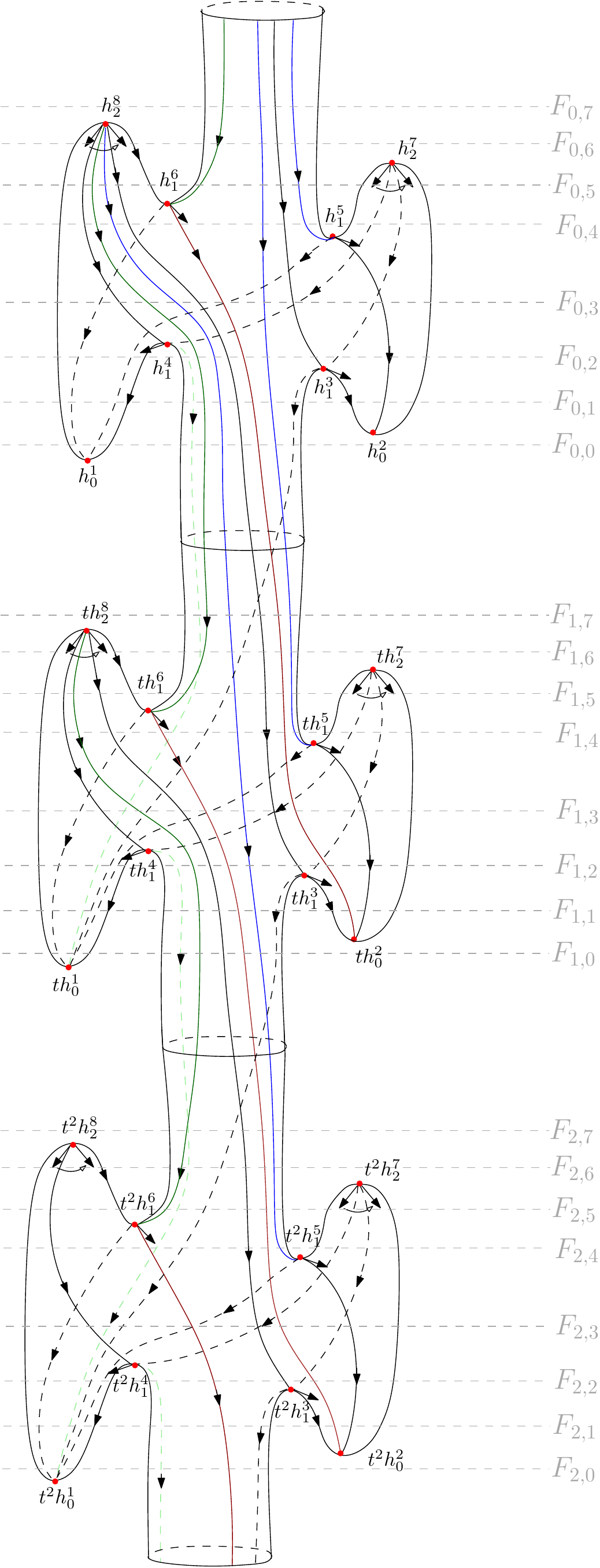}
\vspace{-0.3cm}
\caption{Filtration in the infinite cyclic covering.}
\label{fig:exemplo-art-levantamento}
\end{minipage}
\end{figure}
Let $(\mathcal{N}_{\ast}(f),\Delta)$ be the  Novikov complex associated to $f$ and define  a filtration   $\mathcal{F}=\{\mathcal{F}_{p}\mathcal{N}\}$
 by:
$ F_{0}\mathcal{N} = \mathbb{Z}((t))\big[h^1_0\big] $,
 $ F_{1}\mathcal{N} = \mathbb{Z}((t))\big[h^1_0, h^2_0\big] $, 
$ F_{2}\mathcal{N} = \mathbb{Z}((t))\big[h^1_0, h^2_0,h^3_1\big] $, 
$ F_{3}\mathcal{N} = \mathbb{Z}((t))\big[h^1_0, h^2_0,h^3_1,h^4_1\big] $, 
$ F_{4}\mathcal{N} = \mathbb{Z}((t))\big[h^1_0, h^2_0,h^3_1,h^4_1,h^5_1\big] $, 
$ F_{5}\mathcal{N} = \mathbb{Z}((t))\big[h^1_0, h^2_0,h^3_1,h^4_1,h^5_1,h^6_1\big] $, 
$ F_{6}\mathcal{N} = \mathbb{Z}((t))\big[h^1_0, h^2_0,h^3_1,h^4_1,h^5_1,h^6_1, h^7_2\big] $,
 $ F_{7}\mathcal{N} = \mathbb{Z}((t))\big[h^1_0, h^2_0,h^3_1,h^4_1,h^5_1,h^6_1, h^7_2,h^{8}_2\big] $.
The Novikov differential $\Delta$ is presented  in Figure \ref{fig:ddelta0ex2} and we illustrated the associated covering space enriched with the  finest filtration  $\{\mathcal{F}_{\lambda , p} \}$  in Figure \ref{fig:exemplo-art-levantamento}.

\end{ex}

The next result will provide an algorithmic approach to the dynamical cancellation of critical points of a circle-valued Morse function $f:M\rightarrow S^{1}$ by computing the spectral sequence of the associated Novikov chain complex of $f$.


\begin{teo}[Dynamical Cancellation Theorem via Spectral Sequences in the Novikov setting]\label{main} Let $f:M\rightarrow S^{1}$ be a circle-valued Morse function on a closed $n$-manifold $M$ and  $(\mathcal{N}_{\ast}(f),\Delta)$ be the Novikov chain complex associated to $f$. Let  $\mathcal{F}$ be a finest filtration in $(\mathcal{N}_{\ast}(f),\Delta)$ and let $(E^r, d^r)$ be the associated spectral sequence for this filtered chain complex. 
There exists a  one-to-one correspondence between the algebraic cancellations of the modules $E^r$ with  dynamical cancellations of critical points of $f$, which produces a family of circle-valued Morse functions $f^{r}$, where $f^{r+1}$ is obtained from $f^r$ by removing the pairs of critical points corresponding to the cancelled modules of the $r$-th page of the spectral sequence.
 Moreover, for each $r \geq 0$, the $(r+1)$-th page of  the spectral sequence completely determines  all Novikov incidence coefficients, hence all new  flow lines,  between consecutive critical points of $f^{r+1}$.
\end{teo}

The following example illustrates how the dynamics follows the algebraic unfolding of the spectral sequence.

\begin{figure}[!htb]
\centering
\begin{minipage}[b]{0.5\linewidth}
\includegraphics[width=7.5cm]{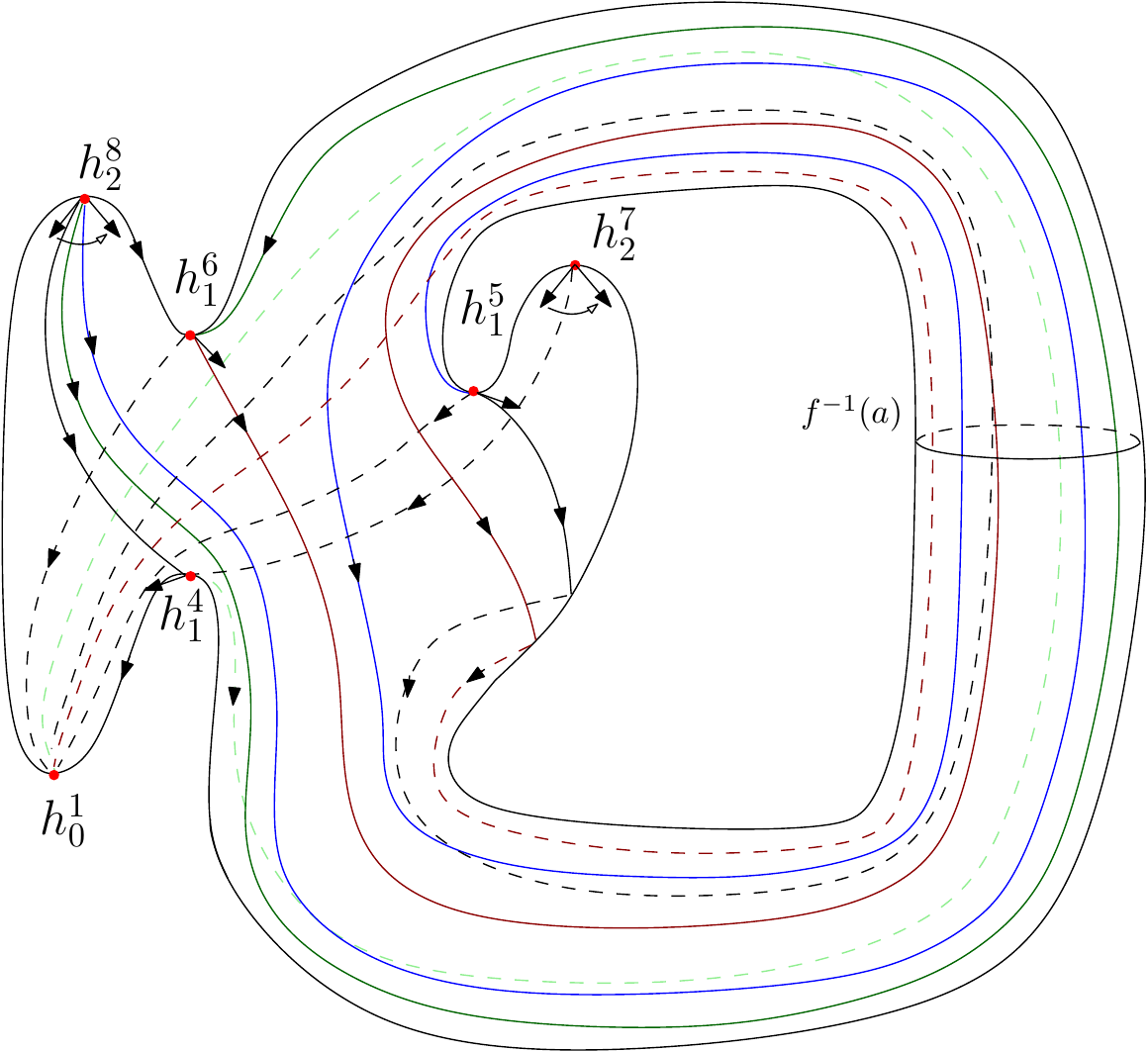} 
\caption{Cancellation of the pair $(h^{2}_{0},h^{3}_{1})$.}
\label{fig:exemplo-art-can1}
\end{minipage} 
\begin{minipage}[b]{0.45\linewidth}
\includegraphics[width=7cm]{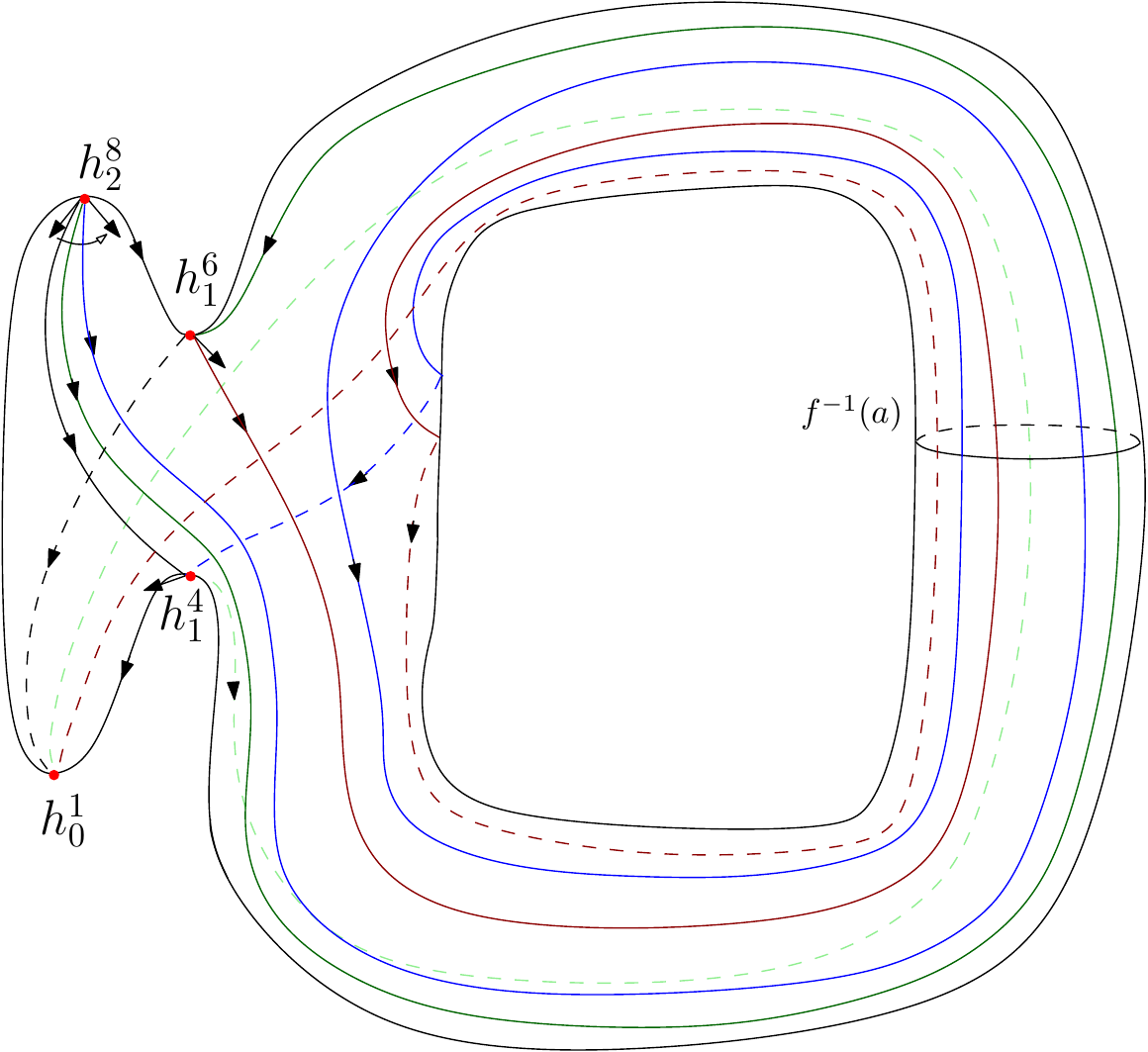}
\caption{Cancellation of  the pair $(h^{5}_{1},h^{7}_{2})$.}
\label{fig:exemplo-art-can2}
\end{minipage}
\end{figure}

\begin{figure}[!htb]
\centering
\begin{minipage}[b]{0.5\linewidth}
\includegraphics[width=7.5cm]{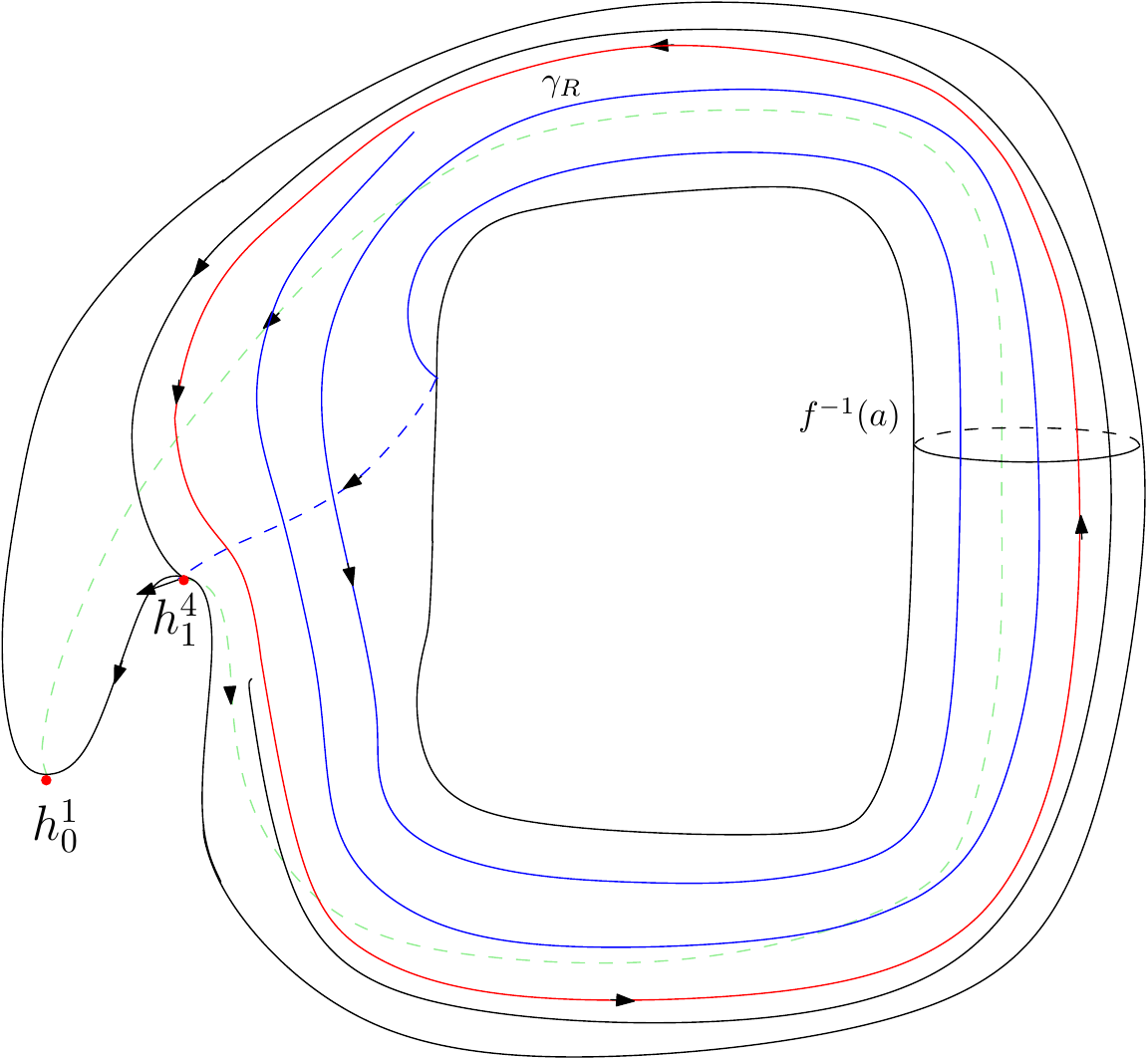}
\caption{Cancellation of the pair $(h^{6}_{1},h^{8}_{2})$.}
\label{fig:exemplo-art-can3}
\end{minipage} 
\begin{minipage}[b]{0.45\linewidth}
\includegraphics[width=6.7cm]{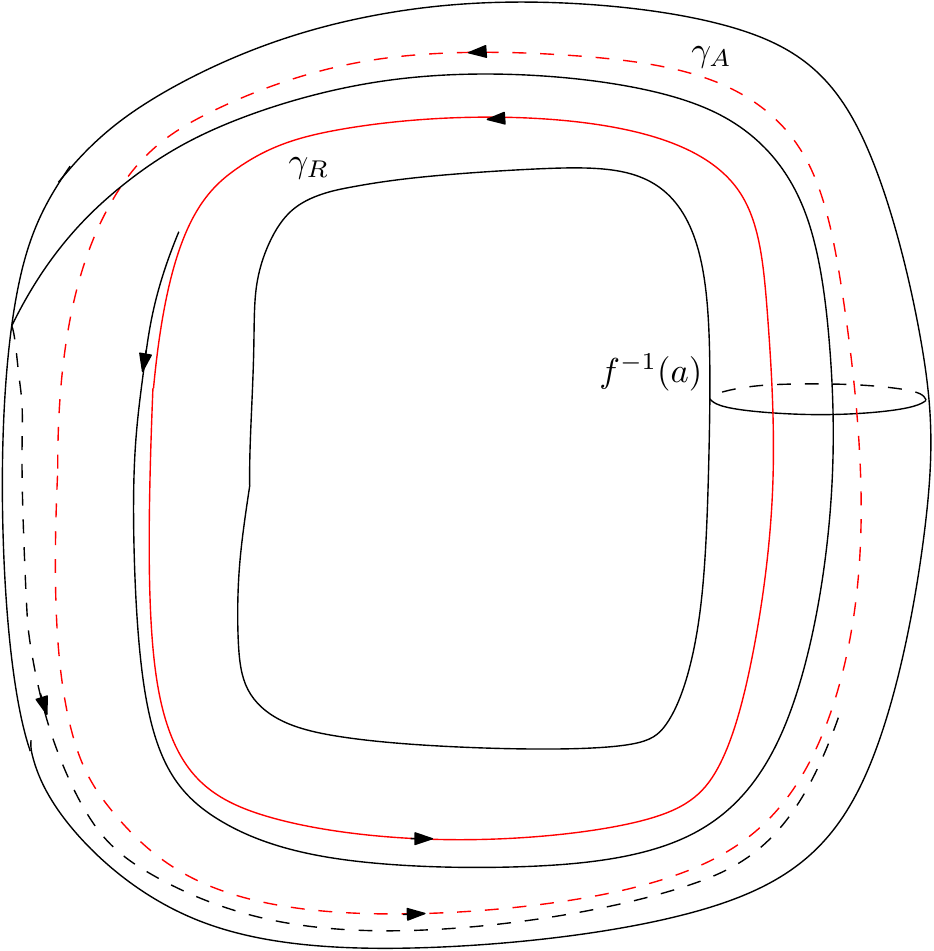}
\caption{Cancellation of the pair $(h^{1}_{0},h^{4}_{1})$.}
\label{fig:exemplo-art-can4}
\end{minipage}
\end{figure}

\begin{ex}\label{exemplocap6}
  Consider the circle-valued Morse function $f$ presented in  Example \ref{exemplo-canc}.
Applying the SSSA to the Novikov differential $\Delta$ in Figure \ref{fig:ddelta0ex2}, one  obtains the sequence of Novikov matrices $\Delta^{1},\cdots,\Delta^{6}$ presented in Figures \ref{fig:delta1ex2}, $\cdots$,  \ref{fig:delta6ex2}, respectively.

\begin{figure}[!htb]
\begin{minipage}[t]{0.5\linewidth}
 \includegraphics[width=7.5cm]{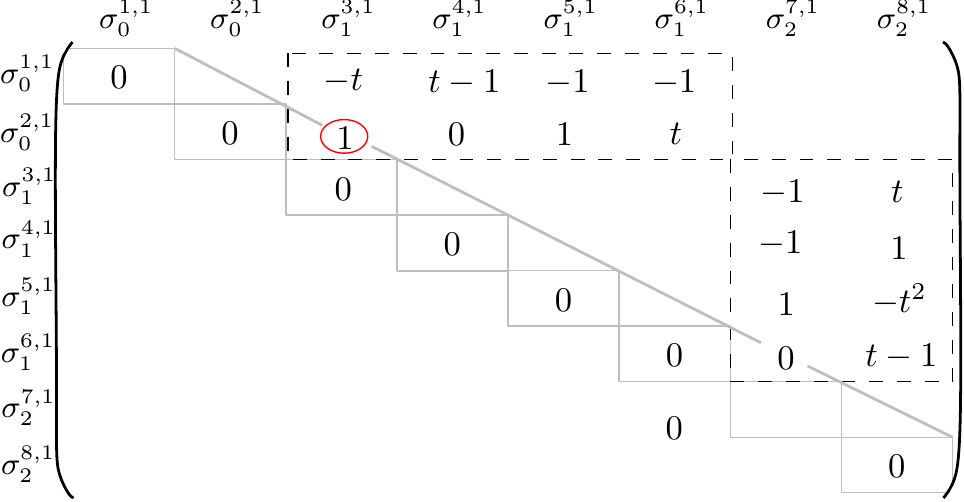}\\
{\footnotesize  $\sigma^{j,1}_{0}\!=\!h^j_{0}$ for $\!j\in \! J_{0}$;   $\sigma^{j,1}_{1}\!=\!h^j_{1}$  for $j\!\in \! J_{1}$;   $\sigma^{j,1}_{2}\!=\! h^j_{2}$  for $j\!\in \! J_{2}$.}\vspace{-0.3cm}
\caption{$\Delta^{1}$, marking primary pivots.}
\label{fig:delta1ex2}
\end{minipage} 
\begin{minipage}[t]{0.5\linewidth}
 \includegraphics[width=7.5cm]{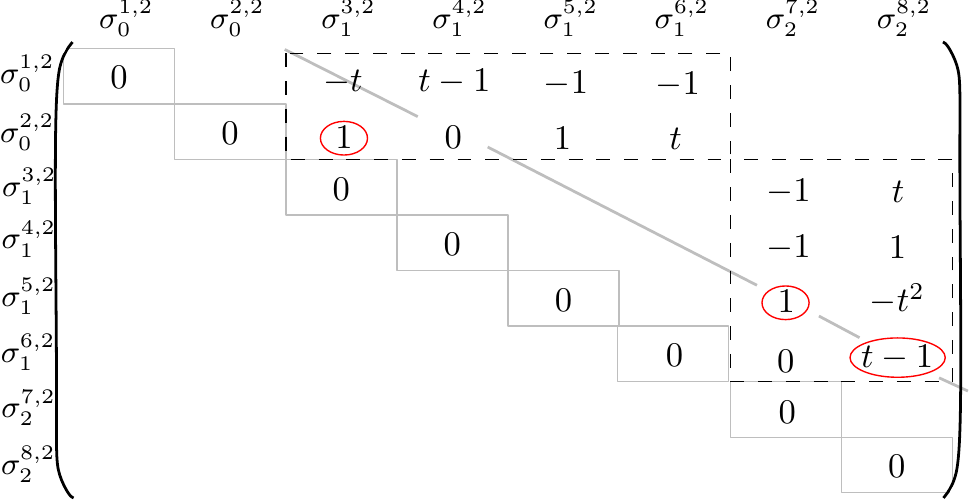}\\
  {\footnotesize    $\sigma^{j,2}_{0}\!=\!h^j_{0}$ for $\!j\in \! J_{0}$;   $\sigma^{j,2}_{1}\!=\!h^j_{1}$  for $j\!\in \! J_{1}$;   $\sigma^{j,2}_{2}\!=\! h^j_{2}$  for $j\!\in \! J_{2}$.}\vspace{-0.3cm}
\caption{$\Delta^{2}$, marking change-of-basis pivot.}
\label{fig:delta2ex2}
\end{minipage}
\end{figure}
\begin{figure}[!htb]
\begin{minipage}[t]{0.5\linewidth}
\includegraphics[width=7.5cm]{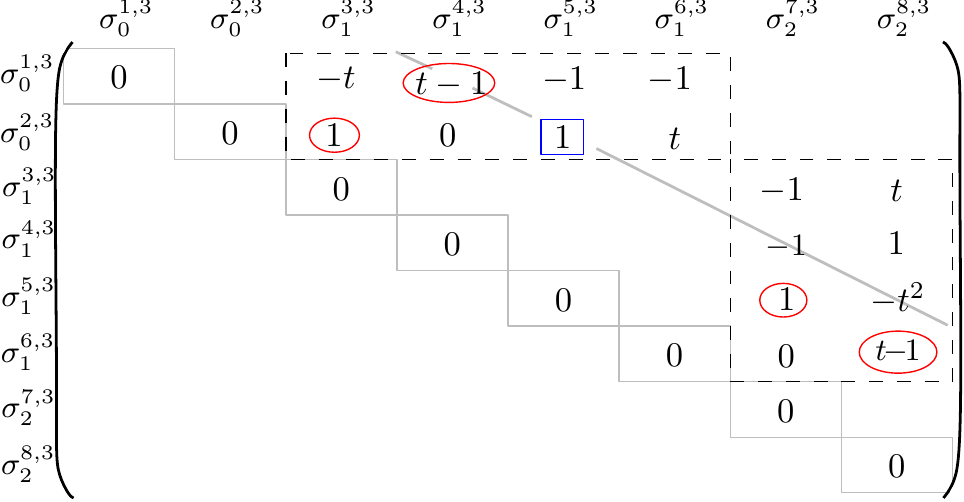}\\
 {\footnotesize    $\sigma^{j,3}_{0}\!=\!h^j_{0}$ for $\!j\in \! J_{0}$;   $\sigma^{j,3}_{1}\!=\!h^j_{1}$  for $j\!\in \! J_{1}$;    \\
$\sigma^{j,3}_{2}\!=\! h^j_{2}$  for $j\!\in \! J_{2}$.}\vspace{-0.3cm}
\caption{$\Delta^{3}$, marking pivots.}
\label{fig:delta3ex2}
\end{minipage} 
\begin{minipage}[t]{0.5\linewidth}
\includegraphics[width=7.5cm]{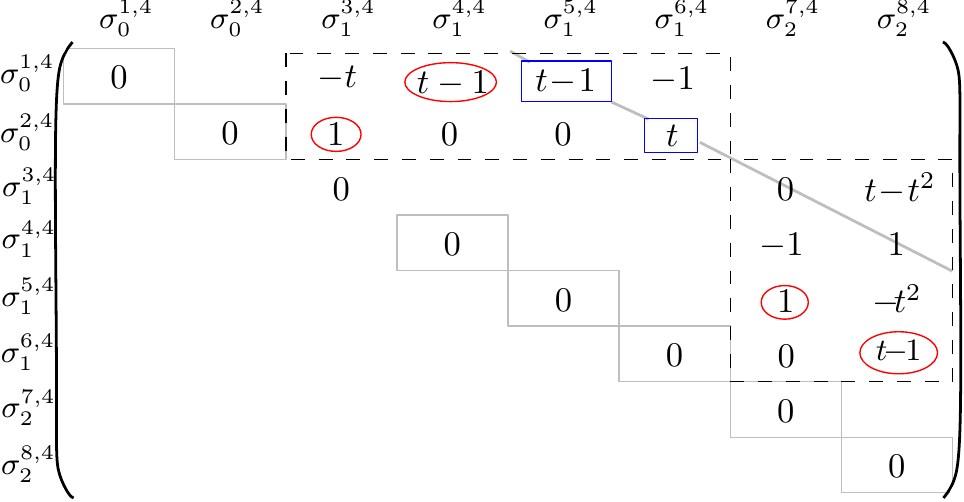}\\
   {\footnotesize    $\sigma^{5,4}_{1}\!=\!h^{5}_{1} - h^{3}_{1}$. \\
  $\sigma^{j,4}_{k} = \sigma^{j,3}_{k}$ for all the  remaining  $\sigma$'s.}\vspace{-0.3cm}
\caption{$\Delta^{4}$, marking change-of-basis pivot.}
\label{fig:delta4ex2}
\end{minipage}
\end{figure}
\begin{figure}[!htb]
\begin{minipage}[t]{0.5\linewidth}
\includegraphics[width=7.5cm]{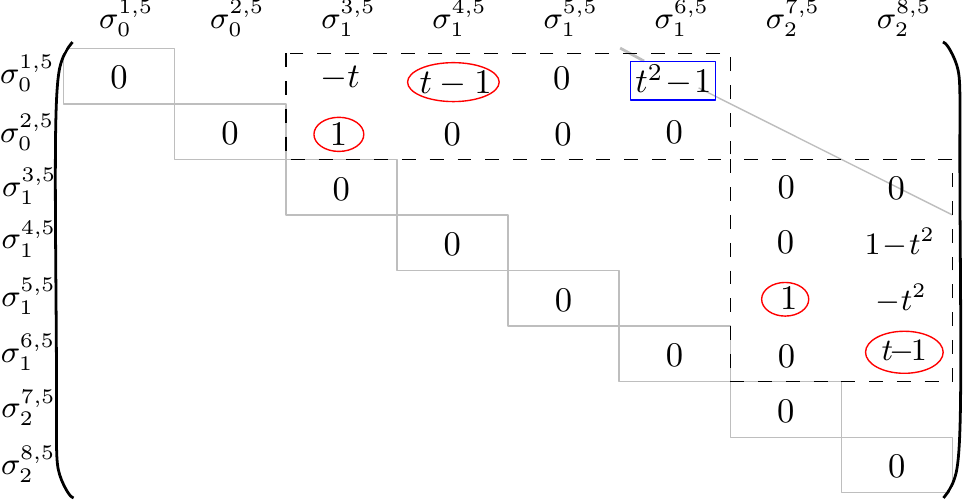}\\
     {\footnotesize    $\sigma^{5,5}_{1}\!=\!h^{5}_{1} - \ h^{4}_{1} - h^{3}_{1}$;    $\sigma^{6,5}_{1}\!=\!h^{6}_{1} -t h^{3}_{1}$;   \\
  $\sigma^{j,5}_{k} = \sigma^{j,4}_{k}$ for all the  remaining  $\sigma$'s.}\vspace{-0.3cm}
\caption{$\Delta^{5}$, marking change-of-basis pivot.}
\label{fig:delta5ex2}
\end{minipage} 
\begin{minipage}[t]{0.5\linewidth}
\includegraphics[width=7.5cm]{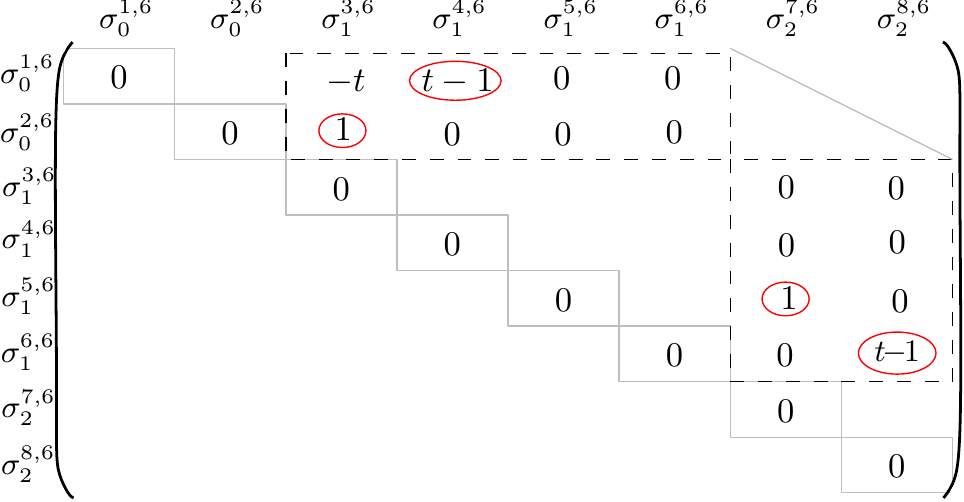}\\
       {\footnotesize    $\sigma^{6,6}_{1}\!=\!h^{6}_{1} +(-t-1) h^{4}_{1} -th^{3}_1$;       \\
  $\sigma^{j,6}_{k} = \sigma^{j,5}_{k}$ for all the  remaining  $\sigma$'s.}\vspace{-0.3cm}
\caption{$\Delta^{6}$ for Example \ref{exemplo}.}
\label{fig:delta6ex2}
\end{minipage}
\end{figure}

 The primary pivots detected by the SSSA induce the differentials of the spectral sequence, as we have shown in Theorem \ref{teo:interpretacaodr}. By Theorem \ref{teo:primarypivots}, these pivots are invertible, hence  responsible for algebraic cancellations of the modules of the spectral sequence. More specifically:
 \begin{itemize}
 \item the primary pivot $\Delta^{1}_{2,3}$ induces the differential $d^{1}_{2}: E^{1}_{2}\rightarrow E^{1}_{1}$ which in turn causes the  algebraic cancellation $E^{2}_{2}=E^2_1=0$;
 \item the primary pivot $\Delta^{2}_{5,7}$ induces the differential $d^{2}_{6}: E^{2}_{6}\rightarrow E^{2}_{4}$  which in turn causes  the algebraic cancellation $E^{3}_{6}=E^3_4=0$;
  \item the primary pivot $\Delta^{2}_{6,8}$ induces the differential $d^{2}_{7}: E^{2}_{7}\rightarrow E^{2}_{5}$  which in turn causes  the algebraic cancellation $E^{3}_{7}=E^3_5=0$;
    \item the primary pivot $\Delta^{3}_{1,4}$ induces the differential $d^{3}_{3}: E^{3}_{3}\rightarrow E^{3}_{0}$  which in turn causes  the algebraic cancellation $E^{4}_{3}=E^4_0=0$.
\end{itemize}

 Figures \ref{fig:exemplo-art-can1},  \ref{fig:exemplo-art-can2},  \ref{fig:exemplo-art-can3} and \ref{fig:exemplo-art-can4}  show the dynamical cancellations of pairs of critical points 
  $(h^{2}_{0},h^{3}_{1})$, $(h^{5}_{1},h^{7}_{2})$, $(h^{6}_{1},h^{8}_{2})$ and  $(h^{1}_{0},h^{4}_{1})$  corresponding to the above algebraic cancellations, as explained in what follows. 
\end{ex}

On the first page of the spectral sequence, the differential $d^{1}_{2}$ determines the algebraic cancellation $E^{2}_{2}=E^2_1=0$.
Corresponding to this algebraic cancellation, we have the dynamical cancellation of the pair  $(h^{2}_{0},h^{3}_{1})$, which determines a circle valued Morse function $f^{2}$ such that $Crit(f^{2}) = Crit(f)\setminus \{ {h^{2}_{0},h^{3}_{1}} \}$ and  is equal (up to isotopy) to $f$ outside a neighborhood of the connection of this critical points. In the corresponding negative gradient flow $\varphi^2$ (see Figure \ref{fig:exemplo-art-can1}), we note  birth of  flow lines from $h^{5}_{1}$   to $h^{1}_{0}$ and from $h^{6}_{1}$ to $h^{1}_{0}$.   Moreover, we have the following relation between the Novikov incidence coefficients
 $$N(h^{5}_{1},h^{1}_{0};f^{2})=
 N(h^{5}_{1},h^{2}_{0};f^{1})t+ N(h^{5}_{1},h^{1}_{0};f^{1}) = 0
$$ and $$  N(h^{6}_{1},h^{1}_{0};f^{2})=
 N(h^{6}_{1},h^{2}_{0};f^{1})t+ N(h^{6}_{1},h^{1}_{0};f^{1}) = t^2-1.$$
The other flow lines between consecutive critical points remain the same.

 On the second page of the spectral sequence, there are two non-zero differentials, namely 
$d^{2}_{6} $ and $d^{2}_{7}$, which determine the algebraic cancellation
$E^{3}_{6}=E^3_4=0$ and $E^{3}_{7}=E^3_5=0$.  
The corresponding pairs of critical points which will  cancel are  $(h^{5}_{1},h^{7}_{2})$ and $(h^{6}_{1},h^{8}_{2})$, respectively. In order to obtain the circle-valued Morse function $f^{3}$, we cancel  the pair  $(h^{5}_{1},h^{7}_{2})$, see Figures \ref{fig:exemplo-art-can2}, and then  the pair  $(h^{6}_{1},h^{8}_{2})$, see Figure  \ref{fig:exemplo-art-can3}.
In the corresponding negative flow $\varphi^3$ (see Figure \ref{fig:exemplo-art-can3}), the cancellation of the pair of critical points $(h^{6}_{1},h^{8}_{2})$ gives rise to the birth of a repeller periodic orbit $\gamma_{{}_{R}}$ and two flow lines whose $\alpha$-limit set is $\gamma_{{}_{R}}$ and whose $\omega$-limit set is the saddle $h^{4}_{1}$. Note that the birth of the periodic orbit $\gamma_{{}_{R}}$ is due to a primary pivot  $\Delta^{2}_{6,8}$ which is a binomial and the period of $\gamma_{{}_{R}}$ is equal to the difference between the exponents of the binomial.    Moreover, we have the following relation between the Novikov incidence coefficient:
 $N(h^{4}_{1},h^{1}_{0};f^{3})=
 N(h^{4}_{1},h^{1}_{0};f^{2})= t-1.
$

On the third page of the spectral sequence, the differential $d^{3}_{3}$ determines the algebraic cancellation $E^{4}_{4}=E^4_0=0$, which corresponds to the dynamical cancellation of the pair  $(h^{1}_{0},h^{4}_{1})$, creating the circle-valued Morse function $f^4$. In the corresponding negative gradient flow $\varphi^4$, see Figure \ref{fig:exemplo-art-can4}, an attractor periodic orbit $\gamma_{{}_{A}}$ is born due to the primary pivot $\Delta^{3}_{1,4}=t-1$. 
Therefore,  the periodic orbits that will appear after the cancellations are codified in the Novikov matrices.

\vspace{0.7cm}

In order to  prove this relation between dynamical and algebraic cancellations,  we need to use an adaptation of the SSSA which is more conducive to dynamical interpretations.
The dynamical inspiration for this adaptation of the SSSA is as follows.
Note that the changes of basis caused by pivots in row $j-r$ reflect all the changes in connecting orbits caused by the cancellation of $h_k^{j}$ and $h_{k-1}^{j-r}$. However, when the pair of critical points $h_k^{j}$ and $h_{k-1}^{j-r}$ is removed, all the connecting orbits between index $k$ critical points and $h_{k-1}^{j-r}$ and also all the ones between $h_{k}^{j}$ and index $k-1$ critical points are immediately removed and new ones take their place. Hence, in order to interpret dynamically the SSSA, we have to perform the changes of basis that occur therein in a different order to reflect the death and birth of connections. More specifically, if $\Delta^{r}_{j-r,j}$ is a primary pivot marked in step $r$ of the SSSA, 
all changes of basis caused by $\Delta^{r}_{j-r,j}$ will be performed in  step $r+1$. 
This modification in the SSSA is 
called {\it Row Cancellation Algorithm} (RCA).
\vspace{0.2cm}

\begin{bf}Row Cancellation Algorithm\end{bf}
\vspace{-0.3cm}

\begin{description}
\item[\textbf{Initialization Step:}]\mbox{}\\
 $\begin{array}{@{}l}
\left[\begin{tabular}{l}
 $r=0$\\
 $\tilde{\Delta}^r=\Delta$\\
 $\tilde{T}^r=I$ ($m\times m$ identity matrix) \end{tabular}\right.\end{array}$

\item[\textbf{Iterative Step:}] Repeat until all diagonals parallel and to the right of the main diagonal have been swept\\
$\begin{array}{@{}l}
\left[\begin{tabular}{l}
\textbf{Matrix $\Delta$ update}\\
\begin{tabular}{@{\hspace{.5cm}}l}
$r\leftarrow r+1$ \\
$\tilde{\Delta}^r = (\tilde{T}^{r-1})^{-1} \tilde{\Delta}^{r-1}\tilde{T}^{r-1}$\\
\end{tabular}\end{tabular}\right.\\
\end{array}$
\\[5pt]
$\begin{array}{@{}l}
\left[\begin{tabular}{l}
\textbf{Markup}\\
\begin{tabular}{@{\hspace{.5cm}}l}
Sweep entries of $\tilde{\Delta}^r$ in the $r$-th diagonal: \\
\textbf{\textsf{If}} $\tilde{\Delta}^r_{j-r,j}\neq 0$ \textbf{\textsf{and}} $\tilde{\Delta}^r_{\textbf{\large .},j}$ does not contain a primary pivot\\
\rule{.5cm}{0pt}\textbf{\textsf{Then}} permanently mark $\tilde{\Delta}^r_{j-r,j}$ as a primary pivot\\
\end{tabular}\end{tabular}\right.\\
\end{array}$
\\[5pt]
$\begin{array}{@{}l}
\left[\begin{tabular}{l}
\textbf{Matrix $\tilde{T}^r$ construction}\\
\begin{tabular}{@{\hspace{.5cm}}l}
$\tilde{T}^r \leftarrow I$\\
\textbf{\textsf{For each}} primary pivot $\tilde{\Delta}^r_{p-r,p}$ such that $j<m$,
change the $p$-th row of $\tilde{T}^r$ as follows\\
\begin{tabular}{@{\hspace{.5cm}}l}
$\tilde{T}^r_{p,\ell}\leftarrow -(1/\tilde{\Delta}^r_{p-r,p})\tilde{\Delta}^r_{p-r,\ell}$, for $\ell=p+1,\dots,m$
\end{tabular}
\end{tabular}\end{tabular}\right.\\
\end{array}$

\item[\textbf{Final Step:}]\mbox{} \\
$\begin{array}{@{}l}\left[
\begin{tabular}{l}
\textbf{Matrix $\Delta$ update}\\
\begin{tabular}{@{\hspace{.5cm}}l}
$r\leftarrow r+1$ \\
$\tilde{\Delta}^r = (T^{r-1})^{-1} \tilde{\Delta}^{r-1} T^{r-1}$\\
\end{tabular}
\end{tabular}\right.
\end{array}$

\end{description}

In \cite{BLMdRS2}, it was proved that, in the case of matrices with entries over $\mathbb{Z}$, the primary pivots on the $r$-th diagonal of $\widetilde{\Delta}^{r}$ marked in the $r$-th step of the RCA coincide in position and value with the ones on the $r$-th diagonal of $\Delta^r$ marked in the $r$-th step of the SSSA.
More details on this algorithm are presented in \cite{BLMdRS2}. The analogous result for the case over $\mathbb{Z}((t))$ is stated in the next lemma  and its proof is completely analogous to the one obtained for coefficients over $\mathbb{Z}$.  

\begin{lem}Let $f:M\rightarrow S^{1}$ be a circle-valued Morse function on a closed $n$-manifold $M$ and  $(\mathcal{N}_{\ast}(f),\Delta)$ be the Novikov chain complex associated to $f$. Let  $\mathcal{F}$ be a finest filtration in $(\mathcal{N}_{\ast}(f),\Delta)$ and let $(E^r, d^r)$ be the associated spectral sequence for this filtered chain complex.
The primary pivots of the matrices $\widetilde{\Delta}^{r}$ of the RCA coincide in position and value with the ones of the matrices $\Delta^r$ of the SSSA.
\end{lem}

\begin{ex} The following set of matrices illustrate the depletion of rows and collumns in the RCA when applied to Example 6.1.

\begin{figure}[!htb]
\begin{minipage}[t]{0.5\linewidth}
 \includegraphics[width=7.5cm]{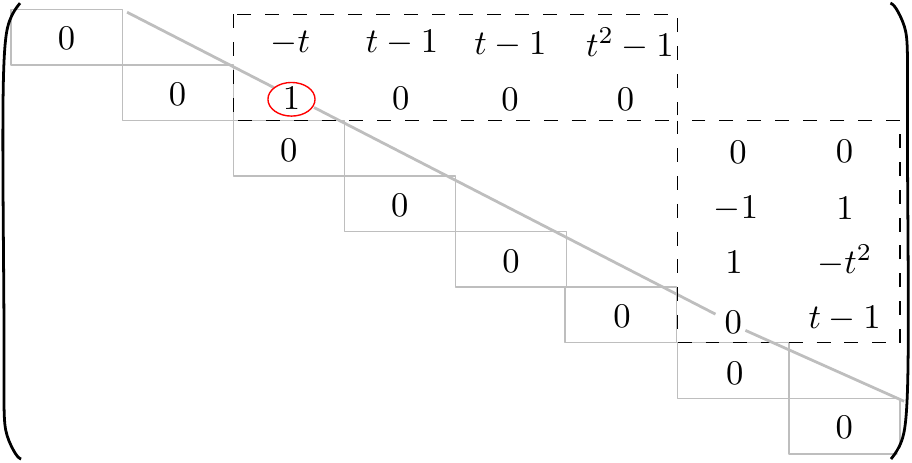}\\ \vspace{-1cm}
\caption{$\widetilde{\Delta}^{1}$ in RCA.}
\label{fig:delta1ex2RCA}
\end{minipage} 
\begin{minipage}[t]{0.5\linewidth}
 \includegraphics[width=7.5cm]{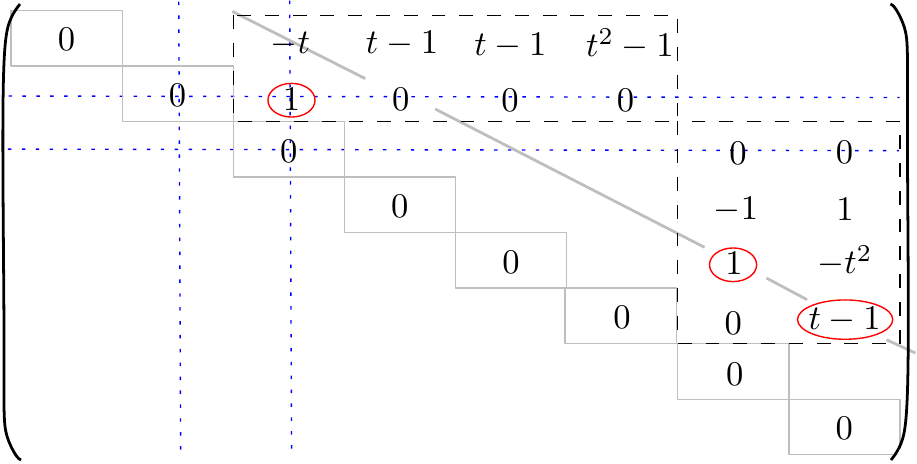}\\ \vspace{-1cm}
\caption{$\widetilde{\Delta}^{2}$ in RCA.}
\label{fig:delta2ex2RCA}
\end{minipage}
\end{figure}
\begin{figure}[!htb]
\begin{minipage}[t]{0.5\linewidth}
\includegraphics[width=7.5cm]{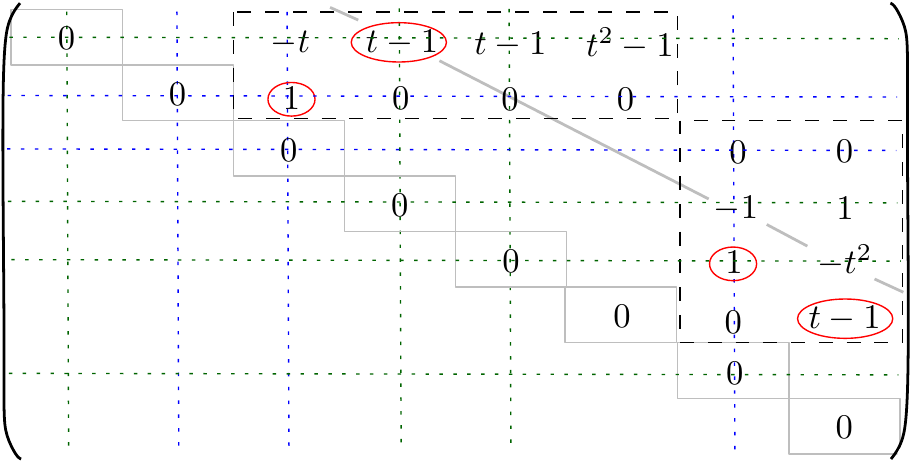}\\ \vspace{-1cm}
\caption{$\widetilde{\Delta}^{3}$ in RCA.}
\label{fig:delta3ex2RCA}
\end{minipage} 
\begin{minipage}[t]{0.5\linewidth}
\includegraphics[width=7.5cm]{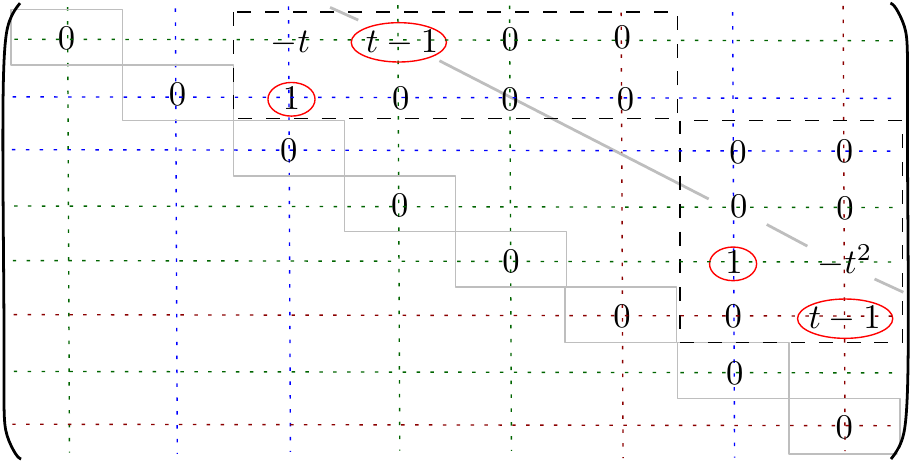}\\ \vspace{-1cm}
\caption{$\widetilde{\Delta}^{4}$ in RCA.}
\label{fig:delta4ex2RCA}
\end{minipage}
\end{figure}
\end{ex}

\vspace{0.2cm}

\noindent \textbf{Proof of Theorem \ref{main}:}
By Theorem \ref{teo:interpretacaodr}, the non zero differentials  of the spectral sequence  are induced by the pivots. Theorem \ref{teo:primarypivots} states that the primary pivots are always invertible polynomials, hence the differentials $d^{r}_p:E^r_{p}\to E^r_{p-r}$ associated to primary pivots are isomorphisms and the ones associated to change of basis pivots always correspond to zero maps.  Consequently, the non-zero differentials  are isomorphisms and this
implies that at the next stage of the spectral sequence they produce
algebraic cancellations, i.e. $E^{r+1}_{p}=E^{r+1}_{p-r}=0$. 

Note that the algebraic cancellation
$E^{r+1}_{p}=E^{r+1}_{p-r}=0$ is associated to the primary pivot
$\Delta^{r}_{p-r+1,p+1}$ on the $r$-th diagonal of
$\Delta^{r}$ in the SSSA. The row $p-r+1$ is associated to
$h_{k-1}^{p-r+1}\in \mathcal{F}_{p-r}C_{k-1}\setminus \mathcal{F}_{p-r-1}C_{k-1}$ and
the column $p+1$ is associated to $h_{k}^{p+1}\in F_pC_k\setminus
F_{p-1}C_k$. Moreover, the primary pivot $\Delta^{r}_{p-r+1,p+1}$ of the SSSA is equal to the primary pivot $\widetilde{\Delta}^{r}_{p-r+1,p+1}$ of the RCA.

We will now prove the theorem by associating dynamical cancellations to the algebraic cancellations of the SSSA via the RCA. Given $\varphi$ the negative gradient flow associated to $f$, we recursively construct a family of negative gradient flows $\{\varphi^r\}$ of circle-valued Morse functions $f^{r}$, $r=0,\ldots , \omega$, where $f^0=f$ and $f^{r+1}$ is obtained from $f^r$ removing the pairs of critical points corresponding to the cancelled modules of the $r$-th page of the spectral sequence, i.e. the pairs associated to the primary pivots on the $r$-th diagonal of $\Delta^r$. In order to do this, we have to prove that whenever a primary pivot $\widetilde{\Delta}^r_{p-r+1,p+1}$ on the $r$-th diagonal of $\widetilde{\Delta}^r$ is marked, it is actually a Novikov incidence coefficient $N(h^{p+1}_{k},h^{p-r+1}_{k-1};f^r)$  between two singularities $h^{p+1}_{k}$ and $h^{p-r+1}_{k-1}$ of consecutive indices of $\varphi^r$.

Without loss of generality, the Novikov differential $\Delta$  associated to $f$ is considered to be the one where all orientations on the unstable manifolds of  critical points of index 2 are the same. In this case, we have the characterization of the Novikov differential $\Delta$ given by Corollary \ref{cor:charac}.

The proof follows by induction in $r$. The base case $r=1$ is trivial, since the differentials $d^0$ are null and hence no pivots are marked in this step. Thus there are no algebraic cancellations, $f^0=f^1$, $\varphi^1=\varphi$ and $\widetilde{\Delta}^{1}=\Delta(M,\varphi^{1})= \Delta(M,\varphi)$, since no columns and rows are removed. Furthermore,
$ N(h^{j}_{k},h^{i}_{k-1};f^{1})=\widetilde{\Delta}^{1}_{ij}$ for all  $h^{j}_{k},h^{i}_{k-1}\in Crit(f^1)$. 

\noindent {\bf Induction hypothesis:}
Suppose that the Morse function $f^{r}$ is well defined and let $\varphi^{r}$ be negative gradient associated to $f^{r}$, where the matrix of the Novikov operator   $\Delta(M,\varphi^{r})$ associated to $f^{r}$ is the submatrix of $\widetilde{\Delta}^{r}$ obtained when we remove all columns and rows  corresponding to all primary pivots marked until diagonal $r-1$ of $\widetilde{\Delta}^{r}$, i.e.,
$ N(h^{j}_{k},h^{i}_{k-1};f^{r})=\widetilde{\Delta}^{r}_{ij}$ for all  $h^{j}_{k},h^{i}_{k-1}\in Crit(f^r)$. 

We now prove the case $r+1$. This proof requires 3 steps. The first step performs dynamical cancellations of the lift of certain critical points of $f^r$ in the covering space. The second step analyses the effect that the dynamical cancellations in the covering space has on the relation between the  Novikov incidence of $f^r$ and $f^{r+1}$.  The third and final step links the algebra of the RCA with the dynamics of $f^{r+1}$, i.e. establishes that $ N(h^{j}_{k},h^{i}_{k-1};f^{r+1})=\widetilde{\Delta}^{r+1}_{ij}$ for all  $h^{j}_{k},h^{i}_{k-1}\in Crit(f^{r+1})$. 

\vspace{0.3cm}

 \begin{bf}Step 1:\end{bf}
 Let $\widetilde{\Delta}_{p-r+1,p+1}^{r}$ be a primary pivot on the $r$-diagonal of $\widetilde{\Delta}^{r}$. By the induction hypothesis,  this primary pivot represents the  Novikov incidence coefficient between two singularities of the flow $\varphi^{r}$, namely  $h^{p+1}_{k}$ and $h^{p-r+1}_{k-1}$. By Theorem \ref{teo:primarypivots}, the primary pivot is either $\pm t^{\ell}$ or $\mp t^{\tilde{\ell}}\pm t^{\ell}$, where $ \ell, \tilde{\ell}\in \mathbb{Z}$, and we assume $\ell < \tilde{\ell}$.
Now, consider the pair of critical points  $t^{\lambda}h^{p+1}_{k}$ and $t^{\lambda + \ell}h^{p-r+1}_{k-1}$ in the covering space $\overline{M}$. Denote by $F^{r}$ the Morse function which makes the following diagram commutative 
\begin{equation}\label{diagramanov}
   \xymatrix@R+2em@C+2em{
  \overline{M} \ar[r]^-{F^{r}} \ar[d]_E & \mathbb{R} \ar[d]^{Exp} \\
  M \ar[r]^-{f^{r}}  & \mathbb{S}^1
  }
 \end{equation}
Since $F^{r}$ is a Morse-Smale function, it  can be perturbed if necessary so that there are no critical points in  $(F^{r})^{-1}[ \ (c_{\lambda+\ell,p-r+1} \ , \  c_{\lambda,p+1}) \ ]$. For a nice proof of this fact see Lemma 3.2 in \cite{Sa2}.

Using the Smale's Cancellation Theorem, we can equivariantly cancel the pairs of critical points  $t^{\lambda}h^{p+1}_{k}$ and $t^{\lambda + \ell}h^{p-r+1}_{k-1}$, for all $\lambda \in \mathbb{Z}$, obtaining a new Morse-Smale function $F^{r+1}:\overline{M}\rightarrow \mathbb{R}$ such that $Crit(F^{r+1})= Crit(F^{r})\setminus \{t^{\lambda}h^{p+1}_{k} \ , \ t^{\lambda + \ell}h^{p-r+1}_{k-1} \mid \lambda \in \mathbb{Z} \}$. Moreover, we have the following correspondence between the intersection numbers with respect to the function $F^{r}$ and $F^{r+1}$, denoted by  $n(\cdot,\cdot;F^{r})$ and $n(\cdot,\cdot;F^{r+1})$, respectively:

\begin{enumerate}
\item[1.] For $k=1$:
\begin{enumerate}
\item[1.1.] If the saddle $t^{\lambda}h^{p+1}_{1}$ connects with exactly two sinks, namely $t^{\lambda+\ell}h^{p-r+1}_{0}$ and  $t^{\lambda+\tilde{\ell}}h^{i_{0}}_{0}$, then given a saddle $t^{\delta}h^{j}_{1} \in Crit(F^{r+1})$, we have that 
\begin{equation}\label{eq:intersection1}
n(t^{\delta}h^{j}_{1},t^{\lambda+\tilde{\ell}}h^{i_{0}}_{0};F^{r+1}) = n(t^{\delta}h^{j}_{1},t^{\lambda +\ell}h^{p-r+1}_{0};F^{r}) + n(t^{\delta}h^{j}_{1},t^{\lambda+\tilde{\ell}}h^{i_{0}}_{0};F^{r}),
\end{equation}
and
\begin{equation}\label{eq:intersection11}
n(t^{\delta}h^{j}_{1},t^{\lambda+\tilde{\ell}}h^{i}_{0};F^{r+1}) =  n(t^{\delta}h^{j}_{1},t^{\lambda+\tilde{\ell}}h^{i}_{0};F^{r}),\,\, i\neq i_{0}
\end{equation} 
 In the particular case that $i_{0}=p-r+1$, since the cancellation is done equivariantly, an eventual flow line in $\varphi^{r}$ from $t^{\delta}h^{j}_{1}$ to $t^{\lambda+\tilde{\ell}}h^{i_{0}}_{0}$ gives place to a flow line starting at  $t^{\delta}h^{j}_{1}$ with  empty $\omega$-limit set.
\item[1.2.]  If the saddle $t^{\lambda}h^{p+1}_{1}$ connects with exactly one sink, namely $t^{\lambda+\ell}h^{p-r+1}_{0}$, then the  intersection numbers between critical points of $F^{r+1}$ remain  the same.
\end{enumerate}

 \item[2.] For $k=2$:
\begin{enumerate}
\item[2.1.] If the saddle $t^{\lambda+\ell}h^{p-r+1}_{1}$ connects with exactly two sources, namely $t^{\lambda}h^{p+1}_{2}$ and  $t^{\lambda+\tilde{\ell}}h^{j_{0}}_{2}$, then given a saddle $t^{\delta}h^{i}_{1} \in Crit(F^{r+1})$, we have that 
\begin{equation}\label{eq:intersection2}
n(t^{\lambda+\tilde{\ell}}h^{j_{0}}_{2},t^{\delta}h^{i}_{1};F^{r+1}) = n(t^{\lambda}h^{p+1}_{2},t^{\delta}h^{i}_{1};F^{r}) + n(t^{\lambda+\tilde{\ell}}h^{j_{0}}_{2},t^{\delta}h^{i}_{1};F^{r}),
\end{equation}
and 
\begin{equation}\label{eq:intersection21}
n(t^{\lambda+\tilde{\ell}}h^{j}_{2},t^{\delta}h^{i}_{1};F^{r+1}) =  n(t^{\lambda+\tilde{\ell}}h^{j}_{2},t^{\delta}h^{i}_{1};F^{r}),\,\, j\neq j_0
\end{equation}
In the particular case that $j_{0}=p+1$, since the cancellation is done equivariantly, an eventual flow line in $\varphi^{r}$ from  $t^{\lambda+\tilde{\ell}}h^{j_{0}}_{2}$ to $t^{\delta}h^{i}_{1}$  gives place to a flow line ending at $t^{\delta}h^{i}_{1}$ with empty $\alpha$-limit set.
\item[2.2.]  If the saddle $t^{\lambda+\ell}h^{p-r+1}_{1}$ connects with exactly one source, namely $t^{\lambda}h^{p+1}_{2}$, 
 then the  intersection numbers between critical points of $F^{r+1}$ remain  the same.

\end{enumerate} 

 \end{enumerate}
 
Formulas  (\ref{eq:intersection1}),  (\ref{eq:intersection11}),  (\ref{eq:intersection2}) and (\ref{eq:intersection21}) are a consequence of the characterization of the Morse boundary operator, when we consider the orientations on the unstable manifolds of index 2 critical points of $F^{r}$ being the same. This follows easily since the characteristic signs of the two flow lines of the unstable (stable) manifold of a saddle are opposite 
as proved in \cite{BLMdRS}.

\vspace{0.3cm}

\begin{bf}Step 2\end{bf}:
Let $f^{r+1}: M\rightarrow S^{1}$ be the smooth map 
such that $f^{r+1}(x)  = Exp \circ F^{r+1} (y)$, where $y \in  E^{-1}(x)$ for all $x\in M$. Since the cancellation was done equivariantly then $f^{r+1}$ is well defined and it is a circle-valued Morse function  such that $Crit(f^{r+1}) = Crit(f^{r})\setminus \{h^{p+1}_{k},h^{p-r+1}_{k-1}\}$ which coincides up to isotopy with $f^{r}$ outside a neighborhood of the flow line joining these critical points.

 As consequence of the formulas in (\ref{eq:intersection1}),  (\ref{eq:intersection11}),  (\ref{eq:intersection2}) and (\ref{eq:intersection21}), we have the following relation between  the Novikov incidence coefficients when considering the functions $f^{r}$ and $f^{r+1}$:

 \begin{enumerate}
\item[3.] If $k=1$:

\begin{enumerate}

\item[3.1.] Suppose that $\widetilde{\Delta}^{r}_{p-r+1,p+1} = \pm t^{\ell} \mp t^{\tilde{\ell}}$ where $\ell <\tilde{\ell}$, then   there is a double connection between $h^{p+1}_{1}$ and $h^{p-r+1}_{0}$ in the flow $\varphi^{r}$. This corresponds to case 1.1 in the covering space  when $i_{0}=p-r+1$. After the equivariant cancellation, new connections are not created between consecutive critical points of $F^{r+1}$, hence $N(h^{j}_{1},h^{i}_{0};f^{r+1}) = N(h^{j}_{1},h^{i}_{0};f^{r})$, for all $h^{j}_{1},h^{i}_{0} \in Crit(f^{r+1})$.

\item[3.2.] Suppose that $\widetilde{\Delta}^{r}_{p-r+1,p+1} = \pm t^{\ell}$.
Given $ h^{i}_{0} \in Crit(f^{r+1})$ with $i<p-r+1$, the entry $\widetilde{\Delta}^{r}_{i,p+1}$ belongs to the submatrix $\Delta(M,\varphi^{r})$, consequently we can use Corollary \ref{cor:charac},  which implies that we have two possible cases: 
\begin{itemize}
\item  
  $\widetilde{\Delta}^{r}_{i,p+1} = 0$ for all    $i<p-r+1$ such that  $ h^{i}_{0} \in Crit(f^{r+1})$. This corresponds to  case  1.2 above. After the equivariant cancellation, there are no changes in the connections between consecutive critical points of $F^{r+1}$, hence $N(h^{j}_{1},h^{i}_{0};f^{r+1}) = N(h^{j}_{1},h^{i}_{0};f^{r})$, for $h^{j}_{1}\in Crit(f^{r+1}) $.
\item 
there exists  $i_{0}<p-r+1$ such that   $ h^{i_{0}}_{0} \in Crit(f^{r+1})$,  $\widetilde{\Delta}^{r}_{i_{0},p+1} = \mp t^{\tilde{\ell}}$ and $\widetilde{\Delta}^{r}_{i,p+1} = 0$ for all    $i<p-r+1$ with $i\neq i_ {0}$  and  $ h^{i}_{0} \in Crit(f^{r+1})$.  In this case, $h^{p+1}_{1}$ connects with exactly two sinks, namely $h^{p-r+1}_{0}$ and  $h^{i_{0}}_{0}$ in the flow $\varphi^{r}$. Then, by case 1.1 above, we have that for all  $h^{j}_{1}\in Crit(f^{r+1}) $, 
\begin{eqnarray}\label{eq:intersection1N}
N(h^{j}_{1},h^{i_{0}}_{0};f^{r+1}) & = &  \ \ \ \sum_{\lambda\in \mathbb{Z}}n(h^{j},t^{\lambda}h^{i_{0}};F^{r+1})t^{\lambda} \nonumber \\ \nonumber
& = & \ \ \ \sum_{\lambda\in \mathbb{Z}}n(t^{\tilde{\ell}}h^{j},t^{\lambda+\tilde{\ell}}h^{i_{0}};F^{r+1})t^{\lambda}\\ \nonumber
&=& \ \ \ \sum_{\lambda\in \mathbb{Z}}\Big[n(t^{\tilde{\ell}}h^{j}_{1},t^{\lambda +\ell}h^{p-r+1}_{0};F^{r}) + n(t^{\tilde{\ell}}h^{j}_{1},t^{\lambda+\tilde{\ell}}h^{i_{0}}_{0};F^{r})\Big]t^{\lambda},\\ \nonumber
&=& \ \ \ \sum_{\lambda\in \mathbb{Z}}n(h^{j}_{1},t^{\lambda +\ell - \tilde{\ell}}h^{p-r+1}_{0};F^{r}) t^{\lambda} + 
\sum_{\lambda\in \mathbb{Z}} n(h^{j}_{1},t^{\lambda}h^{i_{0}}_{0};F^{r})t^{\lambda}\\ \nonumber
&=&  \sum_{\lambda-(\ell-\tilde{\ell})\in \mathbb{Z}}n(h^{j}_{1},t^{\lambda }h^{p-r+1}_{0};F^{r}) t^{\lambda - (\ell-\tilde{\ell})} + 
N(h^{j}_{1},h^{i_{0}}_{0};f^{r})\\ 
& = & N(h^{j}_{1},h^{p-r+1}_{0};f^{r})t^{\tilde{\ell}-\ell}+ N(h^{j}_{1},h^{i_{0}}_{0};f^{r}),
\end{eqnarray} where the third equality follows from (\ref{eq:intersection1}).   Moreover,  $N(h^{j}_{1},h^{i}_{0};f^{r+1}) = N(h^{j}_{1},h^{i}_{0};f^{r})$ for $i\neq i_{0}$.
\end{itemize}
\end{enumerate}

 \item[4.] If $k=2$:
 \begin{enumerate}
 
 \item[4.1.] Suppose that $\widetilde{\Delta}^{r}_{p-r+1,p+1} = \pm t^{\ell} \mp t^{\tilde{\ell}}$ where $\ell<\tilde{\ell}$, then   there is a double connection between $h^{p+1}_{2}$ and $h^{p-r+1}_{1}$ in the flow $\varphi^{r}$. This corresponds to the case 2.1 in the covering space when $j_{0}=p+1$. After the equivariant cancellation new connections are not created between consecutive critical points of $F^{r+1}$, hence
  $N(h^{j}_{2},h^{i}_{1};f^{r+1}) = N(h^{j}_{2},h^{i}_{1};f^{r})$, for all $h^{j}_{2},h^{i}_{1} \in Crit(f^{r+1})$.

\item[4.2.] Suppose that $\widetilde{\Delta}^{r}_{p-r+1,p+1} = \pm t^{\ell}$. 
Given $ h^{j}_{2} \in Crit(f^{r+1})$ with $j>p+1$, the entry $\widetilde{\Delta}^{r}_{p-r+1,j}$ belongs to the submatrix $\Delta(M,\varphi^{r})$, consequently we can use Corollary \ref{cor:charac},  which implies that we have two possible cases: 
\begin{itemize}
\item  
  $\widetilde{\Delta}^{r}_{p-r+1,j} = 0$ for all    $j>p+1$ such that  $ h^{j}_{2} \in Crit(f^{r+1})$. This corresponds to the case  2.2 above. After the equivariant cancellation, there are no changes in the connections between consecutive critical points of $F^{r+1}$, hence $N(h^{j}_{2},h^{i}_{1};f^{r+1}) = N(h^{j}_{2},h^{i}_{1};f^{r})$, for $h^{i}_{1}\in Crit(f^{r+1}) $.
\item 
there exists  $j_{0}>p+1$ such that   $ h^{j_{0}}_{2} \in Crit(f^{r+1})$,  $\widetilde{\Delta}^{r}_{p-r+1,j_{0}} = \mp t^{\tilde{\ell}}$ and $\widetilde{\Delta}^{r}_{p-r+1,j} = 0$ for all    $j>p+1$ with $j\neq j_ {0}$  and  $ h^{j}_{2} \in Crit(f^{r+1})$.  In this case, $h^{p-r+1}_{1}$ connects with exactly two sources, namely $h^{p+1}_{2}$ and  $h^{j_{0}}_{2}$ in the flow $\varphi^{r}$. Then, by case 2.1 above, we have that for all  $h^{i}_{1}\in Crit(f^{r+1}) $, 
\begin{equation}\label{eq:intersection2N}
N(h^{j_{0}}_{2},h^{i}_{1};f^{r+1}) = N(h^{p+1}_{2},h^{i}_{1},f^{r})t^{\tilde{\ell}-\ell} + N(h^{j_{0}}_{2},h^{i}_{1};f^{r}).
\end{equation}
  Moreover,  $N(h^{j}_{2},h^{i}_{1};f^{r+1}) = N(h^{j}_{2},h^{i}_{1};f^{r})$ for $j\neq j_{0}$.

\end{itemize}

\end{enumerate}
 
 \end{enumerate}
 
 \vspace{0.3cm}
 
 \begin{bf}Step 3:\end{bf}
 Now, we link the dynamics with the algebra of the RCA. More specifically, we will show that the matrix of the Novikov operator   $\Delta(M,\varphi^{r+1})$ associated to $f^{r+1}$ is the submatrix of $\widetilde{\Delta}^{r+1}$ obtained when we
remove all columns and rows  corresponding to all primary pivots marked up to the $r$-th diagonal of $\widetilde{\Delta}^{r}$.

Let $\widetilde{\Delta}^{r}_{p-r+1,p+1}$ be a primary pivot and consider an entry $\widetilde{\Delta}^{r+1}_{i,j}$ with $i<p-r+1$ and $j> p+1$. If there is a primary pivot in column/row $i$ or column/row $j$, then the $\widetilde{\Delta}^{r+1}_{i,j}$ is not in $\Delta(M,\varphi^{r+1})$.
 Otherwise,  the entry $\widetilde{\Delta}^{r+1}_{i,j}$
 is given by 
\begin{equation}\label{eq:formulasweeping}
\widetilde{\Delta}^{r+1}_{i,j}  =
     \widetilde{\Delta}^{r}_{i,j}  - \dfrac{\widetilde{\Delta}^{r}_{p-r+1,j} }{\widetilde{\Delta}^{r}_{p-r+1,p+1} } \widetilde{\Delta}^{r}_{i,p+1}.
\end{equation}  
 
From now on consider $\widetilde{\Delta}^{r}_{i,j}$ such that $h^{j}_{k}, h^{i}_{k-1} \in Crit(f^{r+1})$. In this case, there are no primary pivots in column/row $i$ nor in column/row $j$, since $\widetilde{\Delta}^{r}_{i,j}$ is an entry of the Novikov matrix
  $\Delta(M,\varphi^{r})$ of $f^{r}$.

 \begin{enumerate}
\item[5.] If $k=1$:

\begin{enumerate}

\item[5.1.] Suppose that $\widetilde{\Delta}^{r}_{p-r+1,p+1} = \pm t^{\ell} \mp t^{\tilde{\ell}}$. By Corollary \ref{cor:charac}, 
$\widetilde{\Delta}^{r}_{i,p+1} = 0 $ for  $i< p-r+1 $  and  $ h^{i}_{0} \in Crit(f^{r+1})$.  Hence, 
  $\widetilde{\Delta}^{r+1}_{i,j}  =
     \widetilde{\Delta}^{r}_{i,j} $. By the induction hypothesis and   by  item 3.1 above
$\widetilde{\Delta}^{r+1}_{i,j} = N(h^{j}_{1},h^{i}_{0};f^{r+1})$, for all $j$ such that $ h^{j}_{1} \in Crit(f^{r+1})$.

\item[5.2.] Suppose that $\widetilde{\Delta}^{r}_{p-r+1,p+1} = \pm t^{\ell}$. We have two cases to consider.
\begin{itemize}
\item $\widetilde{\Delta}^{r}_{i,p+1} = 0$ for all    $i<p-r+1$ such that  $ h^{i}_{0} \in Crit(f^{r+1})$.  Then  $\widetilde{\Delta}^{r+1}_{i,j}  =
     \widetilde{\Delta}^{r}_{i,j} $. By the induction hypothesis and  by  case 3.2 above 
$\widetilde{\Delta}^{r+1}_{i,j} = N(h^{j}_{1},h^{i}_{0};f^{r+1})$.

\item there exists  $i_{0}<p-r+1$ such that   $ h^{i_{0}}_{0} \in Crit(f^{r+1})$,  $\widetilde{\Delta}^{r}_{i_{0},p+1} = \mp t^{\tilde{\ell}}$ and $\widetilde{\Delta}^{r}_{i,p+1} = 0$ for all    $i<p-r+1$ with $i\neq i_ {0}$  and  $ h^{i}_{0} \in Crit(f^{r+1})$.  In this case, (\ref{eq:formulasweeping}) is equivalent to $\widetilde{\Delta}^{r+1}_{i,j}  =
     \widetilde{\Delta}^{r}_{i,j}  + \widetilde{\Delta}^{r}_{p-r+1,j}  \ t^{\tilde{\ell}-\ell} $. Since by the induction hypothesis $\widetilde{\Delta}^{r}_{i,j}  = N(h^{j}_{1},h^{i}_{0};f^{r}) $   and 
     $\widetilde{\Delta}^{r}_{p-r+1,j}= N(h^{j}_{1},h^{p}_{0};f^{r})$, then by item 3.2  (equation (\ref{eq:intersection1N})) we have $\widetilde{\Delta}^{r+1}_{i,j} = N(h^{j}_{1},h^{i}_{0};f^{r+1})$.
\end{itemize}

\end{enumerate}

 \item[6.] If $k=2$:
 \begin{enumerate}

\item[6.1.] Suppose that $\widetilde{\Delta}^{r}_{p-r+1,p+1} = \pm t^{\ell} \mp t^{\tilde{\ell}}$. 
Given an entry $\widetilde{\Delta}^{r}_{i,j}$ with $j>p+1$ and $h^{j}_{2}, h^{i}_{1} \in Crit(f^{r+1})$, then $\widetilde{\Delta}^{r}_{p-r+1,j} = 0 $, since the submatrix $\Delta(M,\varphi^{r} )$ of $\widetilde{\Delta}^{r}$ satisfies Corollary \ref{cor:charac}.
 Hence, 
$\widetilde{\Delta}^{r+1}_{i,j}  =
     \widetilde{\Delta}^{r}_{i,j} $ and, by the induction hypothesis,  
$\widetilde{\Delta}^{r+1}_{i,j} = N(h^{j}_{2},h^{i}_{1};f^{r+1})$.

\item[6.2.] Suppose that $\widetilde{\Delta}^{r}_{p-r+1,p+1} =\pm t^{\ell}$. 
Given an entry $\widetilde{\Delta}^{r}_{i,j}$ with $j>p+1$ and $h^{j}_{2}, h^{i}_{1} \in Crit(f^{r+1})$, then $\widetilde{\Delta}^{r}_{p-r+1,j} = 0 $ or $\widetilde{\Delta}^{r}_{p-r+1,j} = \mp t^{\tilde{\ell}} $, since the submatrix $\Delta(M,\varphi^{r} )$ of $\widetilde{\Delta}^{r}$ satisfies Corollary \ref{cor:charac}.
In the case that $\widetilde{\Delta}^{r}_{p-r+1,j} =  0 $, then $\widetilde{\Delta}^{r+1}_{i,j}  =
     \widetilde{\Delta}^{r}_{i,j} $ and  
$\widetilde{\Delta}^{r+1}_{i,j} = N(h^{j}_{2},h^{i}_{1};f^{r+1})$.
On the other hand, if  $\widetilde{\Delta}^{r}_{p-r+1,j}  = \mp t^{\tilde{\ell}}$,  (\ref{eq:formulasweeping}) is equivalent to $\widetilde{\Delta}^{r+1}_{i,j}   = \widetilde{\Delta}^{r}_{i,j} + \widetilde{\Delta}^{r}_{i,p+1} \ t^{ \tilde{\ell}-\ell}$. 
By the induction hypothesis $\widetilde{\Delta}^{r}_{i,j}  = N(h^{j}_{2},h^{i}_{1};f^{r}) $   and 
     $\widetilde{\Delta}^{r}_{i,p+1}= N(h^{p+1}_{2},h^{i}_{1};f^{r})$. Then by item 4.2 above (equation (\ref{eq:intersection2N})) we have $\widetilde{\Delta}^{r+1}_{i,j} = N(h^{j}_{2},h^{i}_{1};f^{r+1})$.

\end{enumerate}
 \end{enumerate}

 If there are more primary pivots $\{\widetilde{\Delta}^{r}_{p_{\ell}-r+1,p_{\ell}+1}\}$ on the $r$-th diagonal of $\widetilde{\Delta}^{r}$, we repeat the above construction successively  by determining  neighborhoods $U_{\ell}$ of the closures of the flows lines joining $h^{p_{\ell}+1}_{k}$ and $h^{p_{\ell}-r+1}_{k-1}$.  Cancelling the pairs of critical points $\{ (h^{p_{\ell}+1}_{k}, h^{p_{\ell}-r+1}_{k-1})\}$ in $U_{\ell}$, we obtain a family of flows $\{\varphi^{r}_{\ell}\}$ such that  $\varphi^{r}_{\ell}$ has exactly two less singularities than $\varphi^{r}_{\ell-1}$. Each one of the cancellations is done as above with the dynamics corresponding to  the algebra of the RCA.
 In this case, the flow  $\varphi^{r+1}$ is obtained from $\varphi^r$ after cancelling all pairs of critical points  corresponding to  primary pivots on the $r$-th diagonal.

\cqd

%
%
%

As a consequence of the proof of Theorem \ref{main},  the appearance of periodic orbits in the family $f^{r}$ due cancellations of critical points can be detected by the spectral sequence.

\begin{teo}[Detecting Periodic Orbits]\label{mainperiodic} 
There exists a family of  circle-valued Morse functions  satisfying  Theorem \ref{main} such that, for each cancelled pair of modules $E^r_p$ and $E^r_{p-r}$
corresponding to a binomial primary pivot  $\Delta^{r}_{p-r+1,p+1}=\mp t^{\ell}\pm t^{\tilde{\ell}}$,  where $\ell <\tilde{\ell}$, the associated dynamical cancellation of critical points $h_{k}^{p+1}$ and $h_{k-1}^{p-r+1}$ gives rise to a periodic orbit that crosses the  regular value $f^{-1}(a)$ $\tilde{\ell}-\ell$ times. This  periodic orbit is either an attractor if $k=0$ or a repeller if $k=1$. 
\end{teo}
\dem  In order to perform the cancellation of $h_{k}^{p+1}$ and $h_{k-1}^{p-r+1}$, we consider the critical points  $t^{\lambda}h^{p+1}_k$ and $t^{\lambda +\ell}h^{p-r+1}_{k-1}$ in the covering space $\overline{M}$ and cancel them equivariantly for all $\lambda\in\mathbb{Z}$  using the First Cancellation Theorem in  \cite{M}. 
Note that, for a  fixed $\lambda$, when the pair $t^{\lambda}h^{p+1}_k$ and $t^{\lambda +\ell}h^{p-r+1}_{k-1}$ of critical points of $F^{r}$ is cancelled, the function $F^{r+1}$ conicides with $F^{r}$ (up to isotopy) outside a neighborhood $U$ of  the  connecting orbit between them and a  flow line $\gamma$ from $t^{\lambda-\ell}h^{p+1}_k$ to $t^{\lambda +\tilde{\ell}}h^{p-r+1}_{k-1}$ arises in the negative gradient flow of $F^{r+1}$. One can choose the diffeomorfism  given by the  Assertion 4 in Chapter 5 of \cite{M} such  that $\gamma$  coincides outside of $U$ with the flow lines in the negative gradient flow of  $F^{r}$ going from 
$t^{\lambda-\ell}h^{p+1}_k$ to $t^{\lambda +{\ell}}h^{p-r+1}_{k-1}$ and from
 $t^{\lambda}h^{p+1}_k$ to $t^{\lambda +\tilde{\ell}}h^{p-r+1}_{k-1}$.  Cancelling equivariantly in $\overline{M}$, one obtains a periodic orbit of the negative gradient flow of  $f^{r+1}$ on $M$.
 \cqd

\begin{ex}
We could have considered an alternative finest filtration in Example \ref{exemplocap6}, such as
 $\mathcal{G}=\{G_{p}\mathcal{N}\}$, where
$ G_{0}\mathcal{N} = \mathbb{Z}((t))\big[h^1_0\big] $,
$ G_{1}\mathcal{N} = \mathbb{Z}((t))\big[h^1_0, h^2_0\big] $, 
$ G_{2}\mathcal{N} = \mathbb{Z}((t))\big[h^1_0, h^2_0,h^3_1\big] $, 
$ G_{3}\mathcal{N} = \mathbb{Z}((t))\big[h^1_0, h^2_0,h^3_1,h^4_1\big] $, 
$ G_{4}\mathcal{N} = \mathbb{Z}((t))\big[h^1_0, h^2_0,h^3_1,h^4_1,h^6_1\big] $, 
$ G_{5}\mathcal{N} = \mathbb{Z}((t))\big[h^1_0, h^2_0,h^3_1,h^4_1,h^6_1,h^5_1\big] $, 
$ G_{6}\mathcal{N} = \mathbb{Z}((t))\big[h^1_0, h^2_0,h^3_1,h^4_1,h^6_1,h^5_1, h^8_2\big] $ and
$ G_{7}\mathcal{N} = \mathbb{Z}((t))\big[h^1_0, h^2_0,h^3_1,h^4_1,h^6_1,h^5_1, h^8_2,h^{7}_2\big] $.
 In this case, the dynamical cancellations of pairs of critical points associated to algebraic cancellations of the modules of the spectral sequence through the SSSA will be the following pairs of critical points: $(h_{0}^{2},h_{1}^{3})$, $(h_{1}^{5},h_{2}^{8})$, $(h_{0}^{1},h_{1}^{4})$, $(h_{1}^{6},h_{2}^{7})$. For instance, in the first matrix  $\Delta^{1}$ produced by the SSSA, see Figure \ref{fig:delta1-exefinal}, the entry $\Delta^{1}_{6,7}=-t^2$  is a primary pivot which induces the differential $d^{1}_{6}:E^{1}_6\rightarrow E^{1}_{5}$ of the spectral sequence. Note that by the choice of the finest filtration $\mathcal{G}$,  $E^{1}_6$ is generated by the critical point $h^{8}_{2}$ and $E^{1}_5$ is generated by the critical point $h^{5}_{1}$. Therefore, the algebraic cancellation  of the modules  $E^{1}_6$ and $ E^{1}_{5}$  is associated to the dynamical cancellation of the pair $(h_{1}^{5},h_{2}^{8})$.
Also note that, despite the fact that  $h_{2}^{8}$ and $h_{1}^{5}$ are close with respect to the filtration $\mathcal{G}$, i.e. the gap between them is $1$, the flow line from $h_{2}^{8}$ to $h_{1}^{5}$ is a "long"  orbit both in $M$ and $\overline{M}$ represented by $-t^2$ on the Novikov differential. This orbit in the torus $M$ wraps twice around it. Hence,  the corresponding cancelled pairs $(t^{\lambda +2}h_{1}^{5}, t^{\lambda} h_{2}^{8})$ in the covering space  are also not close compared with the gap in the filtration $\mathcal{G}_{\lambda,p}$, since the gap between them is $19$. The algebraic cancellation corresponding to $-t^2$ also has its dynamical counterpart.

The choice of a different filtration for the Novikov chain complex from the one in Example 6.2, implies in different cancellations of pairs of critical points. The family of negative gradient flows $\varphi^r$ determined by Theorem 6.1 for this example exhibit new orbits that wrap around the torus many times possibly with different orientations in some of the turns. This dynamical behavior is determined by the primary pivots in the Novikov matrices  $\Delta^r$. Once all of the cancellations are performed, the flow associated to the last Novikov matrix is minimal and in this case we are left with no critical points and one attracting and one repelling periodic orbit.

\vspace{0.3cm}

\begin{figure}[!htb]
\centering
\includegraphics[width=7cm]{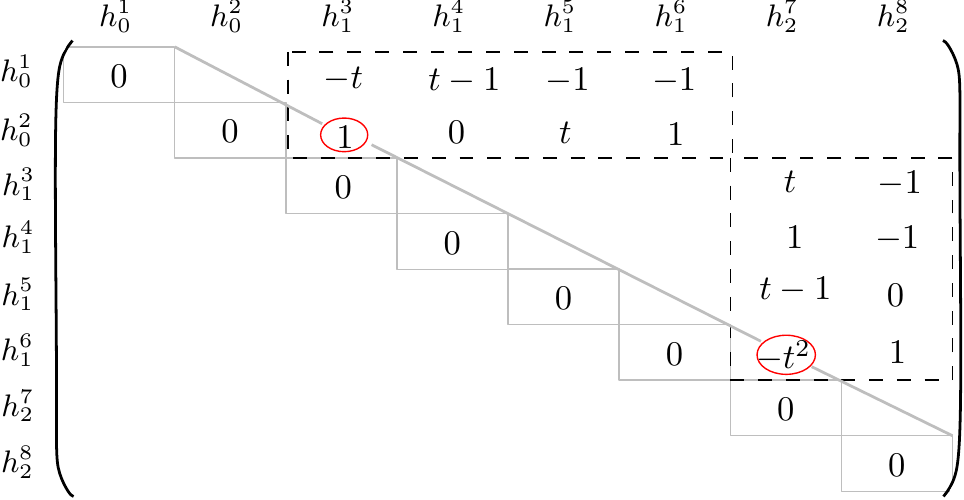}
\vspace{-0.3cm}
\caption{$\Delta^{1}$ for an alternative filtration $\mathcal{G}$.}
\label{fig:delta1-exefinal}
\end{figure}

\end{ex}

\vspace{1cm}

\textsc{Acknowledgments}. \hspace{.1cm} We wish to thank Margarida Mello for her dependable collaboration throughout the years and her invaluable contribution in the implementation of the algorithms presented herein, which facilitates the calculation of the Novikov matrices and thus provides a wealth of examples. Moreover, these implementations are an exceptional platform  for experimental investigation.

\vspace{1cm}

{\sc {\footnotesize
\noindent  D.V.S. Lima - IMECC, Universidade Estadual de Campinas, Campinas, SP, Brazil.
 e-mail: dahisylima@ime.unicamp.br

\noindent  O. M. Neto - ICMC, Universidade de So Paulo, So Carlos, SP, Brazil.
 e-mail: ozimneto@icmc.usp.br

\noindent  K.A. de Rezende - IMECC, Universidade Estadual de Campinas, Campinas, SP, Brazil. e-mail: ketty@ime.unicamp.br

 \noindent   M. R. da
Silveira  - CMCC, Universidade Federal do ABC, Santo Andr{\'e}, SP, Brazil.  e-mail: mariana.silveira@ufabc.edu.br 
}}

\end{document}